\documentclass[a4paper,11pt
]{amsart}
\usepackage[foot]{amsaddr}
\usepackage[margin=2.7cm]{geometry}
\newtheorem{theorem}{Theorem}[section]

\theoremstyle{remark}
\newtheorem{remark}[theorem]{Remark}
\usepackage{systeme}
\usepackage{graphicx} 
\usepackage[]{units}
\usepackage{color}
\usepackage{mathrsfs}
\usepackage{bm}
\usepackage[ansinew]{inputenc}
 \usepackage[breaklinks]{hyperref}
\usepackage{float}
\usepackage{amsmath}
\usepackage{isomath} 
\usepackage{amsthm}
\usepackage{amssymb}

\usepackage{hyperref}

\newcommand{\red}{\color{black}}
\newcommand{\blue}{\color{black}}

\newcommand{\Xgreen}{\color{black}}

\newcommand{\orange}[1]{{\color{black}#1}}

\newcommand{\magenta}[1]{\color{black}#1}
\newcommand{\bblue}[1]{\color{black}#1}

\usepackage{bm}
\usepackage{graphicx,color, amsfonts,amsmath}

\usepackage{mwe}
\usepackage{subcaption,float}

\newcommand{\ldb}{\mathopen{\lbrack\!\lbrack}}
\newcommand{\rdb}{\mathclose{\rbrack\!\rbrack}}

\usepackage{amssymb}

\DeclareMathAlphabet\mathbfcal{OMS}{cmsy}{b}{n}

\begin{document}

\title[{{{ROM for}} geometrically parametrized Cahn-Hilliard system based on CutFEM
}]{ {\magenta{A reduced order model}} for a stable embedded boundary parametrized Cahn-Hilliard phase-field system based on cut finite {{elements}}
}

\author{Efthymios N. Karatzas\textsuperscript{1,2}}
\address{\textsuperscript{1}SISSA, International School for Advanced Studies, Mathematics Area, mathLab, Trieste, Italy.}
\address{\textsuperscript{2} Department of Mathematics, National Technical University of Athens, Greece and FORTH Institute of Applied and Computational Mathematics, Heraclion, Crete, Greece.}
\email{karmakis@math.ntua.gr
}




\author{Gianluigi Rozza\textsuperscript{1}}
\email{grozza@sissa.it}

\subjclass[2010]{78M34, 97N40, 35Q35}

\keywords{Cut Finite Element Method,
          Cahn Hilliard,
          Reduced Order Model,
          POD,
          Stabilization}
  

\dedicatory{}





\begin{abstract}
 In the present work, we {\magenta{investigate 
 a}} cut finite element method {\magenta{for the}} parameterized {\magenta{system of second-order equations stemming from}} the splitting approach of a fourth order nonlinear geometrical PDE, namely the Cahn-Hilliard system. 
%
%
We manage to tackle the instability issues of such methods whenever strong nonlinearities appear and to utilize their flexibility of the fixed background geometry --and mesh-- characteristic, through which, one can avoid e.g. in parametrized geometries the remeshing on the full order level, as well as, transformations to reference geometries on the reduced level. As a final goal, we manage to find an efficient {{global}}, concerning the geometrical manifold, and {{independent}} of geometrical changes, reduced order basis. The POD-Galerkin approach exhibits its strength even with pseudo-random discontinuous initial data verified by numerical experiments.
\end{abstract}
\maketitle
%
\section{Introduction and {\magenta{m}}otivation}\label{sec:intro}
The Cahn-Hilliard model (CH), named after John Cahn and John Hilliard who suggested the system in 1958
, describes a prototype of the process of phase separation, by which the two components of a binary fluid or material impulsively separate and form domains, pure in each component.
These phase-field systems, are very interesting in the scientific community due to their great conservation properties {\magenta{in the sense that the phase separation process conserves the total concentration}}. {\magenta{They can be used to simulate}} many industrial systems, as the two-phase fluid flows for capturing the interface location between two immiscible fluids, \cite{Alpak2016,Anderson98,Israelachvili11}, spinodal decomposition in binary alloys --a process in which a mixture of two fluids or materials {\magenta{decomposes into the pure materials}}-- \cite{GURTIN96,Junseok16}, and the phase diagram for microphase multiscale separation for diblock copolymer-linear chain molecule consisting of two subchains joined covalently to each other, \cite{Choksi09,Jeong2015}. Additionally, we mention the image inpainting, i.e., the filling in of damaged or missing regions of an image with the use of information from surrounding areas, \cite{Bertozzi07,Cherfils17,Reshma14}, micro-structure with elastic inhomogeneity which determines the transformation path and the corresponding microstructure evolution, \cite{Shenyang04} and references therein, tumor growth simulation in order to provide optimal strategies for treatments, \cite{Agosti17,Signori2019}, and topology optimization phase-field approach to the problem of minimizing the mean compliance of a multi-material structure, see e.g. \cite{Zhou2006,Blank2012,REGAZZONI18}.

The investigation of the behavior of such systems, started from the pioneer work \cite{CaHi58} and the early works of \cite{NOVICKCOHEN1984277,Elliott1986,Elliottt1989}. Thereafter, a lower solution space regularity and a second order splitting method {{have been investigated}} in \cite{Elliott1989}, the solution existence and error analysis in \cite{ElioLar92,KaraliNagase14}, for higher order finite element we refer to \cite{CrossGoudenege2012} and for Discontinuous Galerkin in space approach to \cite{WELLS2006860,LI20182100}. Optimal control for non-convective
 or optimal control for the convective case have been studied in \cite{RoccaSprekels15conv_control,conCHZhao2014,convectiveOptimal_Hinze17,controldistrCHiWe12,CoFaGiSpre15,HiKeWe17}, stochastic partial differential equations and stochastic analysis in the very new work of  \cite{FuKovacsLarssonLindgren18}, Navier-Stokes/Cahn-Hilliard systems in \cite{HinzeNS13,HintermullerHinze2013Ellipse}, and Cahn-Hilliard/Allen-Cahn systems in \cite{ANTONOPOULOU20162383,Antonopoulou2013EXISTENCEAR,Zhang2018,Alikakos2019}.

 Throughout this work, an efficient methodology for solving nonlinear systems governed by  Cahn-Hilliard equations is studied within cut finite elements and reduced-order modeling. {\magenta{A stable fully discrete cut finite element numerical scheme for this geometrically parameterized nonlinear fourth-order diffusion system is introduced while a splitting approach, transforming the fourth-order equation into a coupled system of two second-order equations, is considered}}. In this embedded geometry framework, we {\blue{propose}} a model order reduction technique using the advantages of a shape regular background mesh, recently investigated in \cite{KaBaRO18,KaratzasStabileNouveauScovazziRozza2018,KaratzasStabileAtallahScovazziRozza2018}. The combination of unfitted mesh finite element methods and reduced order modeling allows us to obtain a fast evaluation, {\blue{by considering}} the geometrical parametrized system, while we avoid {\blue{remeshing,}} as well as the reference domain formulation, often used in boundary fitted finite element formulations.

This contribution  is organized as follows: 
 In Section \ref{sec:HF} we define the continuous strong formulation of the mathematical problem. The semidiscrete cut elements Nitsche weak formulation and the {\blue{implicit explicit Euler method (IMEX)}} fully discrete problem under consideration is introduced, as well as, the incremental scheme used to solve the full order problem during the offline stage.
 In Section \ref{sec:ROM}, we demonstrate the reduced-order model formulation, {\blue{based on}} the Proper Orthogonal Decomposition {\magenta{(POD)}}, and its main aspects.
 In {\magenta{Section}} \ref{sec:num_exp} the IMEX method and the high fidelity solver efficiency is validated. Subsequently, the proposed ROM technique is tested on a geometrical parametrized problem of a two-phase field problem starting from pseudo-random initial data around an embedded circular domain. Convergence results, errors and reduced execution times are introduced {\blue{and analyzed}}.
 Finally in Section \ref{sec:conclusions}, conclusions and perspectives for future improvements and developments are demonstrated.
To our best knowledge, the results of this work are {{original}} and applicable in many cases of nonlinear time-dependent partial differential equation problems in terms of the prescribed methodology and in the spirit of the geometrical parametrization.

\section{The {\magenta{model}} problem and the {\magenta{full order}} approximation} \label{sec:HF}
\subsection{Strong formulation of the Cahn-Hilliard problem}
We consider the model problem describing a phase {\magenta{field}} time evolution. {\magenta{The unknown {\Xgreen{variable $u$, often}} considered as phase field variable, 
is related to the concentration via $u = (u_A-u_B)/(u_A+u_B)$ with species A,B.}} 
Let us consider an open {\magenta{bounded}} domain $\Omega$ in ${\mathbb R}^d$, with {\magenta{$d=2,3$ and}} Lipschitz boundary $\partial\Omega$. 
{\Xgreen{Let $\mathcal P$ be a $k-$dimensional  parameter space  with}} a parameter vector $\mu\in \mathcal P \subset \mathbb{R} ^k$. We state below, in a time interval $[0,T]$,  the strong form of the evolutionary Cahn-Hilliard phase {\magenta{field}} system of equations with {\magenta{Neumann}} boundary conditions on $\partial\Omega$, geometrically {\blue{parametrized}} by $\mu$. We denote by $\Omega(\mu)$ and $\partial\Omega(\mu)$ the parametrized domain and {\blue{boundary,}} respectively. 

As first suggested by {\magenta{\cite{CaHi58}}}, and thereafter extended in {\magenta{\cite{DeGroot62}}}, if we assume that the mobility is equal to $1$ and $\varepsilon$ is a measure of the size of the interface of two fluids, the mass flux is given by
\begin{equation}\label{eq:flux}
{{\bf{J}}(\mu)}=-\nabla \left( \frac{{1}}{\varepsilon^2}{\bblue{F'}}(u(\mu)) - {{\varepsilon^2}} \Delta u(\mu)\right),
\end{equation}
where $F$ denotes the {\magenta{free energy}}. {\magenta{According}} to \cite{CaHi58}, the Ginzburg--Landau energy becomes
\begin{eqnarray}\label{Eq:free_energy}
E(u(\mu))=\int_{\Omega}\left(F(u(\mu))+ \frac{{\varepsilon^2}}{2}|\nabla u(\mu)|^2\right)d{\bf{x}}.
\end{eqnarray}
An equilibrium state of the considered mixture minimizes the above Ginzburg--Landau energy, subject to the mass conservation 
\begin{eqnarray}\label{Eq:mass_conservation}
{\magenta{\frac{d}{d t}}}  \int_{\Omega(\mu)}u(\mu)d{\bf{x}} = 0.
\end{eqnarray}
Hence, the parametrized Cahn-Hilliard system can be described as:
\begin{eqnarray}\label{Eq:4thOrder}
\frac{\partial{u}(\mu)}{\partial{t}} = -\varepsilon^2 \Delta ^2 u(\mu) + \frac{1}{\varepsilon^2}\Delta F'(u(\mu)) ,  &&\text{ in } {\Omega}(\mu)\times[0,T],
\\\label{Eq:4thOrder_boundary}
\partial_n u(\mu) = \partial _n (-\varepsilon ^2 \Delta u(\mu) + \frac{1}{{\Xgreen{\varepsilon^2}}}F'(u(\mu))) = g_N(\mu), 
&&\text{ on } {\partial \Omega}(\mu)\times[0,T],
\\ \label{Eq:4thOrder_initial}
u(\cdot, 0)= u_0 (\cdot),    &&\text{ in } \Omega(\mu),
\end{eqnarray}
where $n$ is the unit outer normal vector of $\partial\Omega$, and $F$ {\magenta{is a double well free energy often taken as}} a polynomial
function of $u$ of fourth power:  
\begin{equation}\label{double-well}
 F(u(\mu)) = \gamma_2 \frac{u^4(\mu)}{4} + \gamma_1 \frac{u^3(\mu)}{3} +\gamma_0\frac{u^2(\mu)}{2} \text{ with } \gamma_2 > 0{\bblue{,}}
\end{equation}
{{{\magenta{with 
$F'(u(\mu))$ to be the cubic {\Xgreen{expression}} ${\Xgreen{F'(u(\mu))=}}\gamma_2 {u^3(\mu)} + \gamma_1 {u^2(\mu)} +\gamma_0{u(\mu)}$}}, {\magenta{and $g_N(\mu)$ is the Neumann boundary data}}.}}

For more details, the interested reader is referred {\magenta{to \cite{Elliott1989} 
as}} well as for Dirichlet boundary conditions in \cite{LI2013102}.
We remark that the equation (\ref{Eq:4thOrder}) represents the equation of conservation of {\blue{mass,}} with mass flux as described in (\ref{eq:flux}).

\subsection{{\magenta{Full order}} parametrized Nitsche cut elements weak variational formulation}
\subsubsection{A proper continuous weak formulation}
The Cahn-Hilliard equation as it is expressed in {\magenta{equations}} (\ref{Eq:4thOrder})--(\ref{Eq:4thOrder_initial}) is a fourth-order diffusion equation, involving first-order time derivatives, second and fourth-order spatial derivatives. Casting it in a weak form results in second-order spatial derivatives avoiding the fourth-order {\magenta{ones. 
Another}} setback is that if one employs the Nitsche weak boundary enforcement, several integrals in the cut geometry should be calculated including various order derivatives and normal derivatives, \cite{WELLS2006860,Nitsche_4thorder_EmDoHa2010,ZHAO2017177,Harari_Nitsche15}, which in our case of unfitted mesh are avoided due to time expensive integration. 
{\bblue{We also recall that in general, the efficiency of systems arising from CutFEM discretization schemes suffers from a system condition number that depends on mesh and boundary intersection position and cannot be handled by small time -or space- steps. The latter force us to use a proper boundary interface stabilization and to use 
%
a special stabilization term on the cut elements interface for the nonlinearity,}} 
 \cite{CLAUS2019185,BURMAN20101217}. 
To overcome all these difficulties arising from this fourth-order nonlinear system, firstly we employ a splitting method deriving a strong formulation system that requires $H^2$ space regularity. Afterwards, we multiply with a test function and we integrate by parts over $\Omega$ in order to drive the system to a weak form that requires only $H^1$ space regularity. In particular the pair solution $(u(\mu), w(\mu))$ solves the following problem: \orange{find} $u(\mu) \in L^2[0,T;H^1 (\Omega(\mu)\orange{)}]\cap \orange{H^1}[0,T;(H^1 (\Omega(\mu))\orange{)}']$ and $w(\mu) \in  L^2[0,T;H^1 (\Omega(\mu)\orange{)}]$  for all {\magenta{test}} functions $v(\mu) \in H^1 (\Omega(\mu))$ 
such that
\begin{eqnarray}\label{Eq:SplittedCH1}
 {\magenta{\langle}}\frac{\partial{u(\mu)}}{\partial{t}} , v(\mu){\magenta{\rangle}} +  (\nabla w(\mu), \nabla v(\mu)) = \langle g_N(\mu), v(\mu)\rangle_{\partial \Omega(\mu)},&&\qquad
\label{Eq:SplittedCH2}
\\
- (w(\mu), v(\mu)) + \varepsilon^2 (\nabla u(\mu), \nabla v(\mu)) + \varepsilon^{-2} (F'(u(\mu)), v(\mu)) = \langle g_N(\mu), v(\mu)\rangle_{\partial \Omega\orange{(\mu)}},&&\qquad
\label{Eq:SplittedCHr_initial}
\\
u(\cdot, 0;\mu)= u_0 (\cdot;\mu),&&\qquad
\end{eqnarray}
where $w(\mu)$ {\magenta{depends}} on the geometry parameter $\mu$ and usually it is identified as a chemical potential, {\magenta{and the initial concentration}}  $u_0\in {{L^2(\Omega(\mu))}}$. {\magenta{We have used the standard notation $(\cdot,\cdot), \langle\cdot,\cdot\rangle _{\partial\Omega}$ for the $L^2(\Omega)$ and $L^2(\partial\Omega)$ inner products respectively,
while for 
the corresponding duality pairing $\langle{\cdot,\cdot}\rangle _{H^1(\Omega)',H^1(\Omega)}= \langle{\cdot,\cdot}\rangle$}}. In the following, we  are focused on the concentration component $u(\mu)$, and we handle $w(\mu)$ as auxiliary function. Similarly, Dirichlet boundaries will be examined following \cite{LI2013102}.
\begin{remark}
We point out that we are taking into account nonsmooth initial data ${L^2(\Omega(\mu))}$ while the aforementioned regularity $\orange{H^1}[0,T;(H^1(\Omega(\mu))\orange{)}']$ of $u(\mu)$ is needed for the ${\partial{u(\mu)}}/{\partial{t}}$ term\footnote{\orange{namely, ${\partial{u(\mu)}}/{\partial{t}}\in {L^2}[0,T;(H^1(\Omega(\mu)){)}']$.}}, 
see \cite{doi:10.1137/130943108,ChKa2015}. 
\end{remark}
\subsubsection{Discretization, unfitted mesh and stability issues}
{\magenta{A fixed background mesh is applied in the background geometry including the embedded disc area. Although, we do not solve inside the embedded geometry, while a smooth extension of the solution is present in the boundary interface coming from the boundary interface stabilization which causes neighboring ghost elements with values to decrease smoothly to zero. Compare for instance in Section \ref{sec:num_exp}, the cut geometry in Figure \ref{FULL_3D} and the uncut geometry zoomed image in Figure \ref{fig:poisson_zoom}. This approach guarantees a regular solution in the background domain and permits the construction of a reduced basis with better approximation properties, \cite{KaBaRO18}. Finally, we cut the reduced solution onto the truth geometry, which actually is the circular line.
 We clarify that we do not examine the interfacial area between two bulk phases while we experiment with the full order/reduced basis approach for geometrical deformations.
}}

The system (\ref{Eq:SplittedCH1})--(\ref{Eq:SplittedCHr_initial}) in a discrete unfitted mesh formulation needs extra attention. {\magenta{Stability issues appear, in the sense that the solution is exploding in the interior of the domain and/or in the boundary interface, a phenomenon which is strengthened by the fact that the initial condition we use is not energetically favorable. Even if one uses classical finite element methods {\Xgreen{and one strongly}} applies the boundary conditions, to achieve a  stable solution with small errors a quite small time-stepping is necessary, see e.g. \cite[page~9
]{CrossXU2019524}
.}} We highlight that the evolution of the physics of the problem during the first time steps is very fast, although, as time passes it slows down and finally it equilibrates, see also e.g. \cite{WODO20116037,BoPhD16}, and references therein. An adaptive time-stepping approach would appear beneficial, although we will investigate it in a future work, as well as, the way it affects the accuracy and efficiency of the reduced model detailed in Sections \ref{sec:ROM} {{and}} \ref{sec:num_exp}. 
Considering the Nitsche terms one may apply the simple formulation needed for linear systems, {\magenta{i.e without the use of any kind of Newton iterative method}}. 
Under these considerations, we used cut finite element methods to solve the system applying a jump stabilization procedure in the boundary elements interface area.
In the next paragraph, we derive the semi-discrete formulation of the  system, while in paragraph \ref{Time discretization}
the fully discrete system is exploited. 
\subsubsection{Semidiscrete variational formulation} 

\begin{figure}
\centering
\begin{minipage}{\textwidth}
\centering
{\footnotesize{(i)}}\hskip1.0pt
\includegraphics[width=0.2\textwidth]{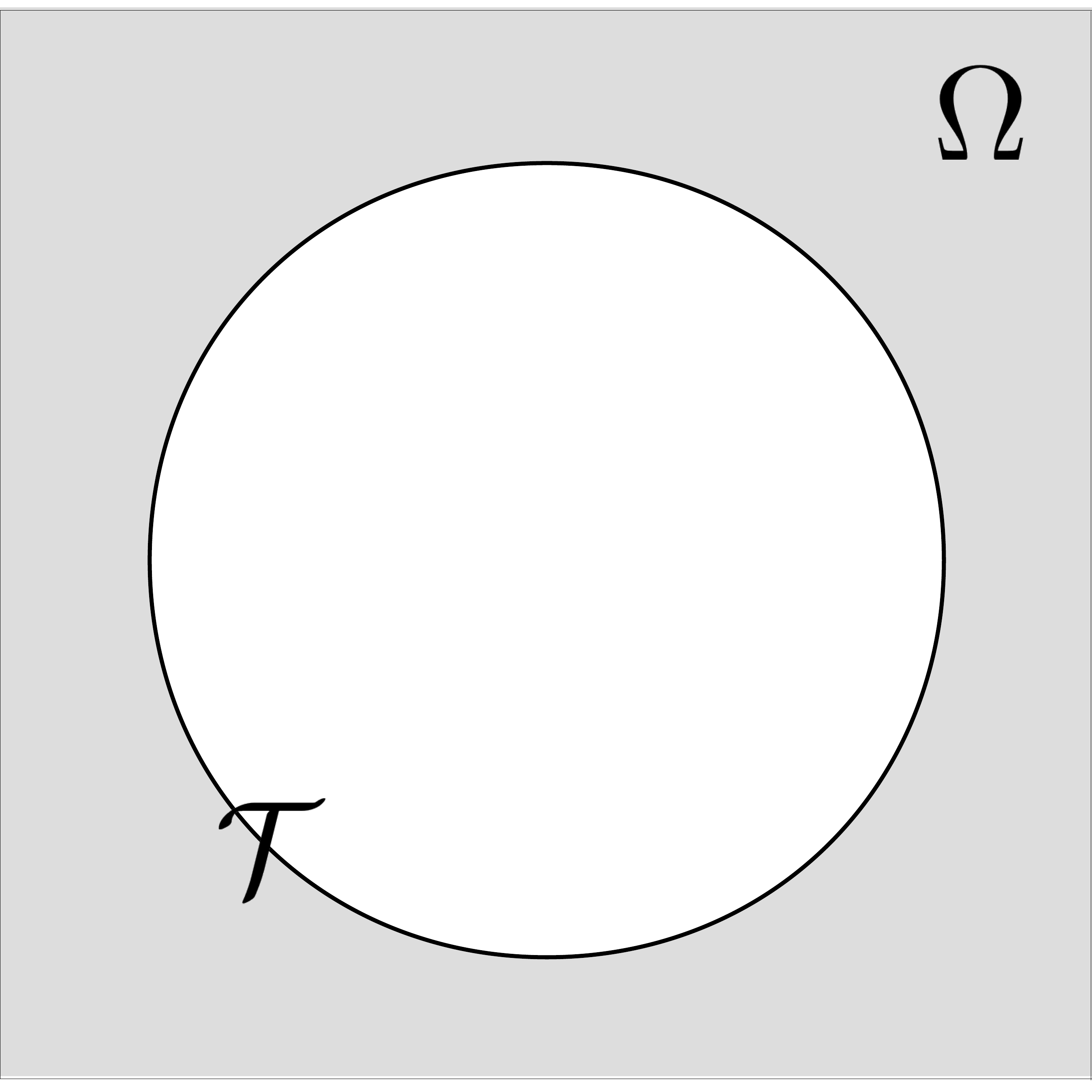}
~~~~~~
\hskip-2.pt{\footnotesize{(ii)}}\hskip0.25pt\includegraphics[width=0.195\textwidth, trim={0 0 0 0.1cm},clip]{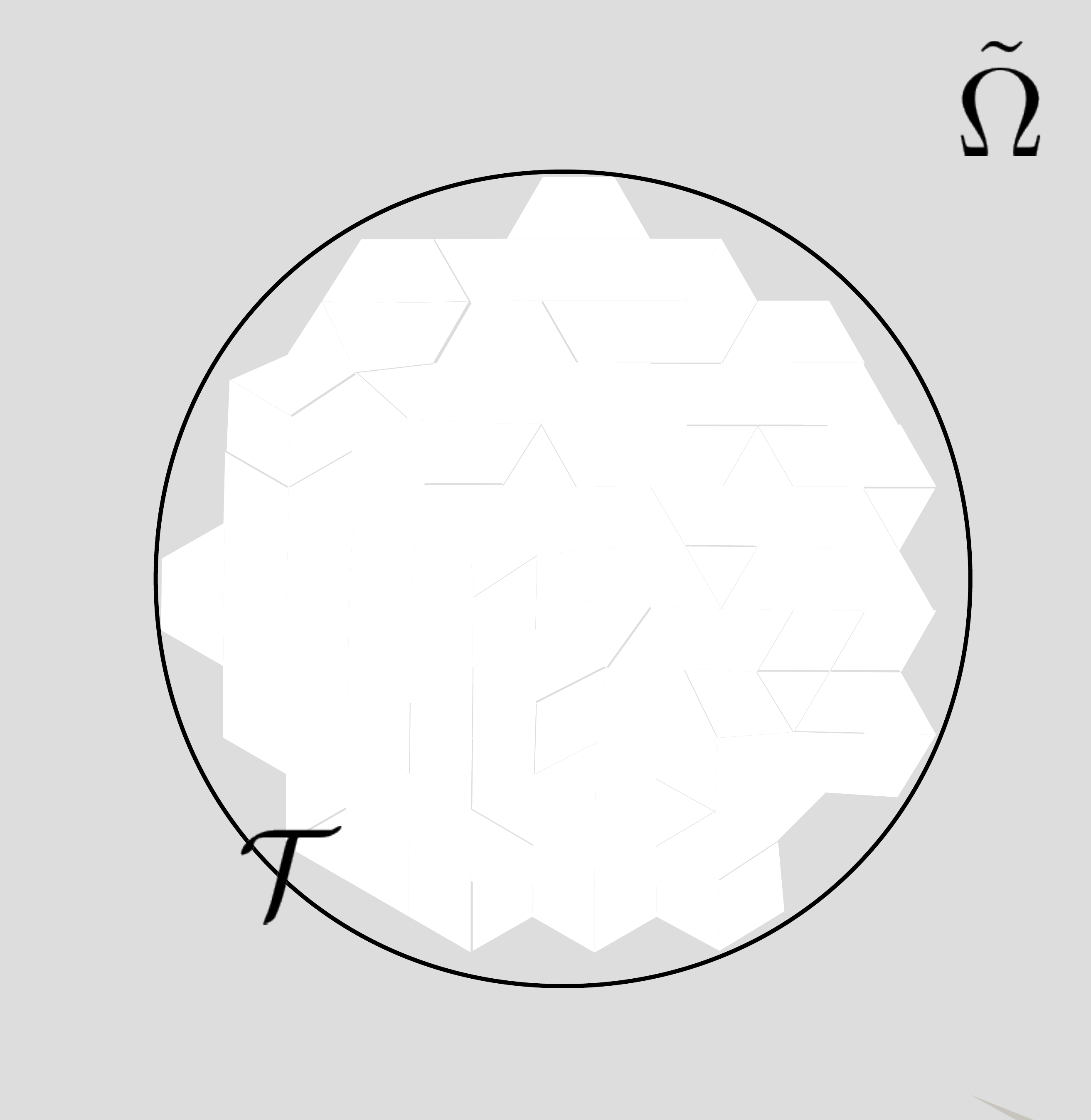} 
~~~~~~
\\
\hskip6.0pt{\footnotesize{(iii)}}\hskip1.0pt
\includegraphics[width=0.2\textwidth, trim={0 0 0 0.8cm},clip]{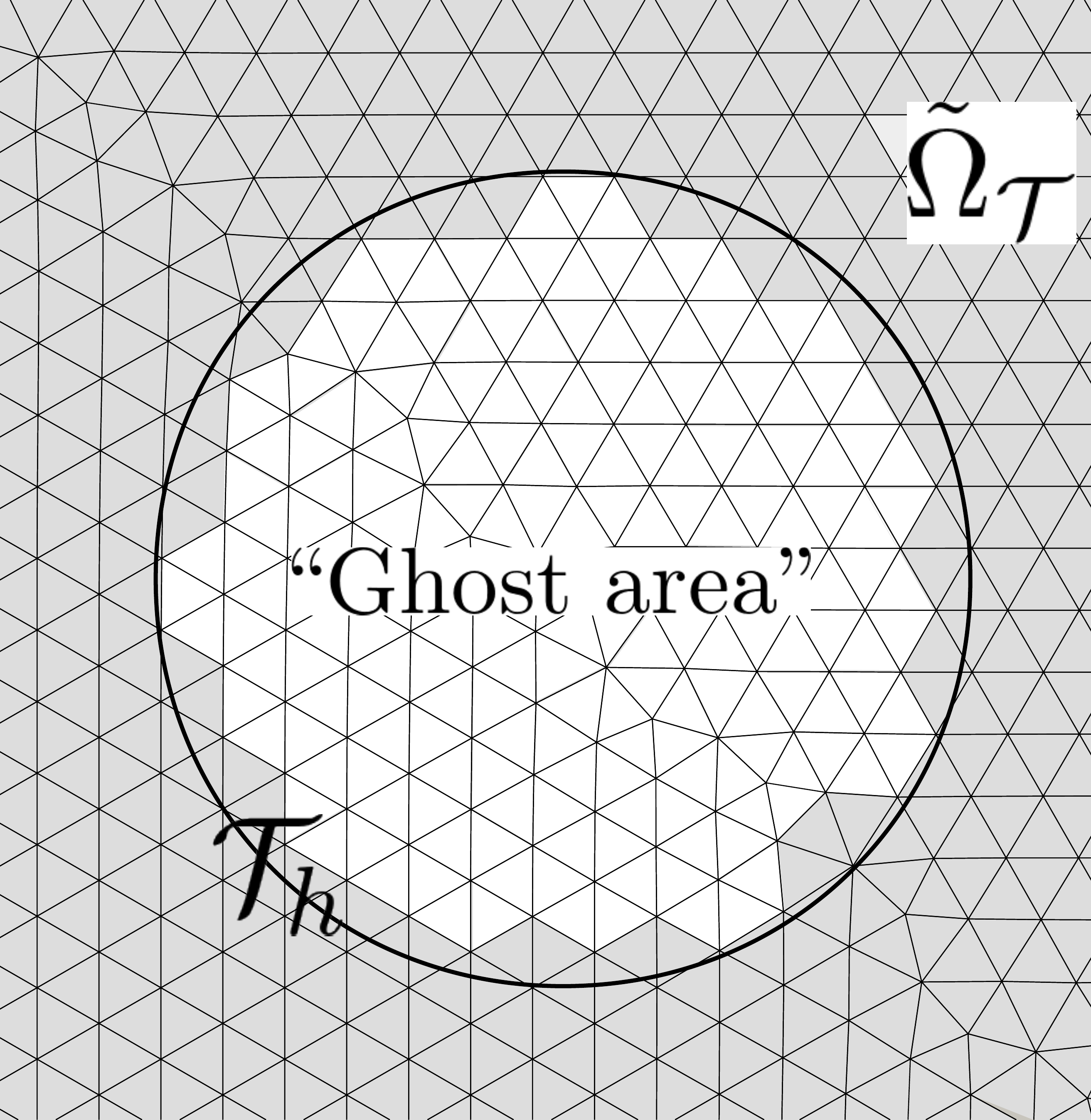}
~~~~~~
\hskip-1.pt{\footnotesize{(iv)}}\hskip-1.pt
\includegraphics[width=0.27\textwidth]{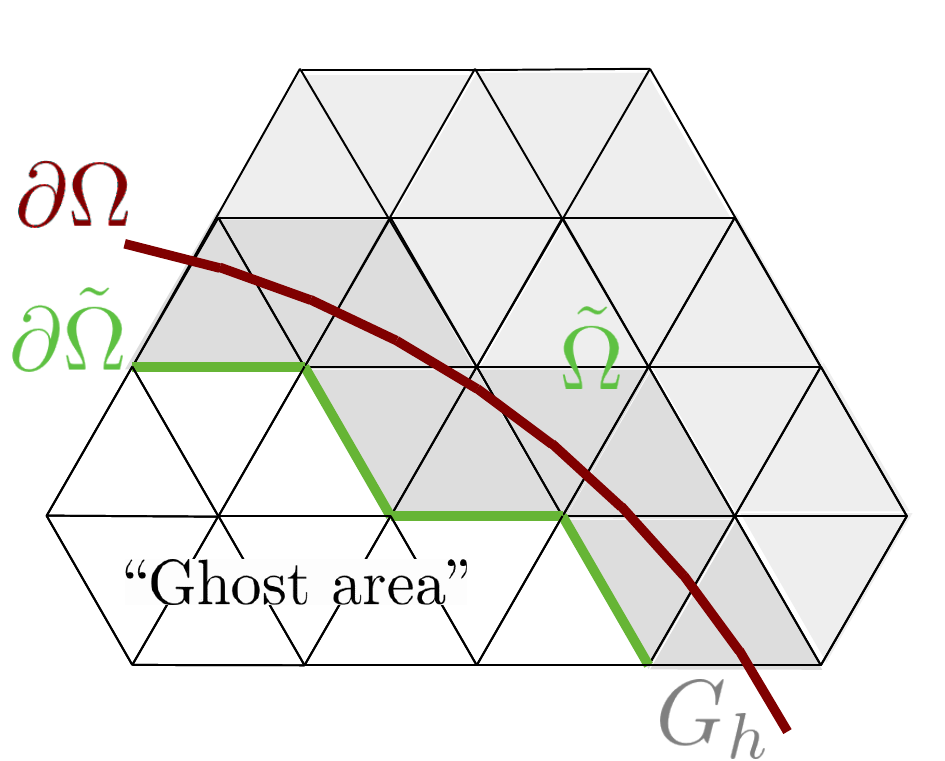} 
\end{minipage}
\caption{(i) The geometry of an embedded disk, (ii) the cut finite element method geometry, (iii) zoom on the background mesh together with the surrogate cut discretized geometry, and (iv) extended mesh and elements intersected by the true boundary.}
\label{SurrogateMesh}
\end{figure}

We denote by $\mathcal{T}$ the background domain, and by $\mathcal{T}_h$ its corresponding mesh, see e.g. Figure \ref{SurrogateMesh}.
%
{\magenta{We consider a family of triangulations (say $\{\mathcal{T}_h\}_{h>0}$) of $\Omega$, defined in the standard way. 
To every element $K\in \mathcal{T}_h$, we associate two parameters $h_K$ and $\rho_K$, denoting the diameter of the set $K$, and the diameter of the largest ball contained in $K$ respectively. The size of the mesh is denoted by $h= \max _{K\in{\mathcal{T}}_h} h_K$. 
The following standard properties of the mesh will be assumed:
(i) There exist two positive constants $\rho_{\mathcal{T}}$ and $\delta_{\mathcal{T}}$ such that $\frac{h_K}{\rho_K} \leq \rho_\mathcal{T}$ and $\frac{h}{h_K}\leq \delta_\mathcal{T}$, $\forall K\in {\mathcal{T}}_h$ and $\forall h >0$. 
(ii) Given $h$, let $\{K_j\}^{N_h}_{j = 1}$ denote the family of triangles belonging to $\mathcal{T}_h$ and having one side included on the boundary $\partial\Omega$. Thus, if the vertices of $K_j\cap \partial\Omega$ are denoted by $x_{j,\partial\Omega},x_{{j+1},\partial\Omega}$ then the straight line $[x_{j,\partial\Omega}, x_{{j+1},\partial\Omega}] \equiv K_j\cap \partial\Omega$. Here, we also assume that $x_{1,\partial\Omega}= x_{N_{h+1},\partial\Omega}$.}}
We recall the standard notation 
$(\cdot, \cdot)
$, 
$\langle \cdot, \cdot \rangle _{\partial\Omega}$ for the $L^2({\Omega})$, 
and $L^2(\partial\Omega)$ inner products {\magenta{associated with}} the truth geometry 
${\Omega}$, 
and $\partial\Omega$, respectively.
 The continuous boundary value problem is next formulated on a domain $\tilde{\Omega}(\mu)$ that contains ${\Omega}(\mu) \subset \tilde{\Omega}(\mu)$, while its mesh $\tilde{\Omega} _{\mathcal{T}}(\mu) := {\mathcal{T}}_h\cap \tilde{\Omega}(\mu)$ is not fitted to the domain boundary {\magenta{$\partial\Omega$}} and ${\tilde{\Omega}}_{\mathcal{T}}(\mu) \subset \mathcal{T}$ 
 for all $\mu \in \mathcal{K}$. Let also $\Omega_{\mathcal{T}}(\mu) := {\mathcal{T}}_h\cap {\Omega}(\mu)$ and ${G}_h(\mu):=\{K\in \mathcal{T}_h(\mu): K\cap\partial\Omega(\mu)\neq\emptyset\}$ be the set of elements that are intersected by the interface. We remark that  ${G}_h(\mu)$ and ${\tilde{\Omega}}_{\mathcal{T}}(\mu)$ depend on $\mu$ through ${\tilde{\Omega}(\mu)}$ (or its boundary), while the background domain $\mathcal{T}$ and its mesh $\mathcal{T}_h$ {do not} depend on $\mu$. 
Furthermore, the set of element faces $\mathcal{F}_G(\mu)$ associated with ${G}_h(\mu)$, is defined as follows: for each face $F\in \mathcal{F}_G(\mu)$, there exist two simplices $K \neq K'${{, such that $F=K\cap K'$,}} and at least one of the two is a member of ${G}_h(\mu)$. Note that the boundary faces of ${\tilde{\Omega}} _{\mathcal{T}}(\mu)$ are excluded from $\mathcal{F}_G(\mu)$.
On a face $F \in \mathcal{F}_G(\mu)$, $F=K\cap K'$, the jump of {\magenta{$v$ and the jump of the gradient of $v$ are defined by 
\begin{eqnarray}
\ldb v \rdb&=& v |_K - v|_{K'}, {\text{ and }}
\label{def:jumps}
\\
\ldb{{\bf n}_F}\cdot\nabla  v \rdb&=&{\bf n}_F\cdot\nabla v |_K -{\bf n}_F\cdot \nabla v|_{K'} \label{def:grad_jumps}
\end{eqnarray}
respectively}}, where ${\bf n}_F$ denotes the outward pointing unit normal vector to $F$ {\magenta{and $v$ is in the same spaces as $u_h,w_h,v_h,q_h$ defined  in equation (\ref{eq:semidiscrete})}}. {\Xgreen{We denote the space of continuous piecewise-linear functions by ${{{\mathcal{V}}_h}} ({{{\tilde{\Omega}}}}(\mu))$
}}
 \begin{equation}
{\mathcal{V}}_h({\magenta{\tilde{\Omega}}}(\mu))=\left \{  \upsilon \in C^0({\overline {\Omega}}_{\mathcal T}(\mu))\,:\,\upsilon |_K  \in P^1(K), \, \forall K\in \mathcal{T}_h(\mu)    \right \}.
\label{eqn:V_poisson}
\end{equation}
For the sake of simplicity, in the next set of equations, we will omit the parameter dependency notation with respect to $\mu$ and the {\magenta{cut finite element method (CutFEM)}} discretization is as follows. 
We seek  ${{u}}_h, w_h\in {{{\mathcal{V}}_h}} ({\orange{\tilde{\Omega}}}(\mu))$
 such that, for {\magenta{all test}} functions $v_h$,  $q_h\in {\mathcal{V}}_h(\orange{\tilde{\Omega}}(\mu))$ {\magenta{it holds}}
\begin{eqnarray} 
({\magenta{(u_h)_t}},v_h) 
+\orange{(\nabla w_h, \nabla v_h)} 
+\orange{\varepsilon^2(\nabla u_h,\nabla q_h)}
-(w_h,q_h)
+
\frac{1}{\varepsilon ^2}( \gamma_2 {u_h^3} + \gamma_1 {u_h^2} +\gamma_0{u}_h,q_h)  
\nonumber\\
+\langle\alpha_N h {\bf n}_{\Omega}\cdot \nabla u_h,{\bf n}_\Gamma\cdot \nabla v_h \rangle_{\partial\Omega}
 +
        \sum_{F\in \mathcal{F}_G}\left(\alpha_1 h^{-2}\ldb u_h\rdb,\ldb v_h\rdb\right)_F \label{eq:semidiscrete}\nonumber \\
 \nonumber =\langle g_N, v_h + q_h +\alpha_N h {\bf n}_\Gamma\cdot \nabla v_h\big\rangle_{\partial\Omega},
\\ \label{eq:semidiscrete}
\end{eqnarray}
where   $\alpha_N$, and $\alpha_1$  are positive penalty parameters related to Nitsche weak imposition of boundary conditions and the boundary interface stabilization term respectively, see for instance \cite{BuHa11}.
\begin{remark}
For the sake of completeness, it would be convenient to emphasize that {both} the aforementioned kind of {\magenta{jumps as they are defined by (\ref{def:jumps}) and (\ref{def:grad_jumps})}}, can be used to apply different type of ghost penalty stabilizations, namely  
$\sum_{F\in \mathcal{F}_G}\left(\alpha_1 h\ldb{{\bf n}_F}\cdot\nabla  u_h\rdb,\ldb{{{\bf n}_F}\cdot\nabla} v_h\rdb\right)_F$ 
and
$\sum_{F\in \mathcal{F}_G}\left(\alpha_1 h^{-2}\ldb u_h\rdb,\ldb v_h\rdb\right)_F$.
Although, we prefer to employ the  {\magenta{projection based}}-jump and not the {\magenta{derivative}}-jump one. 
 All experiments consider first order polynomials. 
Finally, we highlight that the preferred ghost penalty is less computational expensive since it does not involve any kind of derivatives. Concerning literature, most theoretical estimates are related to the {\magenta{derivative-jump}}, nevertheless, they can easily be extended to the {\magenta{projection based}}-jump ghost penalty since the jump can be bounded from above with the derivative jump, for more details we refer to \cite{BURMAN20101217,BuHa14,Schott2016,Lehrenfeld2018L2errorAO}.
\end{remark}
 \subsubsection{{\magenta{Implicit-explicit (IMEX)}} type time discretization}\label{Time discretization}
 The goal here is to find concentration $u_h = u_h(x, t)$ satisfying {\magenta{equation}} (\ref{eq:semidiscrete}) for all time instances $t$ in a time interval $[0,T] \subset {\mathbb{R}}$ and space positions $x \in \Omega \subset {\mathbb{R}}^{d}$.
We fully discretize the system by an IMEX approach, see e.g. \cite{RokhzadiPhD18} and references therein, as well as \cite{KaKaTra20} for similar types of nonlinearities. Approximations will be constructed on a time partition $0=t^0 < t^1 < \ldots < t^{\magenta{N_t}}=T$  on which, each interval $I_n : = (t^{n},t^{n+1}]$  is of length $\tau_n=t^{n+1}-t^{n}$, {\magenta{$n=0,...,{\magenta{N_t}}-1$}}, starting from {\magenta{the}} initial condition $u^0$.
Therefore, 
we apply a splitting method onto time integration level --often called operator splitting method-- which means that the differential operator is rewritten as the sum of two complementary operators.  The latter extensively attracted attention  IMEX method technique treats on the nonlinear term explicitly, allowing it to act as a forcing term in the $w$-equation (\ref{Eq:SplittedCH2}). In this way, we avoid stability issues caused by the nonlinearity {\magenta{and complicated Nitsche boundary enforcement on the boundary interface area, see for example \cite[page~153]{Schott2016} and references there in}}.
All the other terms have been handled implicitly for {\bblue{increased stability}}. 
Hence, the IMEX method results in the set of fully-discrete Cahn-Hilliard equations:

Find $(u^{n+1}_h,w^{n+1}_h)\in {\mathcal{V}}_h\times {\mathcal{V}}_h$, {\magenta{such that}} for all $(v^n_h,q^n_h)\in {\mathcal{V}}_h\times {\mathcal{V}}_h$
 $$
 A(u^{n+1}_h) + \tau_n \cdot L(u^{n+1}_h,w^{n+1}_h) = A(u^n_h) - \tau_n\cdot N(u^n_h) + B(g_N) \orange{\cdot \tau_n},
 $$
 where
\begin{eqnarray*}
 A(u^{n+1}_h)&=&(u^{n+1}_h,v^{n+1}_h)={\bf{A}}u^{n+1}_h,\\
 L(u^{n+1}_h,w^{n+1}_h)&=&(\nabla w^{n+1}_h,\nabla v^{n+1}_h) + \varepsilon^2(\nabla u^{n+1}_h, \nabla q^{n+1}_h) -(w^{n+1}_h,q^{n+1}_h) \\
&&+\langle\alpha_N h {\bf n}_{\Omega}\cdot \nabla u^{n+1}_h,{\bf n}_\Gamma\cdot \nabla v^{n+1}_h \rangle_{\partial\Omega}
 + \sum_{F\in \mathcal{F}_G}\left(\alpha_1 h^{-2}\ldb u_h\rdb,\ldb v_h\rdb\right)_F
 \\
 &&= {\bf{L}_1}u^{n+1}_h + {\bf{L}_2}w^{n+1}_h + {\bf{C_{\mathcal{F}_G}}}u^{n+1}_h,
 \\
B(g_N) &=& \langle g_N, v^{n}_h + q^{n}_h +\alpha_N h {\bf n}_\Gamma\cdot \nabla v^{n}_h\big\rangle_{\partial\Omega} = {\bf{B}}_v + {\bf{B}}_q,
 \\
N(u^n_h)&=&\frac{1}{\varepsilon^2}(\gamma_2 {\magenta{{{(u^n_h)}^3}}} + \gamma_1 {\magenta{{(u^n_h)}^2}} +\gamma_0{u^n_h},v^n_h)={\bf{N}}(u^{n}_h)u^{n}_h.
\end{eqnarray*}
{\bblue{In the above formulation, we denote by ${\bf{L}}_1$ and ${\bf{L}}_2$ the operators related to the concentration $u^{n+1}_h$ and the auxiliary quantity $w^{n+1}_h$ respectively, while we have seperated the ghost penalty stabilization which is denoted by ${\bf{C_{\mathcal{F}_G}}}u^{n+1}_h$. The term $\alpha_N h {\bf n}_\Gamma\cdot \nabla v^{n}_h$ in the right hand side is the Nitsche penalty and it is related to the operators ${\bf{B}}_v$, ${\bf{B}}_q$.}}
{\magenta{Next, we derive the matrix form related to that system}}. We define the {\magenta{parameter-dependent}} Cahn-Hilliard operator 
$$
G(U^{n}_h(\mu)):= G\left( \begin{bmatrix} u^{n}_h(\mu) \\ w^{n}_h(\mu) \end{bmatrix}
  \right)
  =
  \begin{bmatrix} \orange{ 
  {\bf{L}_1}} +  {{\bf{N}}}({{u}}^{n}_h(\mu)) & {\bf{L}_2}  
  \\ {\bf{C_{\mathcal{F}_G}}}
  & 0  
  \end{bmatrix}
  \begin{bmatrix} {{u}}^{n}_h(\mu) \\ w^{n}_h(\mu) \end{bmatrix}{\blue{,}}
$$
and the right hand side consists of the forcing boundary data related to stabilization and Nitsche weak enforcement boundary terms $ F_N(\mu):= \begin{bmatrix} {\bf{B}}_v (\mu) \\ {\bf{B}}_q(\mu) \end{bmatrix}$. 
These definitions result in the following residual
$$
R(U^{n}_h(\mu))=G(U^{n}_h(\mu))
 - F_N(\mu), 
$$
which {\magenta{yields}} the following algebraic system of equations for the increment $\delta U^{n+1}_h(\mu) = 
\orange{\begin{bmatrix} \delta u^{n}_h(\mu) \\ \delta w^{n}_h(\mu) \end{bmatrix}}
= U_h^{n+1}(\mu)-U_h^{n}(\mu)$:
\begin{eqnarray}\label{eq:system_linear0} 
\orange{\begin{bmatrix} {{\bf{A}}} + 
  \tau_n{\bf{L}_1}  & \tau_n{\bf{L}_2}  
  \\ \tau_n{\bf{C_{\mathcal{F}_G}}}
  & 0  
\end{bmatrix}}
 \delta U^{n+1}_h(\mu) = -\tau_nR(U^{n}_h(\mu)),
\end{eqnarray}
%
noting that in the above formulation the nonlinear term is treated only explicitly and the remaining part implicitly. 
\begin{remark}
IMEX method approach is beneficial since we avoid the extra computation of iterative approaches, e.g. Newton-Raphson method, and simultaneously we avoid the small-time stepping of an explicit time integration method, needed for a stable solution in order to minimize the dispersion error associated with these schemes {\magenta{\cite{eyre_1998,RokhzadiPhD18}}}.
In the above system of equations, it is important to underline that the discretized differential operators $\bm{A}$, ${\bf{L}_1}$, ${\bf{L}_2}$ and ${\bf{C_{\mathcal{F}_G}}}$ are parameter-dependent, a feature of great importance as we will see in Section \ref{sec:ROM} and the ROM basis construction. Also, a pre-assembling technique for the involved matrices will be employed in Section \ref{sec:num_exp}, {\blue{by minimizing}} the time-consuming integration of several involved inner products. 
\end{remark}
\section{Reduced order model with a POD-Galerkin method}\label{sec:ROM}
{{In this paragraph, a POD-Galerkin approach is briefly recalled as in \cite{HeRoSta16,Rozza2004ReducedBM}. We emulate the high fidelity model system with a reproduced one, which allows predictive errors, within {\blue{the aim of}} a reduced computational {\blue{cost}} and solution {\blue{time}} in a way adjusted to embedded-immersed boundary finite element methods. This reduced-order system has been approved advantageous when geometrically deformed systems appear and in comparison with traditional finite element methods and/or reduced-order modeling, see for instance~\cite{BaFa2014,KaBaRO18,KaratzasStabileAtallahScovazziRozza2018,KaratzasStabileNouveauScovazziRozza2018}. In particular, we employ a projection-based reduced order model which consists of the projection of the governing equations onto the reduced basis space constructed on a fixed background mesh.

Following the literature, one could see reduced basis (RB) methods applied to linear elliptic equations in \cite{Rozza2008}, to linear parabolic equations in \cite{grepl2005} and to non-linear problems in \cite{Veroy2003,Grepl2007}. Although the number of works on reduced-order models with classical {\magenta{finite element methods (FEM)}} are now significant big including Cahn-Hilliard systems, see e.g. \cite{convectiveOptimal_Hinze17,HeRoSta16,Rozza2004ReducedBM,Rozza2008,grepl2005,Veroy2003,Grepl2007}
and references therein, to the best of the authors' knowledge, only very few research works \cite{BaFa2014,KaBaRO18,KaratzasStabileAtallahScovazziRozza2018,KaratzasStabileNouveauScovazziRozza2018,KaratzasStabileNouveauScovazziRozzaNS2019} can be found concerning embedded boundary methods on linear systems and ROMs and much fewer for nonlinear, \cite{KaratzasStabileNouveauScovazziRozzaNS2019,KaNoBaRo20}}. 
{\magenta{In this work}}, we investigate and focus on 
how we can achieve a stable solution for a {\magenta{geometrically parametrized system of second-order equations stemming from the splitting approach of}} the fourth order evolutionary Cahn-Hilliard PDE system, in a Full Order Method (FOM) and in a Reduced Order Model (ROM) framework, {\magenta{and}} within an embedded finite element method namely in a cut finite elements setting. 
%
%
%
%
%
The new techniques of \cite{KaBaRO18,KaratzasStabileAtallahScovazziRozza2018,KaratzasStabileNouveauScovazziRozza2018}, which are based on 
%
{\magenta{the combined use of a fixed background mesh for all geometrical deformations
and}} 
a  proper orthogonal {\magenta{decomposition strategy}} will be employed. The key feature of our approach is that we avoid the remeshing effort or/and the need of a map of all the deformed geometries to reference geometries often used in fitted mesh finite element methods, see e.g. \cite{HeRoSta16,RoVe07,Rozza2009,ballarin2015supremizer,RoHuMa13,Rozza2008,BeOhPaRoUr17}. 

Regarding the reduced-order modeling, we investigate how ROMs can be applied to time dependent cut finite element methods and generally, to embedded boundary methods simulations considering time dependent nonlinear systems and Cahn-Hilliard systems. The main interest is to generate ROMs on parametrized geometries. The cut elements unfitted mesh finite element method with levelset geometry description is used to apply parametrization and the reduced order techniques (offline-online). 
An important aspect is also to test the efficiency of a geometrically parametrized reduced-order nonlinear model without the usage of the transformation to reference domains, which is an important advantage of embedded methods relying on fixed background meshes. 

Before going into thorough, we specify some basics for reduced basis modeling. {\magenta{We start}} by the generation of a set of full order solutions of the parametrized problem under a parameter values random choice. The final objective of RB methods is to emulate any member of this solution set with a low number of basis functions and this is based on a two-stage procedure, the offline and the online stage, \cite{quarteroniRB2016,RoVe07,Haasdonk2008}. 
\paragraph{\bf Offline stage} In order to derive reduced-order solutions emulating the full order system, a low dimensional reduced basis is constructed based on a specific number of full order solves. This reduced basis will be able to approximate any member of the solution set to a predictive error accuracy.
{\magenta{
Predictive in the sense that the mean relative errors for a specific number of random samples which are not the same as any of the samples used in the training stage, allow these mean relative errors to remain the same after repeating the procedure and for different samples. 
}}
Hence, it is possible to project the FOM differential operators, describing the governing equations, onto the reduced basis space, applying a Galerkin projection technique and to create a reduced system of equations. The offline stage is computationally very expensive, nevertheless, it is executed only once. 

\paragraph{\bf Online stage} Thereafter, during this stage, one can, even with very few computational resources, \orange{calculate} a reproduced reduced system of equations involving any new value of the input parameters. For more details and applications, we refer to~\cite{BeOhPaRoUr17,ChinestaEnc2017} and the references therein. 

\subsection{{\magenta{Proper orthogonal decomposition}}}
The {\magenta{full order}} model, as illustrated in paragraph \ref{Time discretization},
is solved for each $\mu^k \in \mathcal{K}=\{ \mu^1, \dots, \mu^{N_k}\} \subset \mathcal{P}$ where $\mathcal{K}$ is a finite-dimensional training set of parameters chosen inside the parameter space $\mathcal{P}$. The considered problem can be simultaneously parameter and time-dependent. In order to collect snapshots for the generation of the reduced basis spaces, one needs to consider both the time and parameter dependency. For this reason, discrete-time instants $t^k \in \{ {\magenta{t^{k_1},\dots,,{\Xgreen{k_{N'_t}}}\}}} \subset [0,T]$ {\magenta{with $k_1,...,{\Xgreen{k_{N'_t}}} \in \{0,...,N\}$}} belong in a finite-dimensional training set, which is a subset of the simulation time interval and are considered as parameters. The total number of the full order snapshots is then equal to $N_s = N_k\cdot N_t$.  The {\magenta{snapshot}} matrices $\bm{\mathcal{S}_u}$ and $\bm{\mathcal{S}_w}$, 
are then given by $N_s$ {\magenta{full order}} snapshots:
\begin{gather}
\bm{\mathcal{S}_u} = [{u}(\mu^1,{\magenta{t^{k_1}}}),\dots,{u}(\mu^{{\magenta{N_k}}},{\magenta{,{\Xgreen{k_{N'_t}}}}})] \in \mathbb{R}^{N_u^h\times N_s},
\label{Su}\\
\bm{\mathcal{S}_w} = [w(\mu^1,{\magenta{t^{k_1}}}),\dots,w(\mu^{{\magenta{N_k}}},{\magenta{,{\Xgreen{k_{N'_t}}}}})] \in \mathbb{R}^{N_w^h\times N_s},
\label{Sw}
\end{gather}
where ${N_u^h}$ and ${N_w^h}$ are the number of degrees of freedom for the discrete {\magenta{full order}} solution for the concentration $u$ and the auxiliary variable $w$, respectively, {\magenta{and in our case with the fixed background mesh and linear polynomials holds ${N_u^h}={N_w^h}$}}.
%
{\magenta{
As we will see below, we can derive an efficient reduced-order method based on the sets $S_u$ and $S_w$ creating for each of them a separate basis.
}}
In order to generate the reduced basis spaces, for the projection of the governing equations, one can find in literature several techniques such as the Proper Orthogonal Decomposition (POD), the Proper Generalized Decomposition (PGD) and the Reduced Basis (RB) with a greedy sampling strategy. For more details about the different strategies the reader may see  \cite{HeRoSta16,Rozza2008,ChinestaEnc2017,Kalashnikova_ROMcomprohtua,quarteroniRB2016,Chinesta2011,Dumon20111387}. In this work, the POD strategy is applied onto the full {\magenta{snapshot}} matrices.
{\magenta{
We clarify that in this work, we do a first investigation proving numerically that the reduced basis approach for such systems in an embedded finite element framework and a fixed background mesh is stable and with acceptable errors. Whether such kind of study with the proper orthogonal decomposition gives good results it is strongly promising 
that we can go to the next step applying the aforementioned approaches.
}}
The aforementioned procedure includes both time and parameter dependency. In the case of parametric and time-dependent problems also other approaches are available such as the POD-Greedy approach \cite{Haasdonk2008} or the nested POD {{approach,}} where the POD is applied first in the time domain and then on the parameter space.
Given a 
${u}(t;{\magenta{\mu}})$ in $\mathbb{R}^d${\magenta{, $d=2,3$}}, with a certain number of {\magenta{snapshots}} ${u}_1,\dots, {u}_{N_s}$, the POD problem consists in finding, for each value of the dimension of POD space $N_{POD} = 1,\dots,N_s$, the scalar coefficients $a_1^1,\dots,a_1^{{\magenta{N_{POD}}}},\dots,a_{N_s}^1,\dots,a_{N_s}^{{\magenta{N_{POD}}}}$ and functions ${\varphi}_1,\dots,{\varphi}_{{\magenta{N_{POD}}}}$ minimizing the quantity:
\begin{gather}\label{eq:pod_energy}
E_{N_{POD}} = {\magenta{\frac{1}{N_s}}}\sum_{i=1}^{N_s}||{u}_i-\sum_{k=1}^{N_{POD}}a_i^k {\varphi_k}||_{\magenta{L^2(\Omega)}} \hspace{0.5cm}\forall\mbox{ } N_{POD} = 1,{\blue{\dots,{\magenta{N_s}},}}\\
\mbox{ with } {\magenta{(}}{{\varphi}_i,{\varphi}_j}{\magenta{)}} = \delta_{ij} \mbox{\hspace{0.5cm}} \forall \mbox{ } i,j = 1,\dots,{\magenta{N_{POD}}} .
\label{eq:pod_energy2}
\end{gather}
{\magenta{
 The unknown coefficients ${\mathbf{a}}= [a_1^1,...,a_{N_s}^{{N_{POD}}}]
 $ are obtained through a Galerkin projection  of  the  governing  equations  onto  the  reduced  basis space for any parameter $\mu$ as we will see in equations (\ref{eq:system_linear}) and (\ref{eq:system_linear_reduced}). 
 }}
 {\magenta{We indicate here that we leave out the description for how to determine a POD basis for $w$, however, we compute a separate POD basis for the snapshot sets $S_u$ and $S_w$. 
Also, we describe above the POD basis computation for specifically $S_u$, although it is applicable for any set of snapshots.}}
 It can be shown, \cite{Kunisch2002492}, that the minimization problem of equation~\eqref{eq:pod_energy} is equivalent {\magenta{to}} solving the following eigenvalue problem:
\begin{gather}
\bm{\mathcal{C}^u}\bm{Q}^u = \bm{Q^u}\bm{\lambda^u},\\
\mathcal{C}^u_{ij} = {\magenta{(}}{{u}_i,{u}_j}{\magenta{)}} \mbox{\hspace{0.5cm} for } i,j = 1,\dots,N_s ,
\end{gather}
where $\bm{\mathcal{C}^u}$ is the correlation matrix obtained {\magenta{associated with}} the snapshots $\bm{\mathcal{S}_u}$, $\bm{Q^u}$ is a square matrix of eigenvectors and $\bm{\lambda^u}$ is a vector of eigenvalues. 
{\magenta{In practice we do not really solve a minimization problem for each $N_{POD}$ but rather we consider the decay of the eigenvalues or singular values and then decide on one specific $N_{POD}$.
}}
The basis functions can then be obtained with: 
\begin{equation}
{\varphi_i} = \frac{1}{N_s{\magenta{{\lambda_i^u}^{1/2}}}}\sum_{j=1}^{N_s} {u}_j Q^u_{ij}.
\end{equation}
The POD spaces are constructed using the aforementioned methodology resulting in the spaces:
\begin{equation}
\begin{split}
&{\bm{{\mathcal{B}}}}_u = \text{span}\{{\varphi_1}, \dots , {\varphi_{N_u^r}}\} \in \mathbb{R}^{N_{u}^h \times N_u^r},
\end{split}
\end{equation}
and similarly for the auxiliary variable: 
\begin{equation}
\begin{split}
&{\bm{{\mathcal{B}}}}_w = \text{span}\{{\chi_1}, \dots , {\chi_{N_w^r}}\} \in \mathbb{R}^{N_{{\magenta{w}}}^h \times N_w^r},
\end{split}
\end{equation}
where $N_u^r$, $N_w^r < N_s$ are chosen according to the eigenvalue decay of the vectors of eigenvalues $\bm{\lambda}^u$ and $\bm{\lambda}^w$. 

Once the POD functional spaces are set, the reduced quantities fields can be approximated with:
\begin{equation}\label{eq:aprox_fields}
{u^r} \approx \sum_{i=1}^{N_u^r} a_i(t,\mu) {\varphi_i}(\bm{x}), \mbox{\hspace{0.5 cm}}
w^r\approx \sum_{i=1}^{N_w^r} b_i(t,\mu)\chi_i{\blue{(\bm{x}),}}
\end{equation}
{\blue{w}}here the coefficients $a_i$ and $b_i$ {\it{depend only}} on the time and parameter spaces and the basis functions $\bm{\varphi}_i$ and $\bm{\chi}_i$ {\it{depend only}} on the physical space and not on the parametrized geometry. 

By denoting ${\bm{{\mathcal{B}}}}=\begin{bmatrix} {\bm{{\mathcal{B}}}} _u & \bm{0} \\ \bm{0} & {\bm{{\mathcal{B}}}}_w \end{bmatrix}$ and ${\bm{{\mathcal{B}}}}^T =\begin{bmatrix} {{\bm{{\mathcal{B}}}}}^T_u & {\bm{0}} \\ {\bm{0}} & {\bm{{\mathcal{B}}}}^T_w \end{bmatrix}$, the {\magenta{unknown 
coefficients}} $V=\begin{bmatrix} \bm{a} \\ \bm{b} \end{bmatrix}$ then can be obtained through a Galerkin projection of the {\magenta{full order}} system of equations onto the POD reduced basis spaces with the solution of a consequent reduced iterative algebraic system of equations for the increment $\delta V_h^{n+1}(\mu) = V_h^{n+1}(\mu) -V_h^{n}(\mu)$,
\begin{equation}\label{eq:system_linear}
{{\bm{{\mathcal{B}}}}}^T\orange{\begin{bmatrix} {{\bf{A}}} + 
  \tau_n{\bf{L}_1}  & \tau_n{\bf{L}_2}  
  \\ \tau_n{\bf{C_{\mathcal{F}_G}}}
  & 0  
\end{bmatrix}}
{\bm{{\mathcal{B}}}}\delta {V_h^{n+1}}(\mu) = - \tau_n {\bm{{\mathcal{B}}}}^T  R({\bm{{\mathcal{B}}}}V_h^{n}(\mu)),
\end{equation}
which leads to the following algebraic reduced system:
\begin{equation}\label{eq:system_linear_reduced}
 \orange{\begin{bmatrix} {{\bf{A}}} + 
  \tau_n{\bf{L}_1}  & \tau_n{\bf{L}_2}  
  \\ \tau_n{\bf{C_{\mathcal{F}_G}}}
  & 0  
\end{bmatrix}^r}
\delta {V_{h,r}^{n+1}}(\mu) = - \tau_n R^r(V_{h,r}^{n}(\mu)).
\end{equation}
{\bblue{We underline, that in the aforementioned reduced level stage
we employ once more the IMEX approach which treats the nonlinear term explicitly, allowing it to act as a forcing term on the right-hand side, avoiding stability issues caused by the nonlinearity.
}}
{\magenta{We also clarify that during the POD procedure we have to assemble all matrices including the one related to the nonlinear term. Then we project them in the reduced basis space and we solve the reduced system. With this method we gain execution time during the latter stage, and not in the assembling stage as a DEIM/EIM would allow.}}
%
%
\begin{remark}
The initial conditions for the ROM system of equation \eqref{eq:system_linear_reduced} are obtained performing a Galerkin projection of the initial {\magenta{full order}} condition {\magenta{${u}(\cdot,0;\mu)$ 
onto}} the POD basis spaces.
For an efficient ROM basis construction, we follow some ideas of the authors as demonstrated in \cite{KaBaRO18}. In particular, concerning the full order method snapshots extension, and the extension of the solution to the surrogate domain into the ghost area,  we use the solution values as they have been computed using the cut finite element method smooth mapping from the true to the unfitted mesh domain. {\bblue{The stabilization term $\sum_{F\in \mathcal{F}_G}\left(\alpha_1 h^{-2}\ldb u_h\rdb,\ldb v_h\rdb\right)_F$, which depends on $\alpha_1$, extends the coercivity
from the physical domain $\Omega(\mu)$ to the extended mesh domain $\Omega_\mathcal{T}(\mu)$}}. This allows a smooth extension of the boundary solution to the neighbouring ghost elements with values which are decreasing smoothly to zero, see for instance the zoomed image in \autoref{fig:poisson_zoom}.
 This approach provides a regular solution in the background domain and permits, therefore, the construction {\magenta{of a 
 reduced basis with better approximation properties. {\magenta{For the instability issues related to the investigation, visualization, as well as, relative error reporting and comparison between several stabilization methods we refer to the detailed numerical investigation in the work  \cite{KaBaRO18}}} where we test several approaches, namely snapshots' zero extension, natural smooth extension, harmonic extension and we examine how we can avoid instabilities and achieve a good and computationally cheap reduced basis.}}
\end{remark}
\section{Numerical experiments}\label{sec:num_exp}
In the present section, we test the presented methodology considering numerical experiments for the evolutionary Cahn-Hilliard system while natural homogeneous Neumann boundary conditions, and/or zero Dirichlet ones are present.
We start by testing the robustness of the full order model (FOM) {\magenta{for two settings of a classical benchmark test case}} and we continue with numerical examples in which geometrical parametrization for the embedded domain is considered.
The background domain in all experiments is the rectangle $[-0.5, 0.5] \times [-0.5, 0.5]$, and {\magenta{the parameters in the double well free energy}} $\gamma_0 $,  $\gamma_1$ and $\gamma_2$ {\magenta{in}} formula (\ref{double-well}) {\magenta{for the ROM experiments}} are taken as $\gamma_0 = 2 $,  $\gamma_1 = 9$ and $\gamma_2 = 4$, with 
{\magenta{$\varepsilon=10^{-2}$  proportional to the two fluids interface size}}.
The results for all test problems have been obtained with  mesh size $h = 1/48$ 
{\magenta{and simulation time points inside the interval $[0, 100\tau_n]$ and for the time step size $1.5625\cdot10^{-6}$}}
unless otherwise stated, Nitsche parameter $a_N = 10$ and jump stabilization parameter $a_1 = 0.1^3$.
%
%
{\magenta{For the numerical experiments, we carried out using the software package ngsxfem extension  of ngsolve \cite{ngsolve,ngsxfem,schoeberl} for the full order solution part, and the software RBniCS \cite{rbnics} for the reduced order solution one. 
 We have used a device with an Intel\textsuperscript{\textregistered} Core\textsuperscript{TM}  i7-4770HQ 3.70GHz CPU. 
 %
%
%
%
%
%
%
%
}}
\subsection{Robustness of the FOM solver}
In this first numerical example, we test the validity of the full order solver and the IMEX method. {\magenta{Two settings of a classical 
benchmark test are examined, 
 for a cross-shaped initial condition 
phase-field interface: 
i) The classical setting $\gamma_2=1$, $\gamma_1=0$, $\gamma_0=-1$, adding the unit value, \cite{Cha16,CrossXU2019524}, and 
ii) the suggested in the present work setting $\gamma_0 = 2 $,  $\gamma_1 = 9$ and $\gamma_2 = 4$
}}
applying {\magenta{the}} time step {\magenta{size}} $\tau_n={\mathcal{O}(10^{-8})}$,  for $10000$ time steps respectively and without any embedded geometry, see Figure \ref{background_mesh_Clean} (i).
 \begin{figure} \centering
   (i) \hskip-3pt \includegraphics[width=0.305\textwidth]{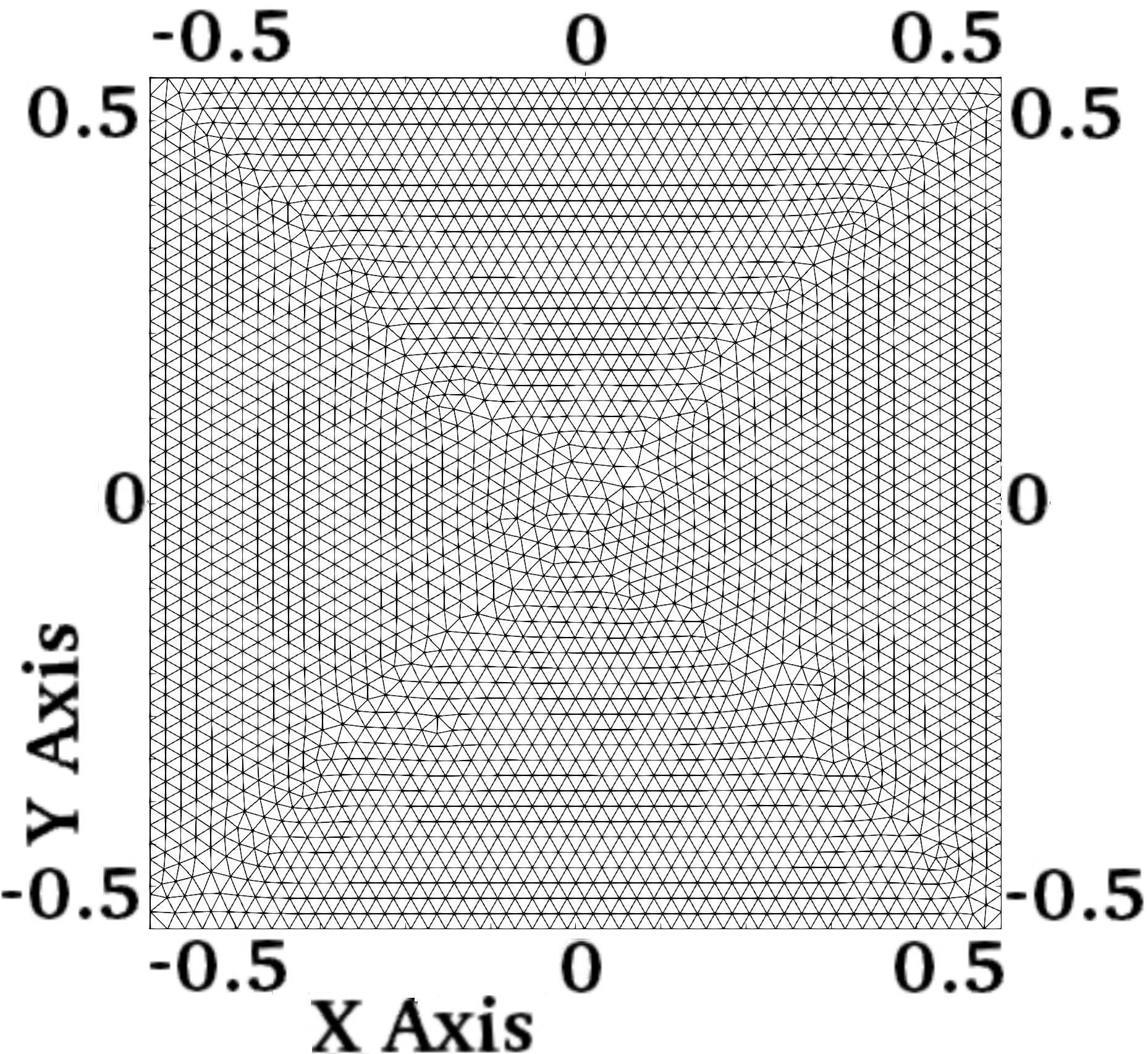}
   \qquad
    \includegraphics[width=0.33\textwidth]{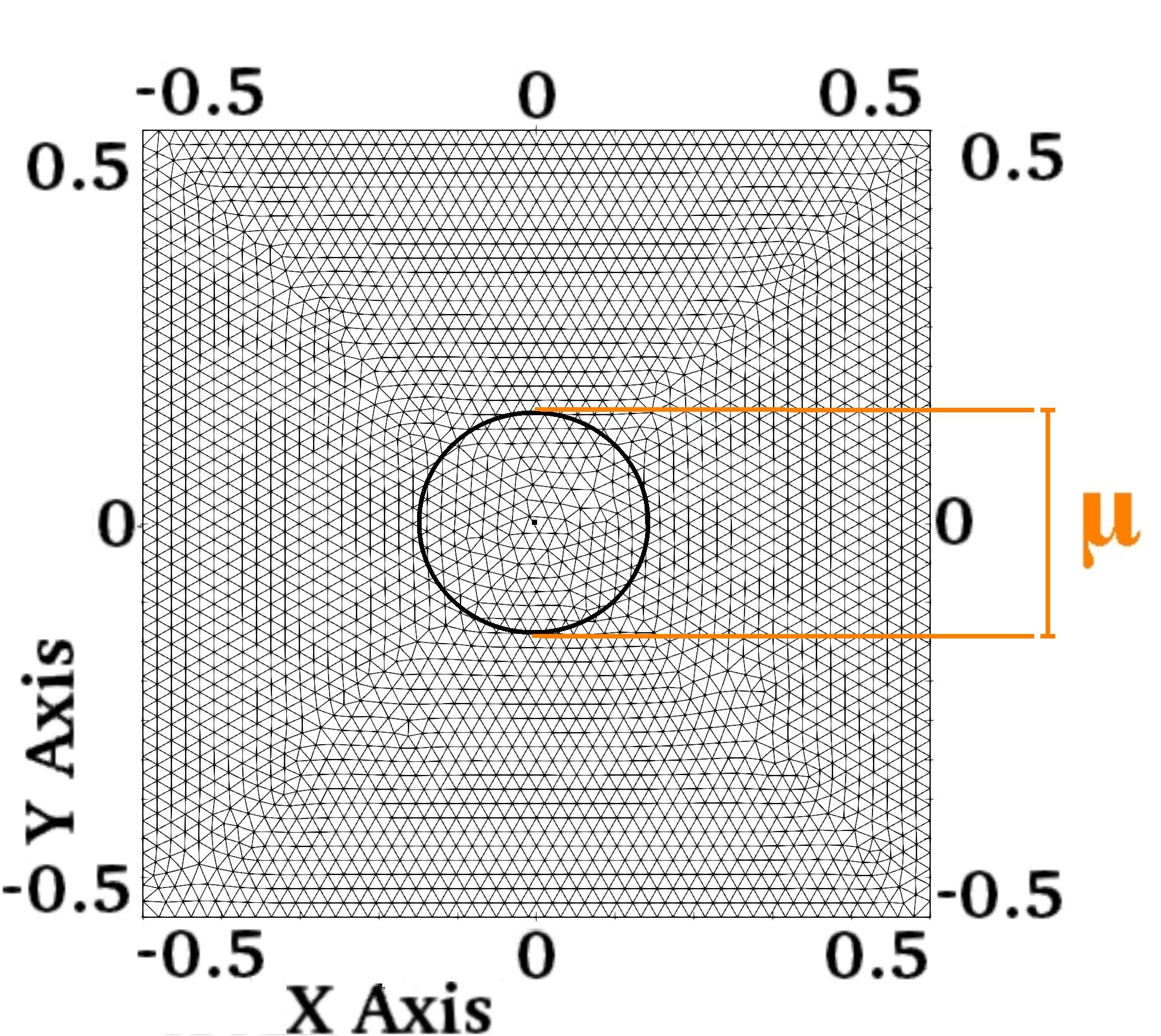} \hskip-10pt(ii)
\caption{The background mesh (i), and a sketch of the embedded domain and the parameters considered in the numerical examples (ii).}
  \label{background_mesh}\label{background_mesh_Clean}
\end{figure}
\begin{figure}
\centering
\includegraphics[width=0.33\textwidth]{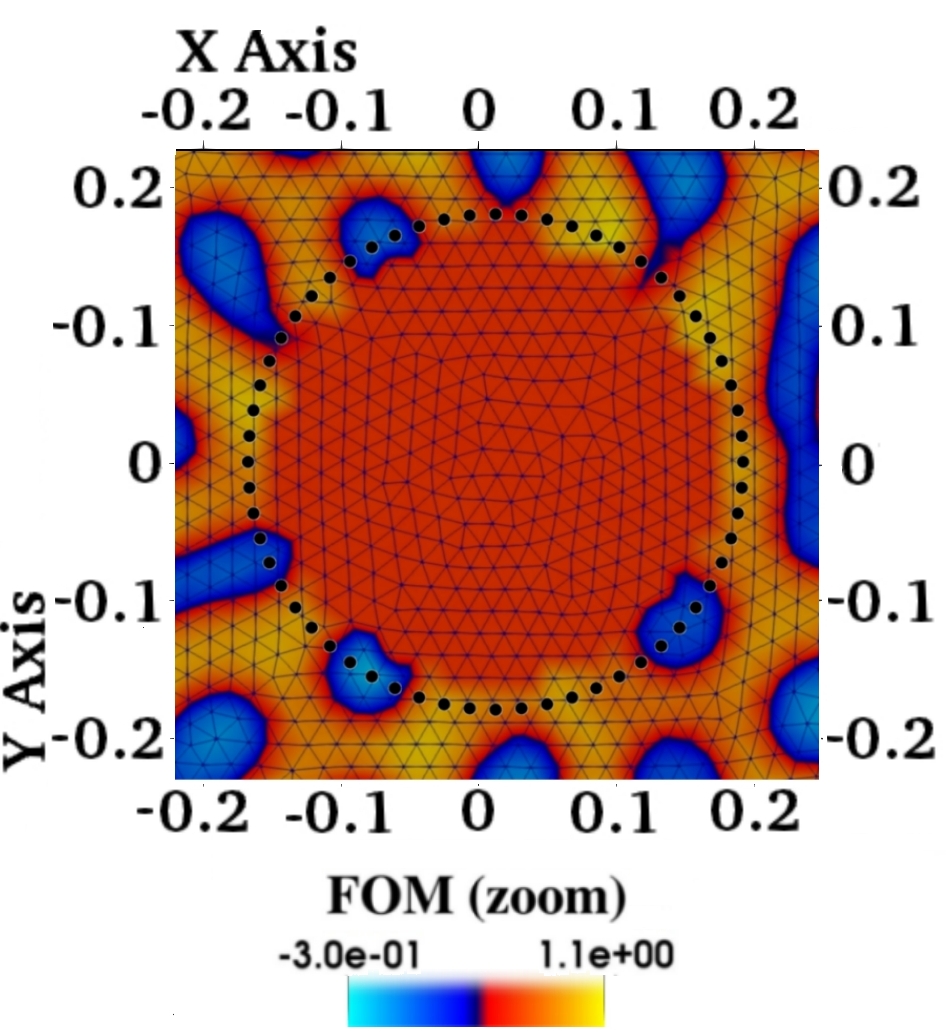}
\caption{A zoom onto the embedded circular domain
of the numerical example shows the smoothing natural smooth extension procedure employed by the cut finite element method  inside the ghost area.}
\label{fig:poisson_zoom}
\end{figure}
{\magenta{We specify that comparing the ii) choice with the traditional in the literature $u^3-u$ nonlinearity i) case as in \cite{Cha16,CrossXU2019524} numerical tests, 
shown much faster dynamics at the beginning of the evolution of the physics of the problem for the ii) setting. Due to the latter, and after experimental verification, we needed to choose the aforementioned order of time step size to achieve a good approximation of the conservation of mass property even if we have employed the IMEX approach.}}
For the first experiment, we consider  the initial state

\begin{eqnarray}
   u_0(x, y) = \left\{\begin{array}{lr}
        {\magenta{0.95}},  & \text{if } 5|(y -0.5) - \frac{2}{5}(x - 0.5)| + |\frac{2}{5}(x - 0.5) - (y - 0.5)|\leq 1,\\
    {\magenta{0.95}},  &  \text{if } 5|(x -0.5) - \frac{2}{5}(y - 0.5)| + |\frac{2}{5}(y - 0.5) - (x - 0.5)|\leq 1,\\
        {\magenta{-0.95}}, &
        \text{otherwise},
        \end{array}\right.
\end{eqnarray}
setting as in \cite{Cha16,CrossXU2019524}, 
while for the second
{\magenta{we have substitute the initial values $-95$, $95$ to $0$, $0.6$ respectively.
}}
%
In both cases, with sharp corners in a cross-shaped interface case to be present, we observe {\magenta{the}} evolution {\magenta{toward a}} circular interface, see Figure \ref{CrossBench}. Moreover mass conserves 
with error of order $10^{-15}$ with respect to time evolution, {\magenta{see Figure \ref{fig:coserv_mass_crossBench}.  
For}} the sake of shortness we have visualized only the ii) experiment while for i) we noticed 
similar results as well as the conservation of mass. 
\begin{figure} \centering
\begin{minipage}{\textwidth}
\centering
\makebox[\textwidth][c]{%
\hskip4pt \begin{minipage}{0.225\textwidth}
  \includegraphics[width=\textwidth]{{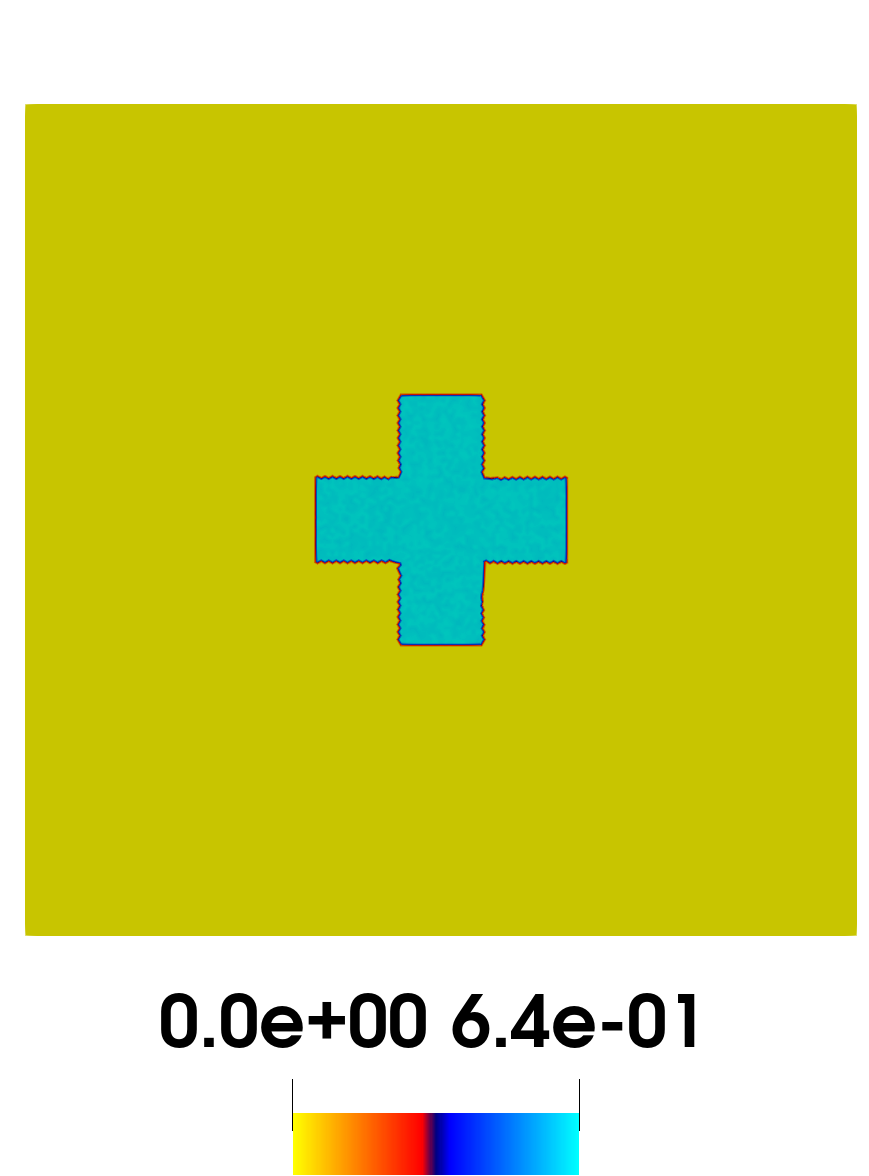}} 
\end{minipage}
 \qquad \qquad
 \begin{minipage}{0.225\textwidth}
  \includegraphics[width=\textwidth]{{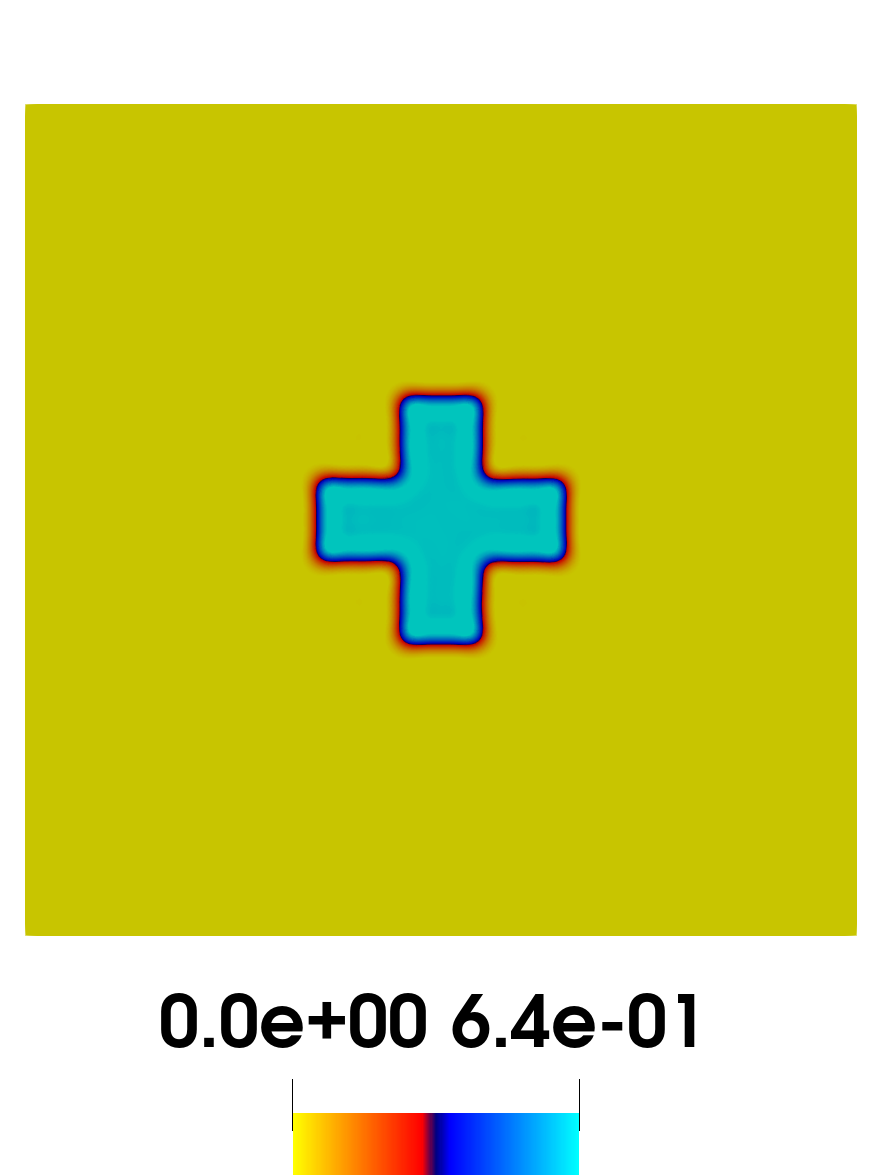}}
\end{minipage}
 \qquad \qquad
 \begin{minipage}{0.225\textwidth}
  \includegraphics[width=\textwidth]{{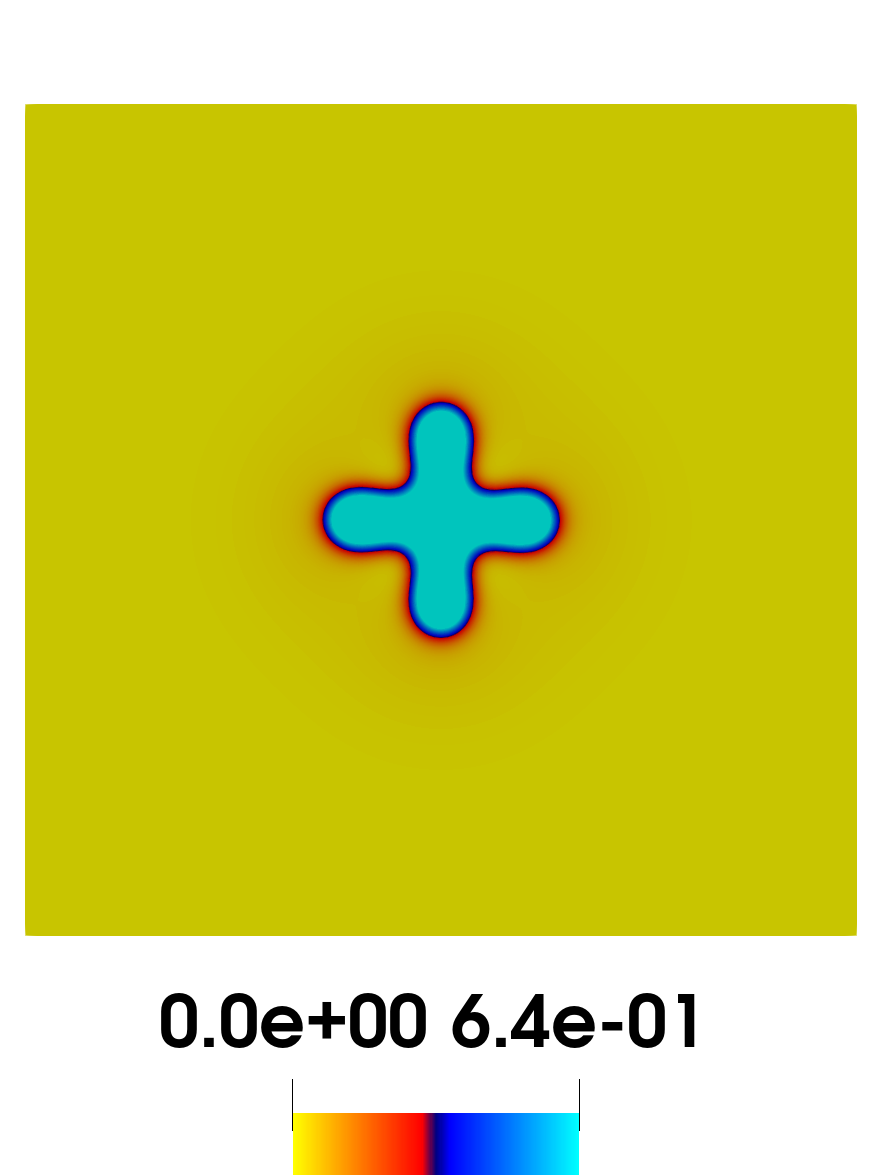}}
\end{minipage}
}
\newline
\makebox[\textwidth][c]{
\begin{minipage}{0.225\textwidth}
  \includegraphics[width=\textwidth]{{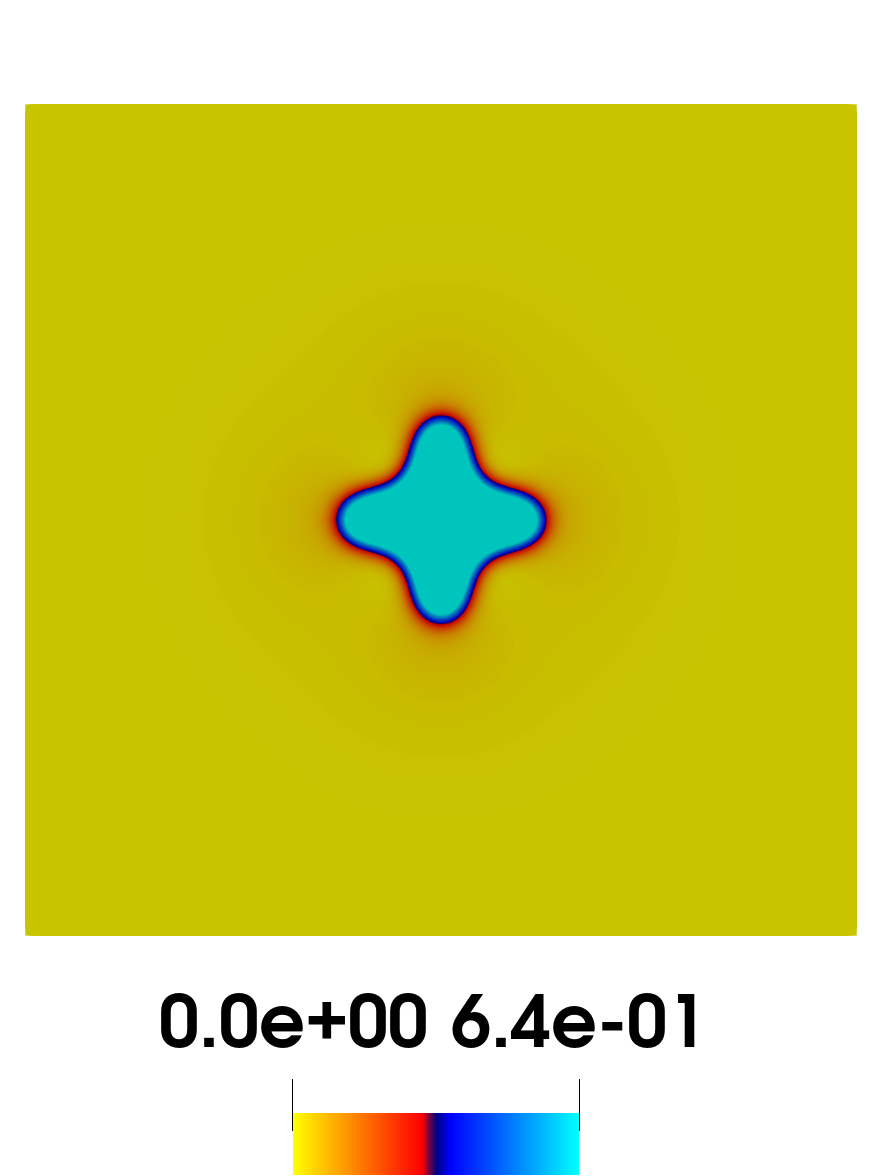}}
\end{minipage}
 \qquad \qquad
\begin{minipage}{0.225\textwidth}
  \includegraphics[width=\textwidth]{{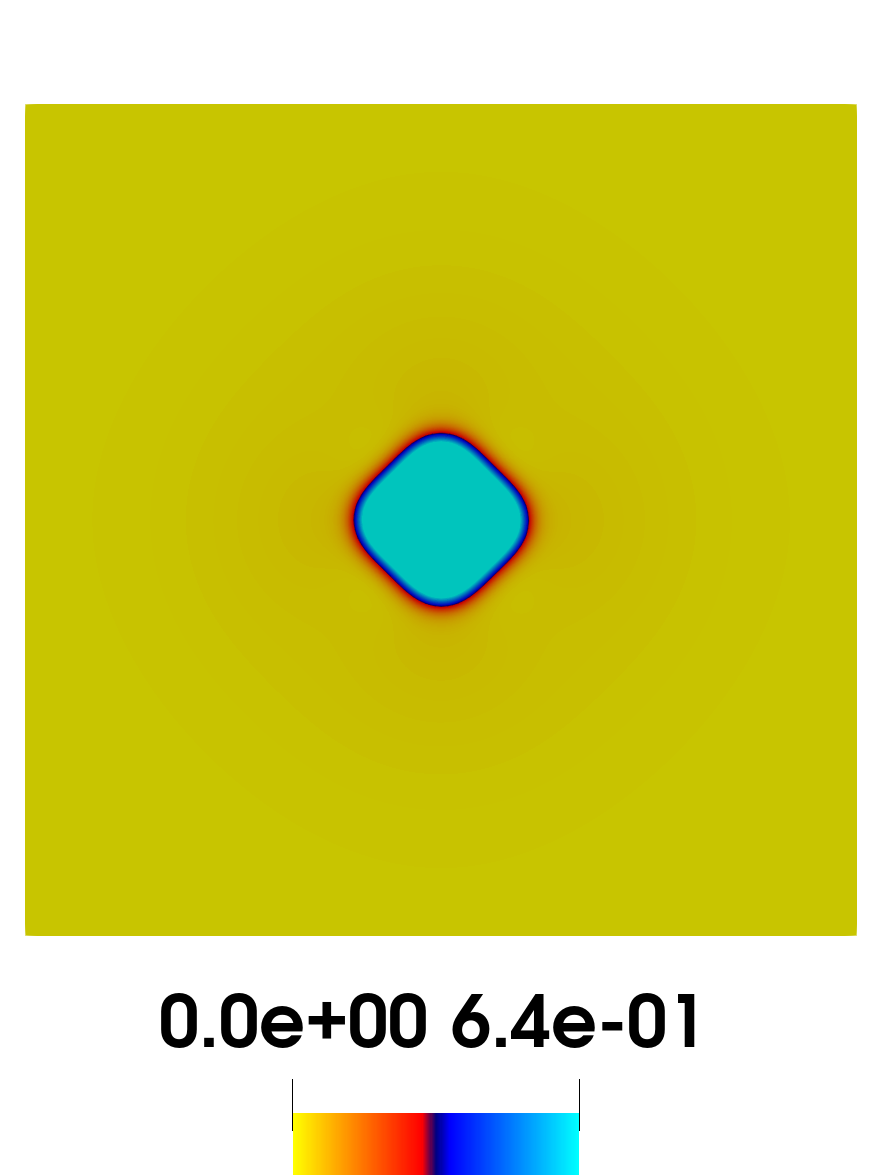}}
\end{minipage}
 \qquad \qquad
 \begin{minipage}{0.225\textwidth}
  \includegraphics[width=\textwidth]{{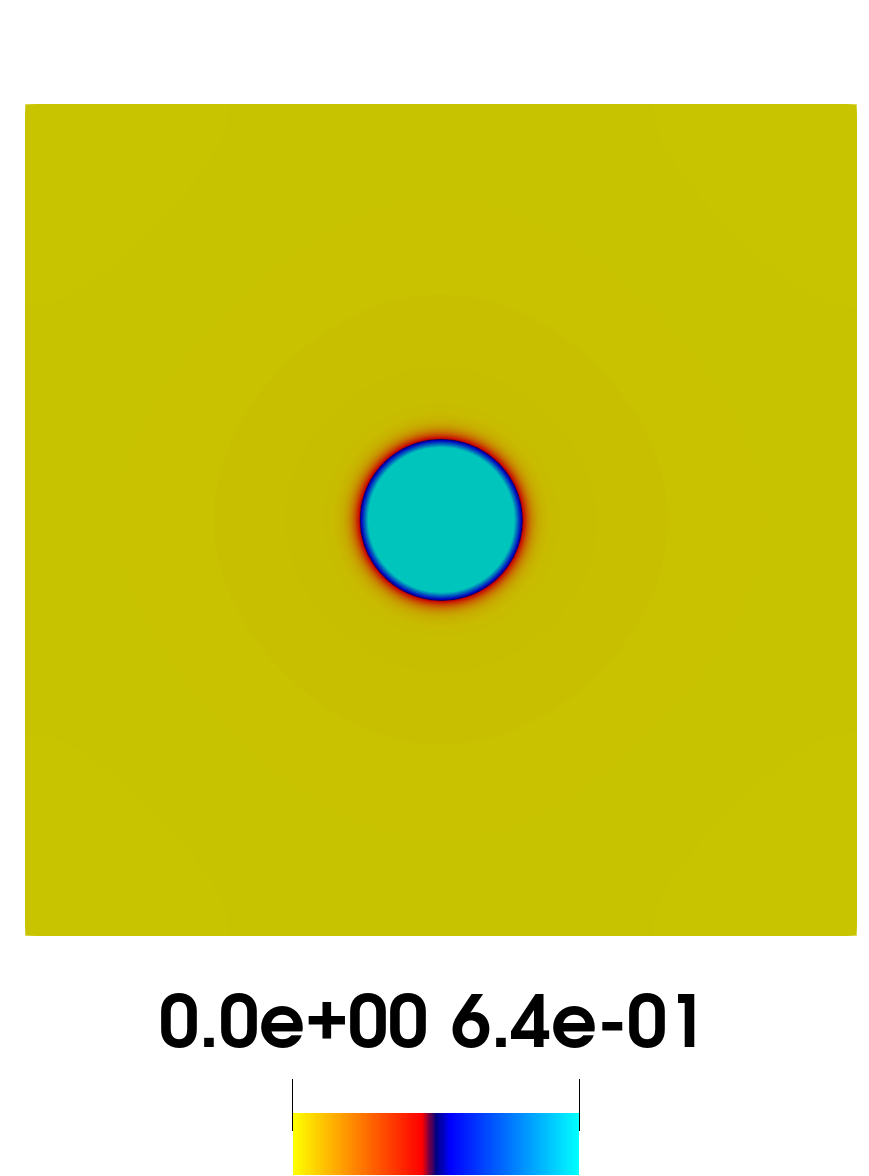}}
\end{minipage}
}
\end{minipage}
  \caption{Testing the IMEX method: Cross type initial data and the phase field results for times $t=[0,10,1300,3000,5000,10000]\tau_n$ which visualizes that it {\magenta{evolves 
  to}} a circular interface.}
  \label{CrossBench}
\end{figure}
\begin{figure} \centering
\begin{minipage}{0.66\textwidth}
  \includegraphics[width=\textwidth]{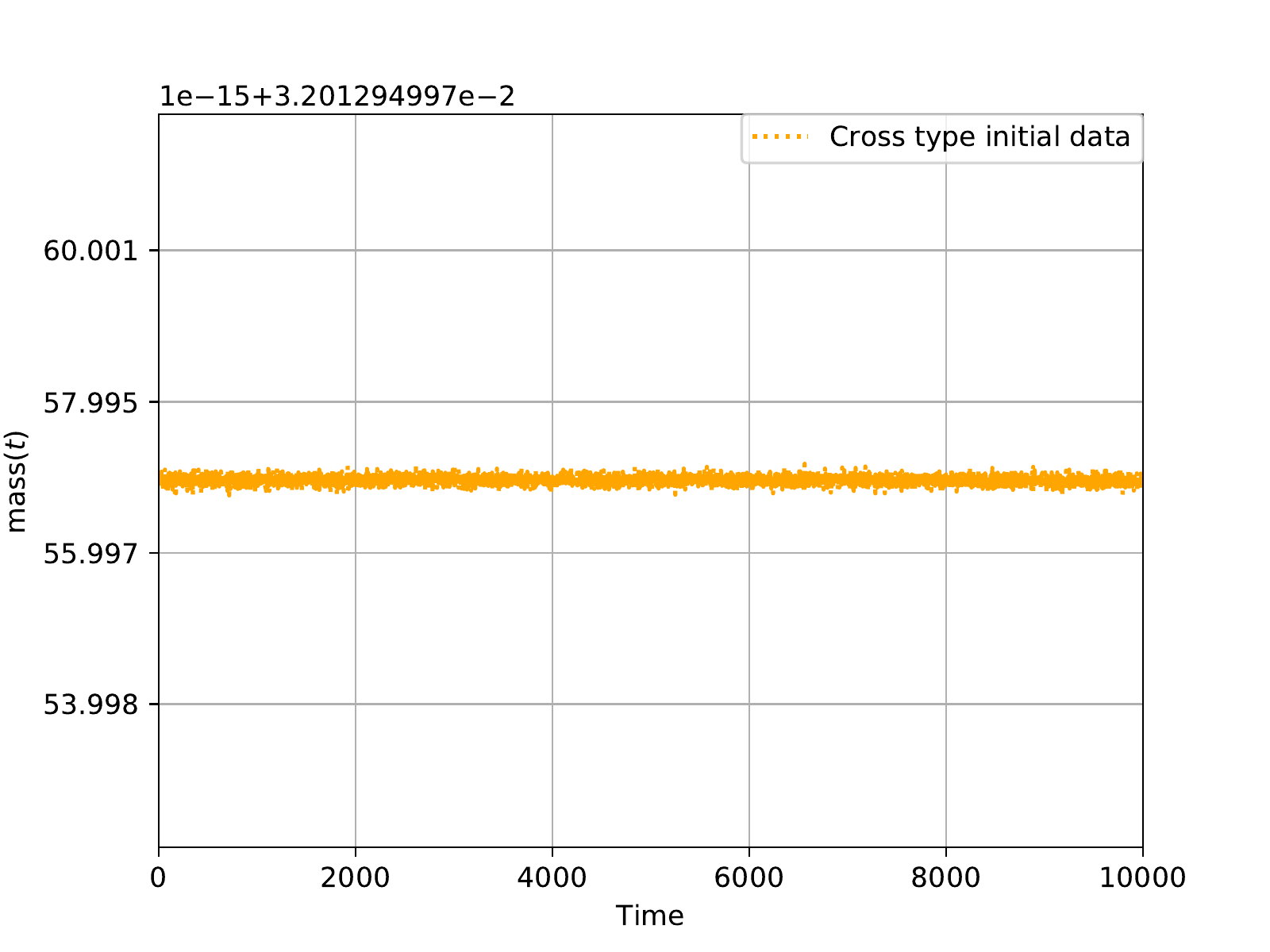}
\end{minipage}
  \caption{Testing the IMEX method: Cross type initial data and the mass conservation as time evolves for times $t=[0,...,10000]\tau_n$
  .}
  \label{fig:coserv_mass_crossBench}
\end{figure}
%
%
\begin{remark}
{\magenta{We also remark that for the FOM solver and as a second step after the initial cross-type experiment, we tested numerically an example considering a long time interval with Neumann and Dirichlet embedded geometries. {\Xgreen{This has shown a}}  stable behavior as time evolves, while the conservation of mass is fulfilled, emphasizing the good results at the final times, e.g. $500\tau_n$. For this experiment we used space-time discretization mesh sizes $h=1/96$ and $\tau_n= 
\mathcal{O}{(10^{-7})}$, see e.g. Figure \ref{FULL_3D}.}}
\end{remark}
\begin{figure}
\centering
\begin{minipage}{\textwidth}
\centering
\begin{minipage}{0.225\textwidth}
  \includegraphics[width=\textwidth]{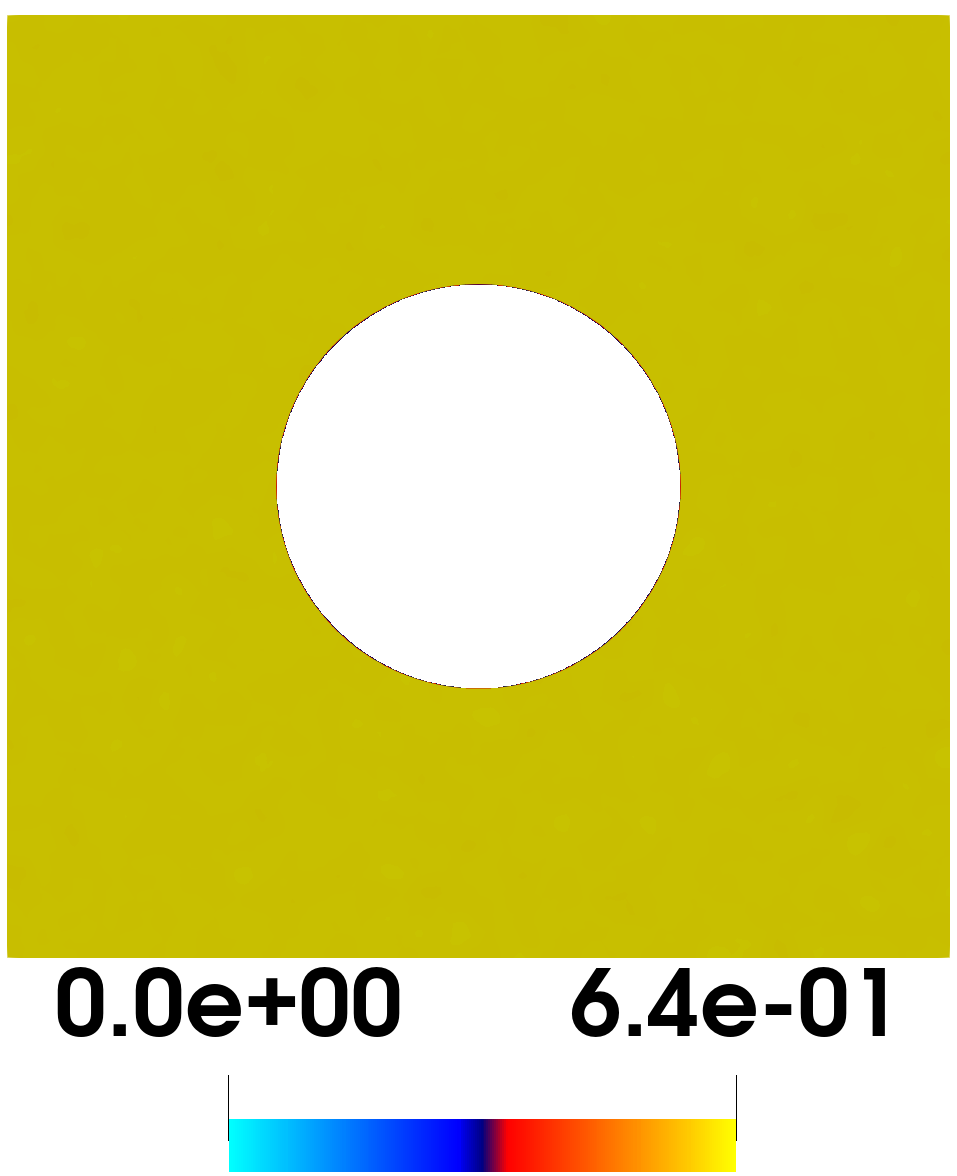}
\end{minipage}
\hskip18pt
\begin{minipage}{0.225\textwidth}
  \includegraphics[width=\textwidth]{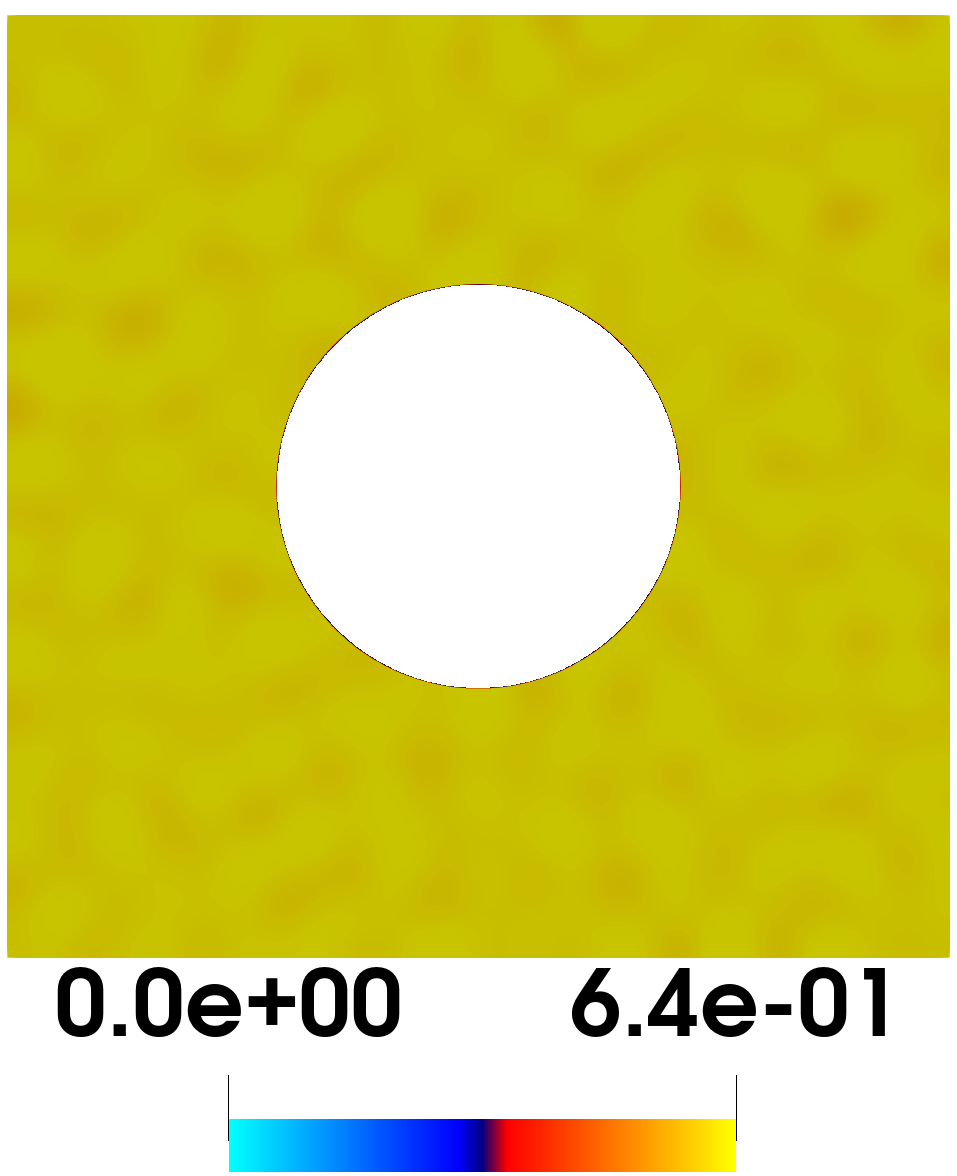}
\end{minipage}
\hskip18pt
\begin{minipage}{0.225\textwidth}
  \includegraphics[width=\textwidth]{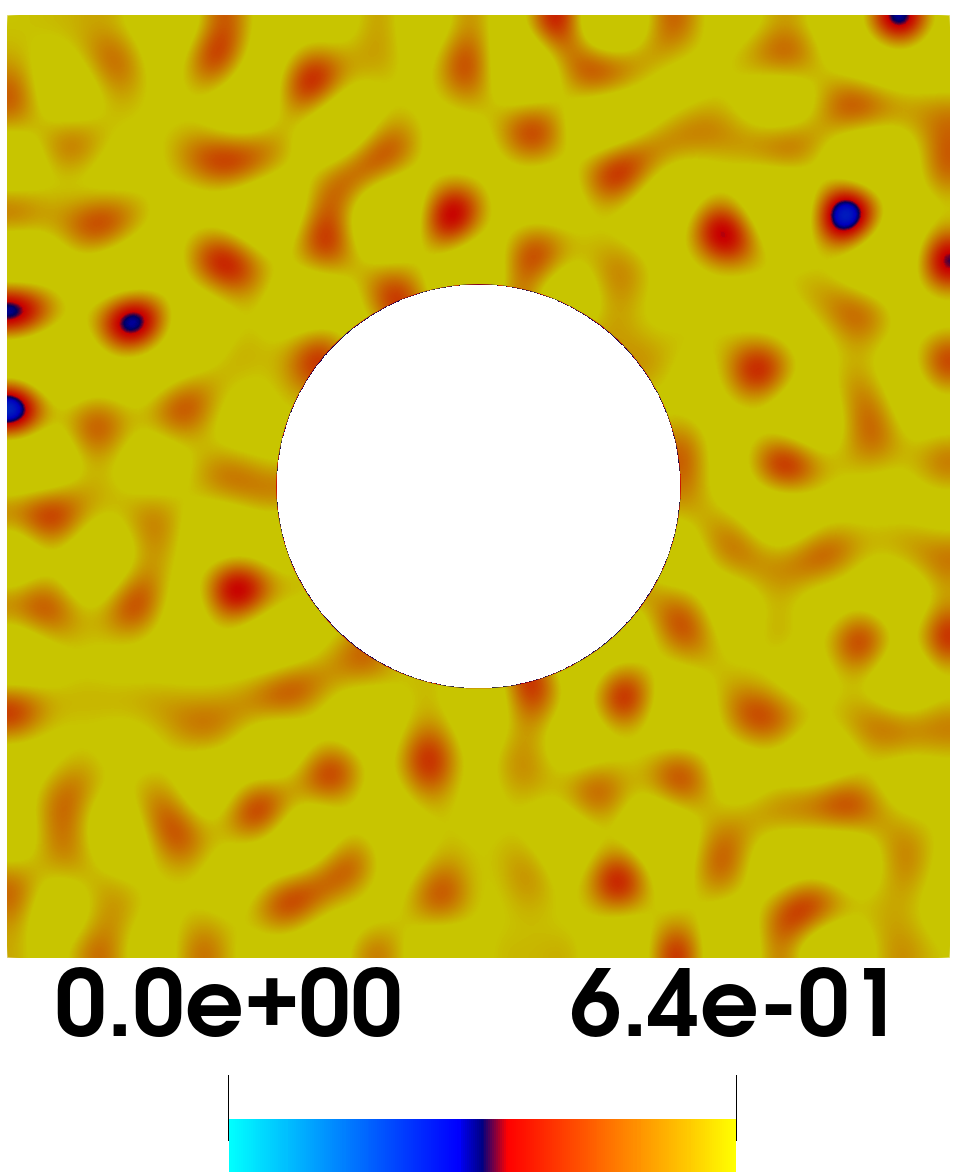}
\end{minipage}
\vskip20pt
\begin{minipage}{0.225\textwidth}
  \includegraphics[width=\textwidth]{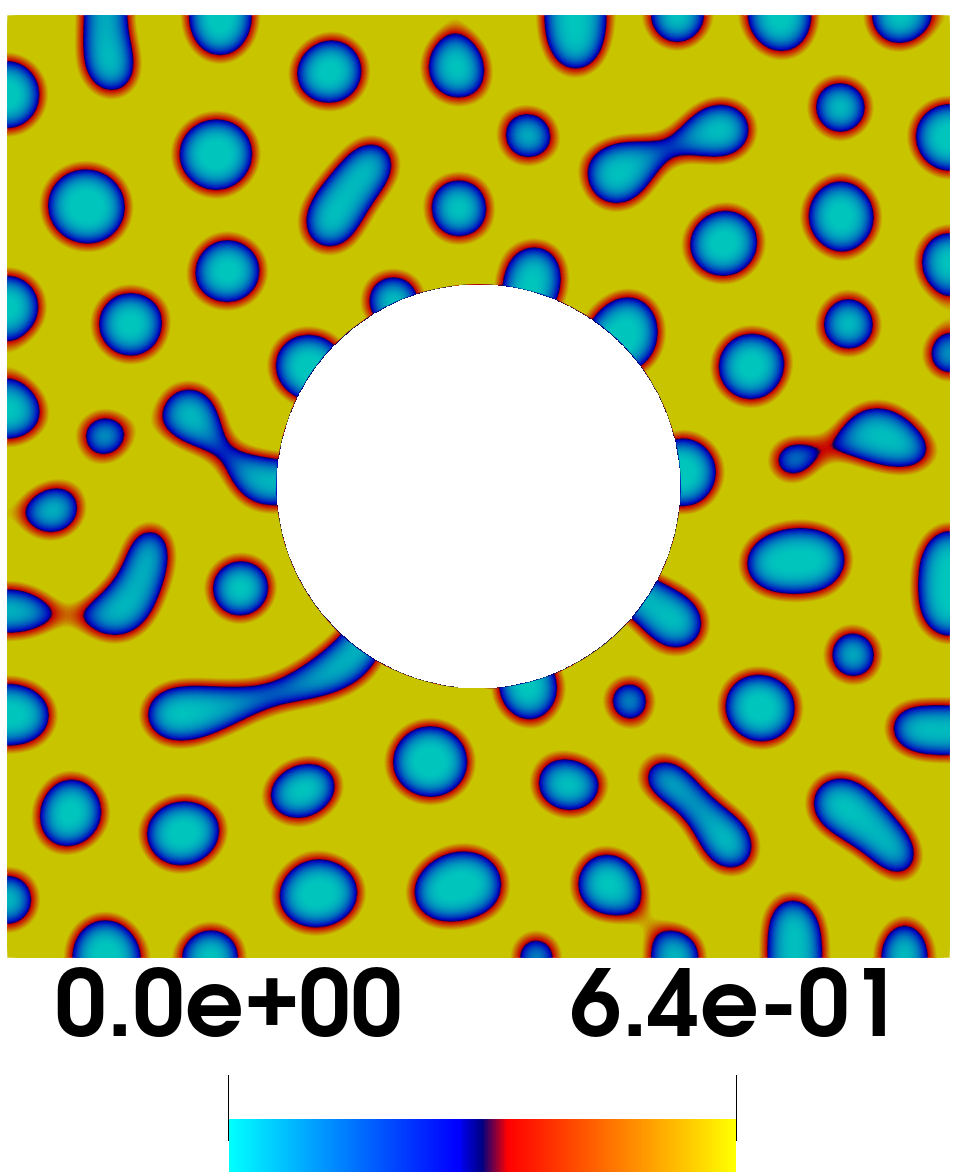}
\end{minipage}
\hskip18pt
\begin{minipage}{0.225\textwidth}
  \includegraphics[width=\textwidth]{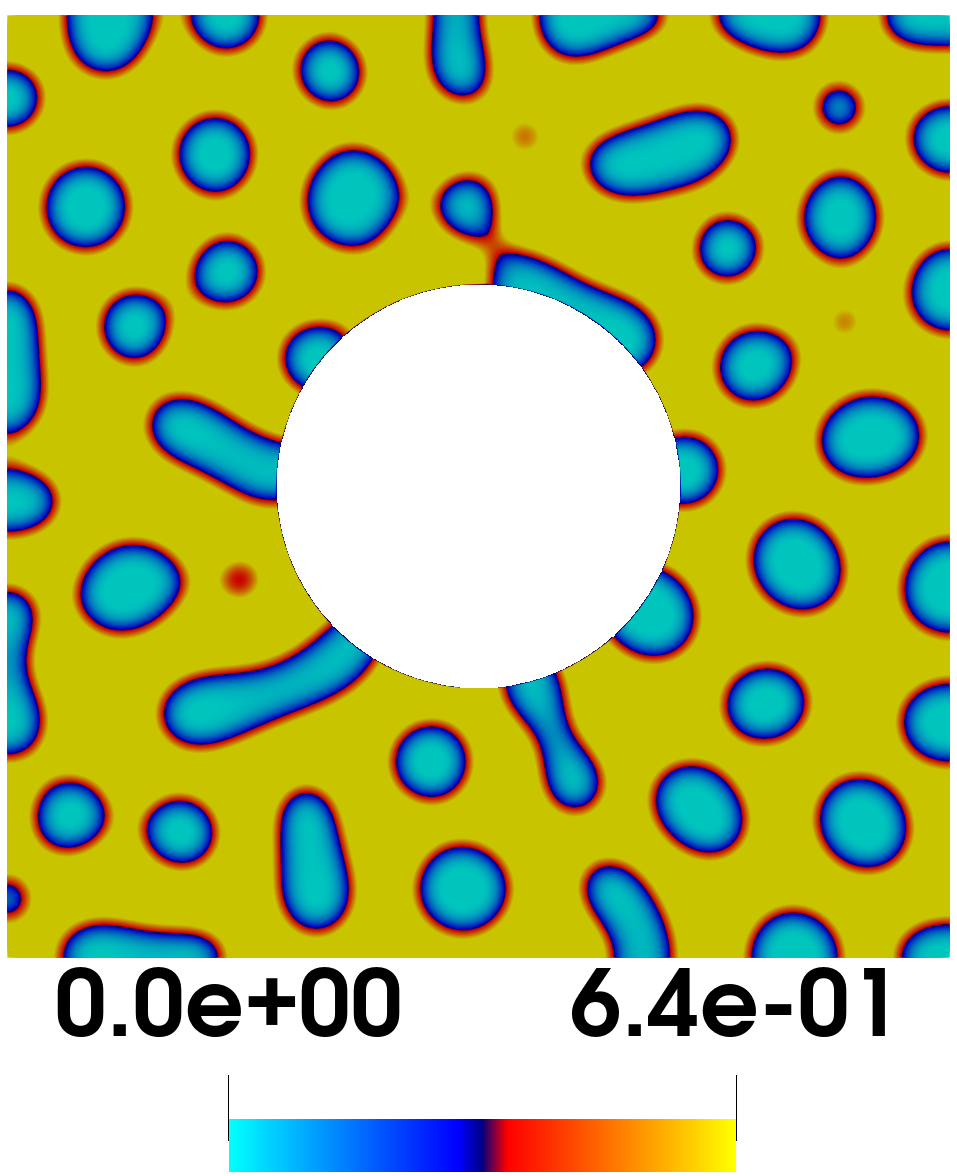}
\end{minipage}
\hskip18pt
\begin{minipage}{0.225\textwidth}
  \includegraphics[width=\textwidth]{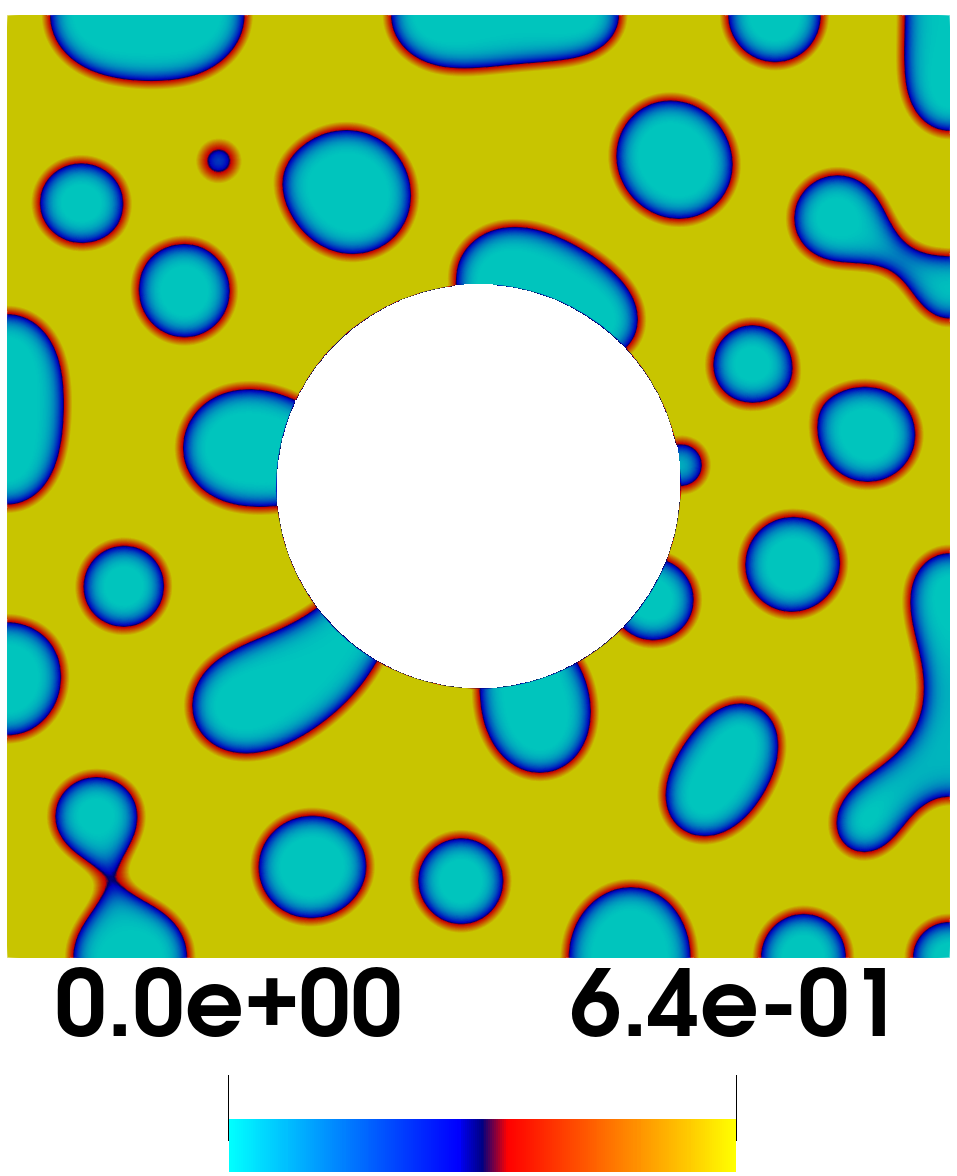}
\end{minipage}
\end{minipage}
\vskip20pt
\begin{minipage}{\textwidth}
\centering
\begin{minipage}{0.225\textwidth}
  \includegraphics[width=\textwidth]{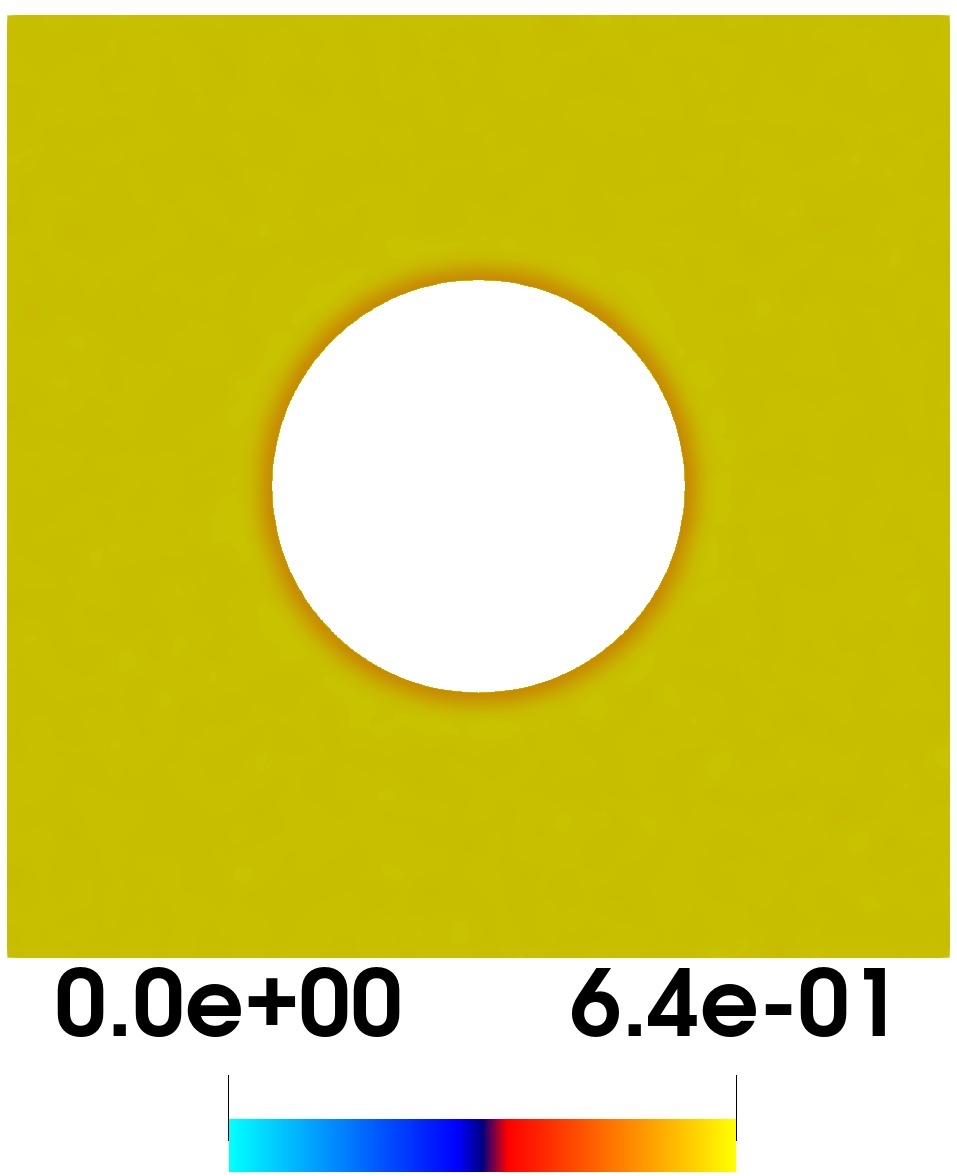}
\end{minipage}
\hskip18pt
\begin{minipage}{0.225\textwidth}
  \includegraphics[width=\textwidth]{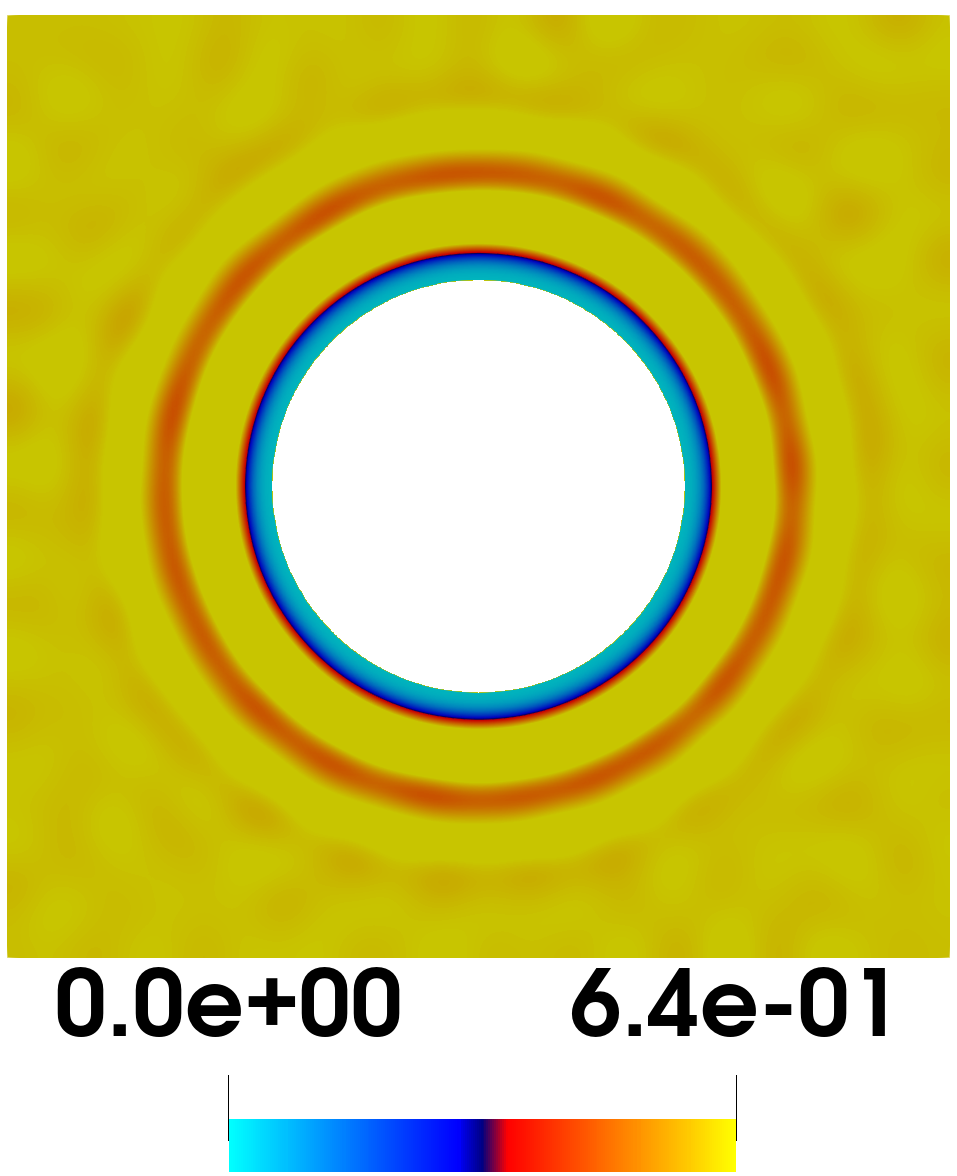}
\end{minipage}
\hskip18pt
\begin{minipage}{0.225\textwidth}
  \includegraphics[width=\textwidth]{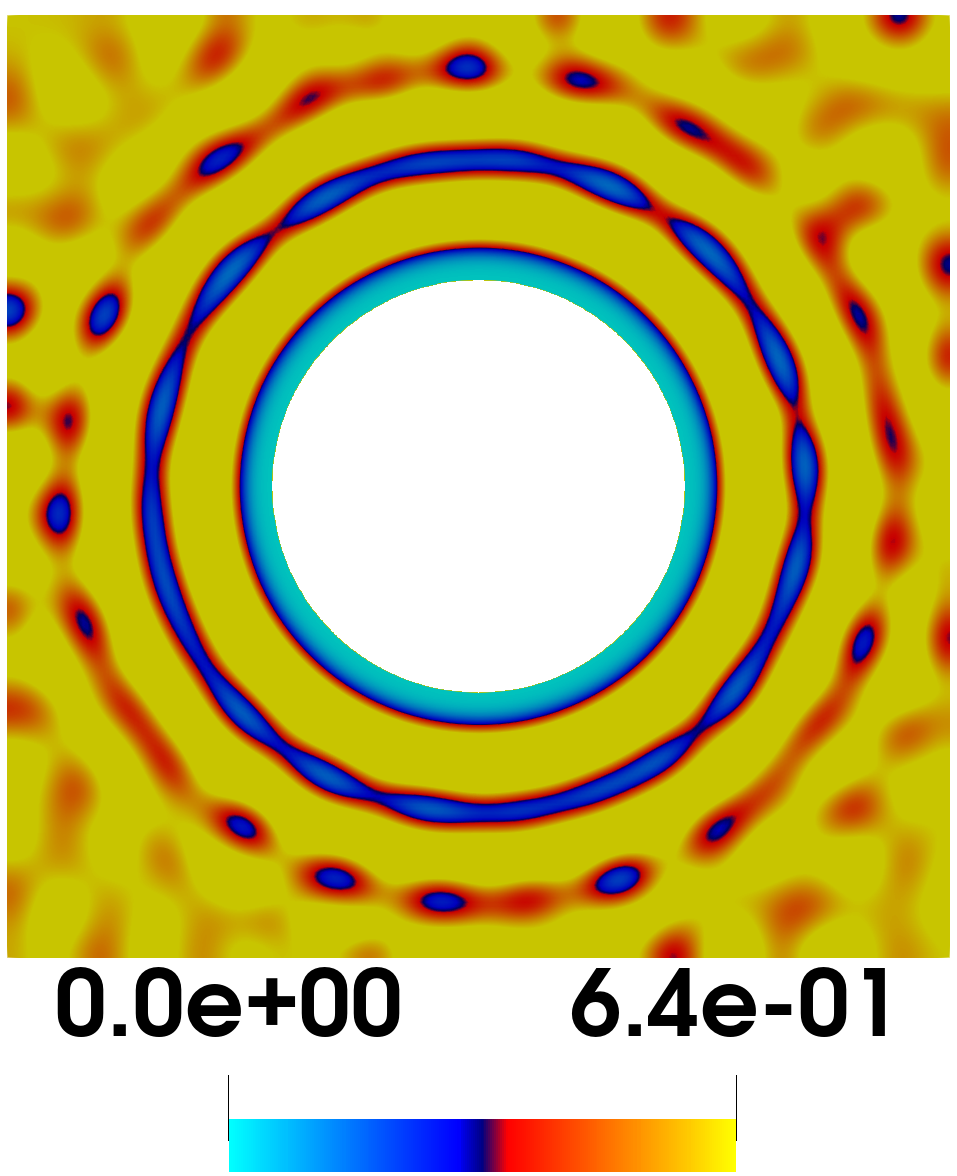}
\end{minipage}
\vskip20pt
\begin{minipage}{0.225\textwidth}
  \includegraphics[width=\textwidth]{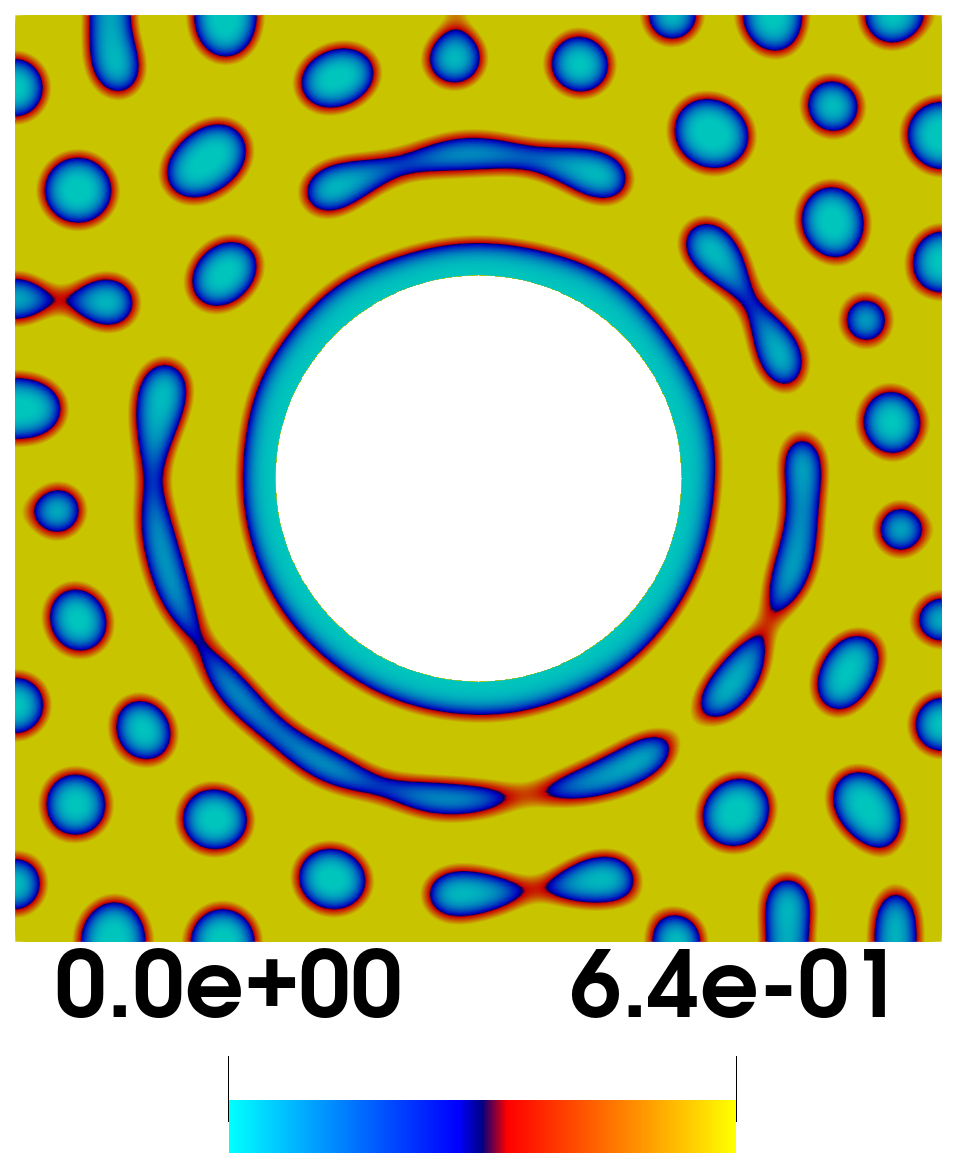}
\end{minipage}
\hskip18pt
\begin{minipage}{0.225\textwidth}
  \includegraphics[width=\textwidth]{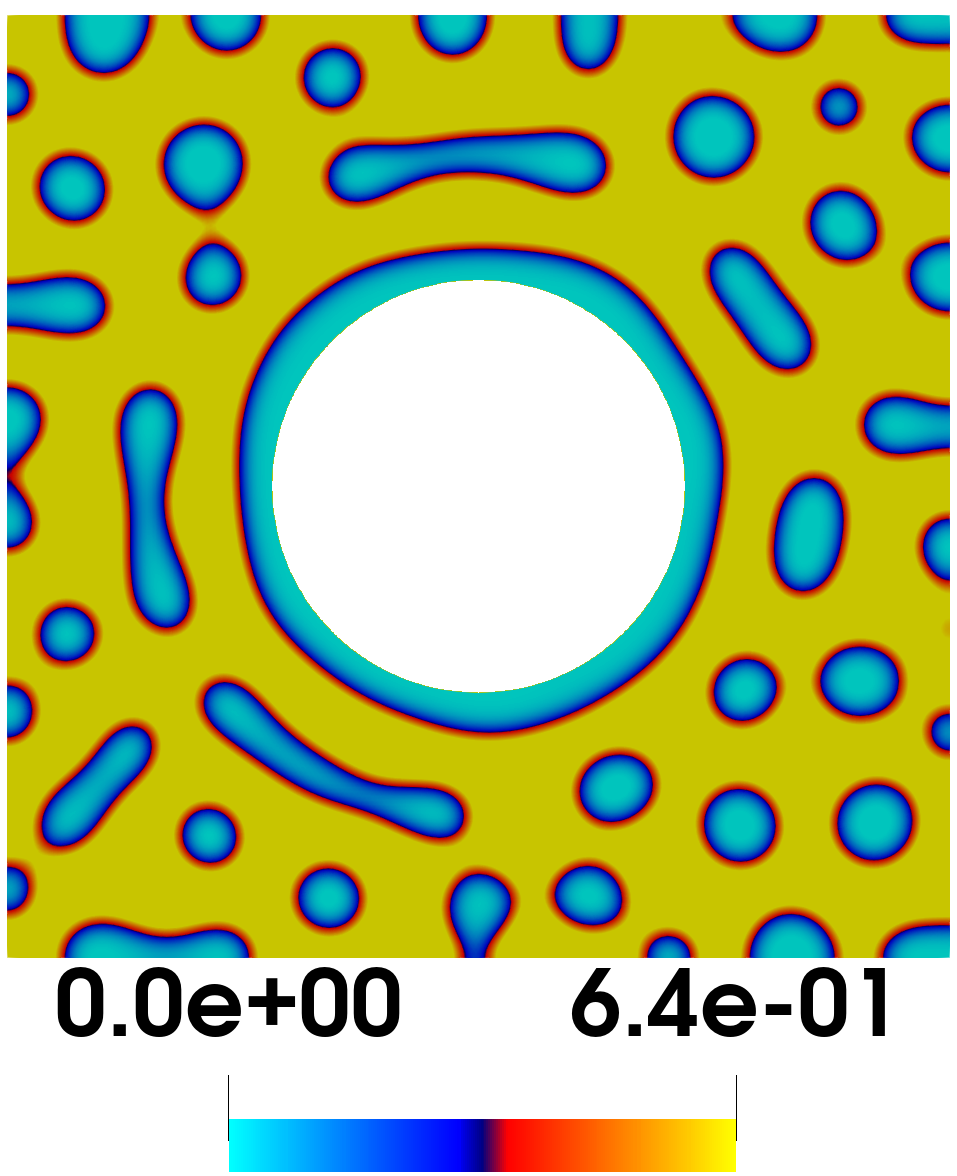}
\end{minipage}
\hskip18pt
\begin{minipage}{0.225\textwidth}
  \includegraphics[width=\textwidth]{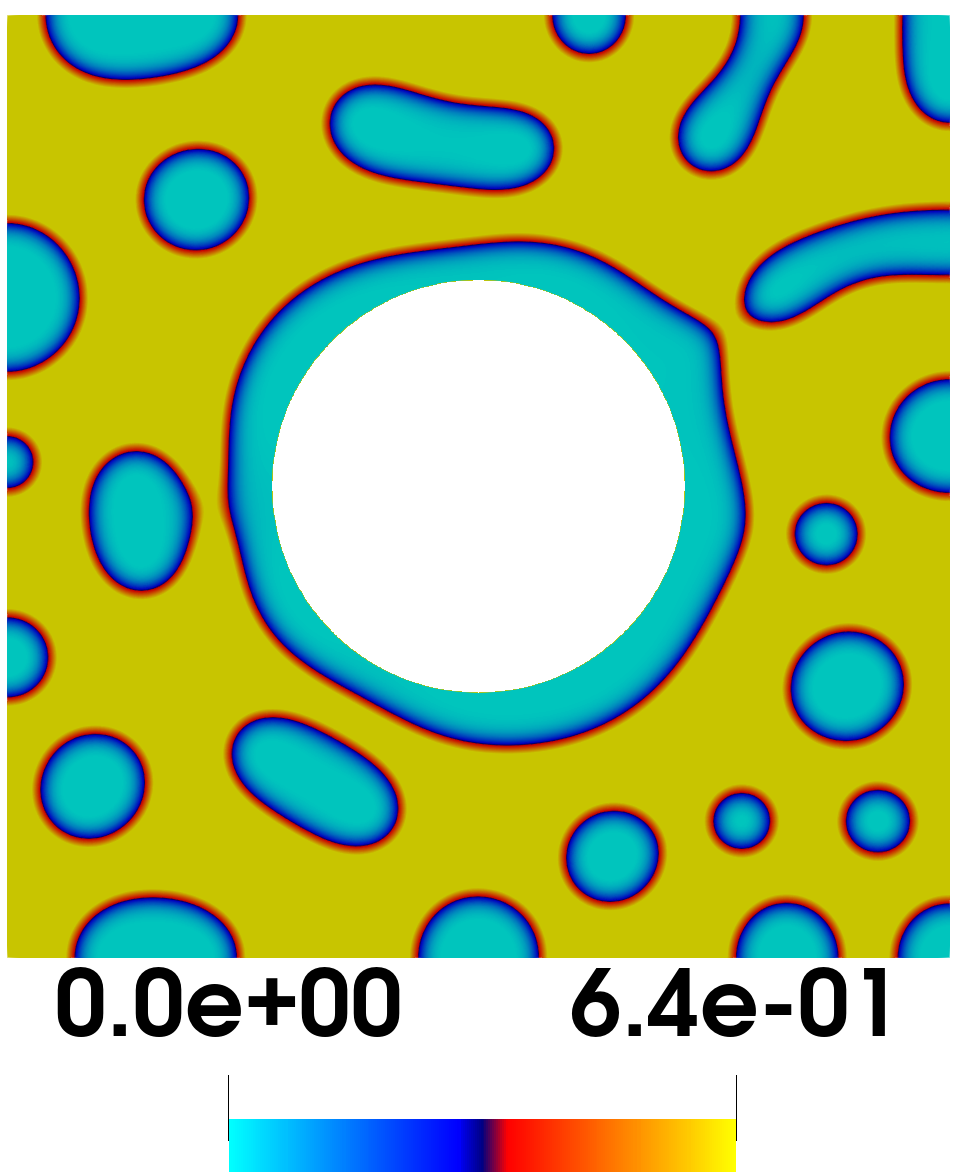}
\end{minipage}
\end{minipage}
\caption{{{The full order solution of the concentration field}} 
(whole background geometry) for a fixed parameter $\mu_{\text{test}} = 0.436$ at times $t=[1, 25, 50, 120, 200, 500]\tau_n$ and mesh size $h=1/96$.
} \label{FULL_3D}
\end{figure} 

\subsection{Geometrical parametrization}\label{exp:Geometrical parametrization}
The numerical examples consider a geometrical parameterization on the embedded domain. The embedded domain is in fact parametrized through $\mu$ according to the expression:
$$
x^2 + y^2 \le \mu ^2/4,
$$
where the parameter $\mu$ describes the diameter size of the circular embedded domain located {\magenta{at}} the center $(0,0)$ of the background domain, see e.g. \autoref{background_mesh} (ii).
The experiments focus on the initial evolution period when the phase-field changes fast, namely in the time interval $0-100$ $\tau_n$  for time step size $\tau_n={\mathcal{O}(10^{-6})}$. We tested {\magenta{a}} pseudo-random initial concentration {\magenta{$u_0$}}. 
For all parameters, we experimented the same pseudo-random initial state applied on the whole background mesh and afterward restricted onto the active parametrized geometry. 
%
{\magenta{We recall that the aforementioned pseudorandom number generated sequence is not truly random, because it is completely determined by the ``seed": an initial value,  which may include truly random values, \cite{pseudo1,pseudo2}. These sequences that are closer to truly random can be generated using hardware random number generators. This pseudorandom set of values in our case is important in practice for their speed in number generation and their reproducibility. Actually, we use this generator to create a not favorable and non-smooth initial condition
using the aforementioned pseudorandom values between the values of the bulk phases. 
In this way, we can reproduce exactly the experiments using the aforementioned pseudorandom sequences.}}
In the first experiment the boundaries are free ({\magenta{Neumann}} boundary condition) everywhere, while in the second one we consider Dirichlet {\magenta{conditions}} only for the embedded geometry, \cite{LI2013102}. Though, only the {\magenta{circle}} is treated as embedded. Linear {\magenta{${P}^1\times{P}^1$}} polynomials have been employed for the discretization.
\begin{figure} \centering
\begin{minipage}{\textwidth}
\centering
\makebox[\textwidth][c]{
\begin{minipage}{0.225\textwidth}
  \includegraphics[width=\textwidth]{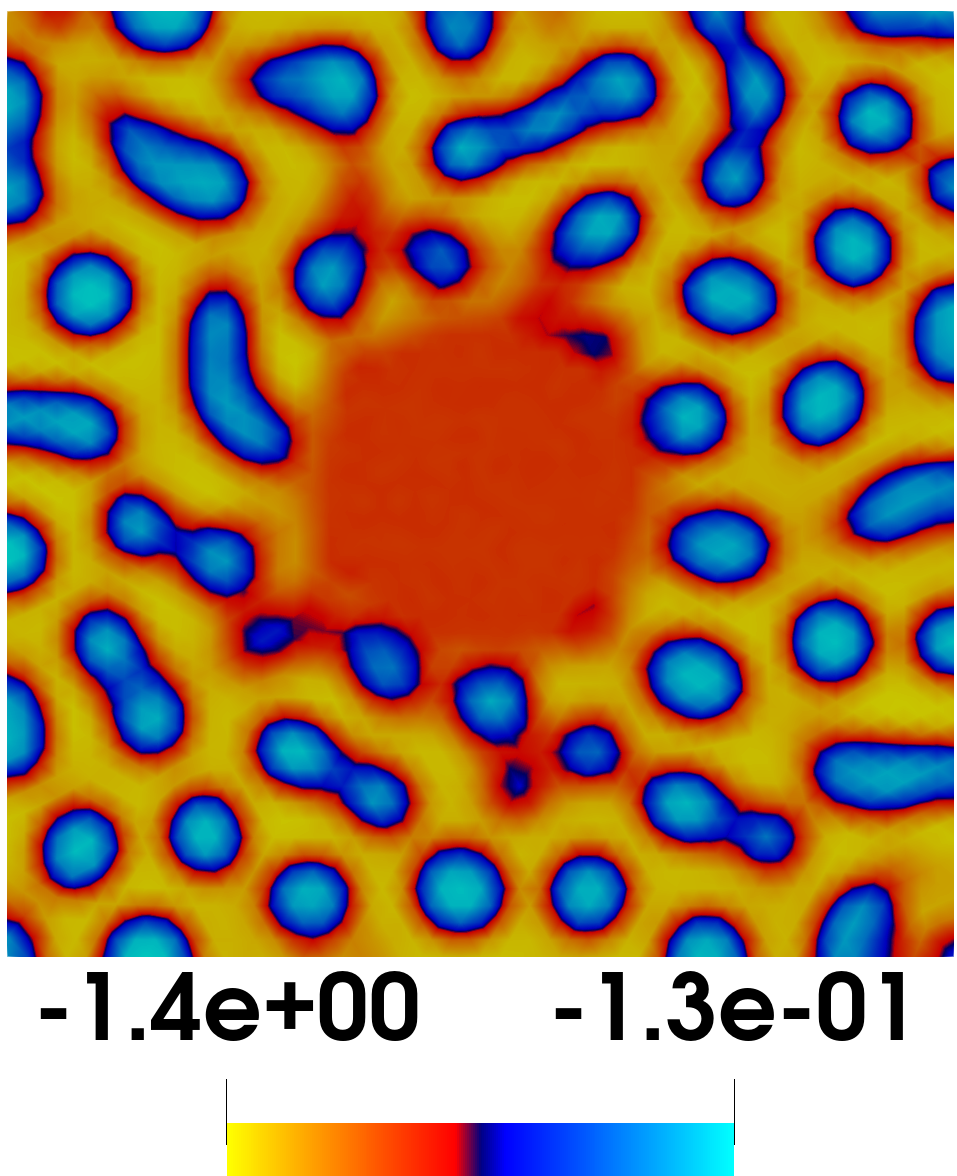} 
\end{minipage}
\hskip18pt
\begin{minipage}{0.225\textwidth}
  \includegraphics[width=\textwidth]{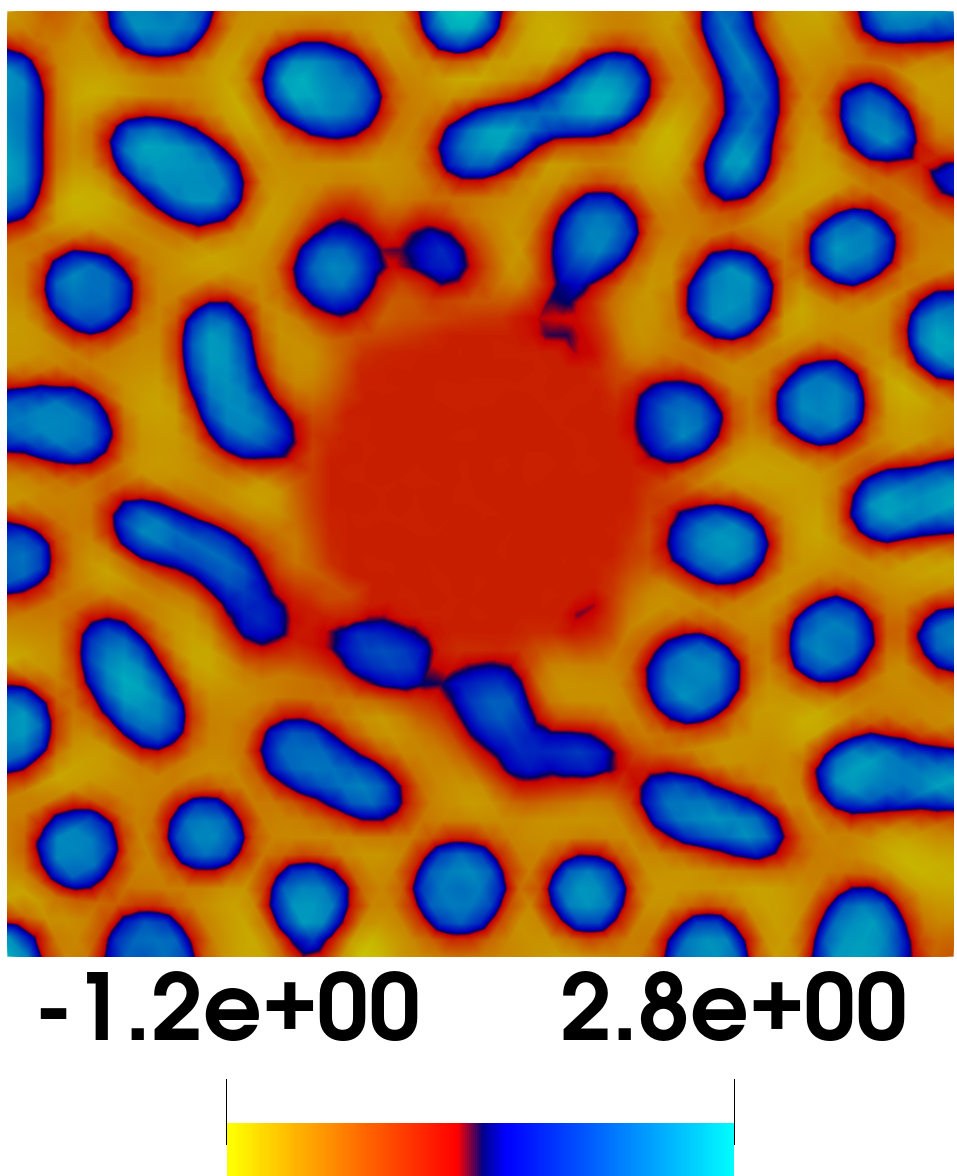}
\end{minipage}
\hskip18pt
\begin{minipage}{0.225\textwidth}
  \includegraphics[width=\textwidth]{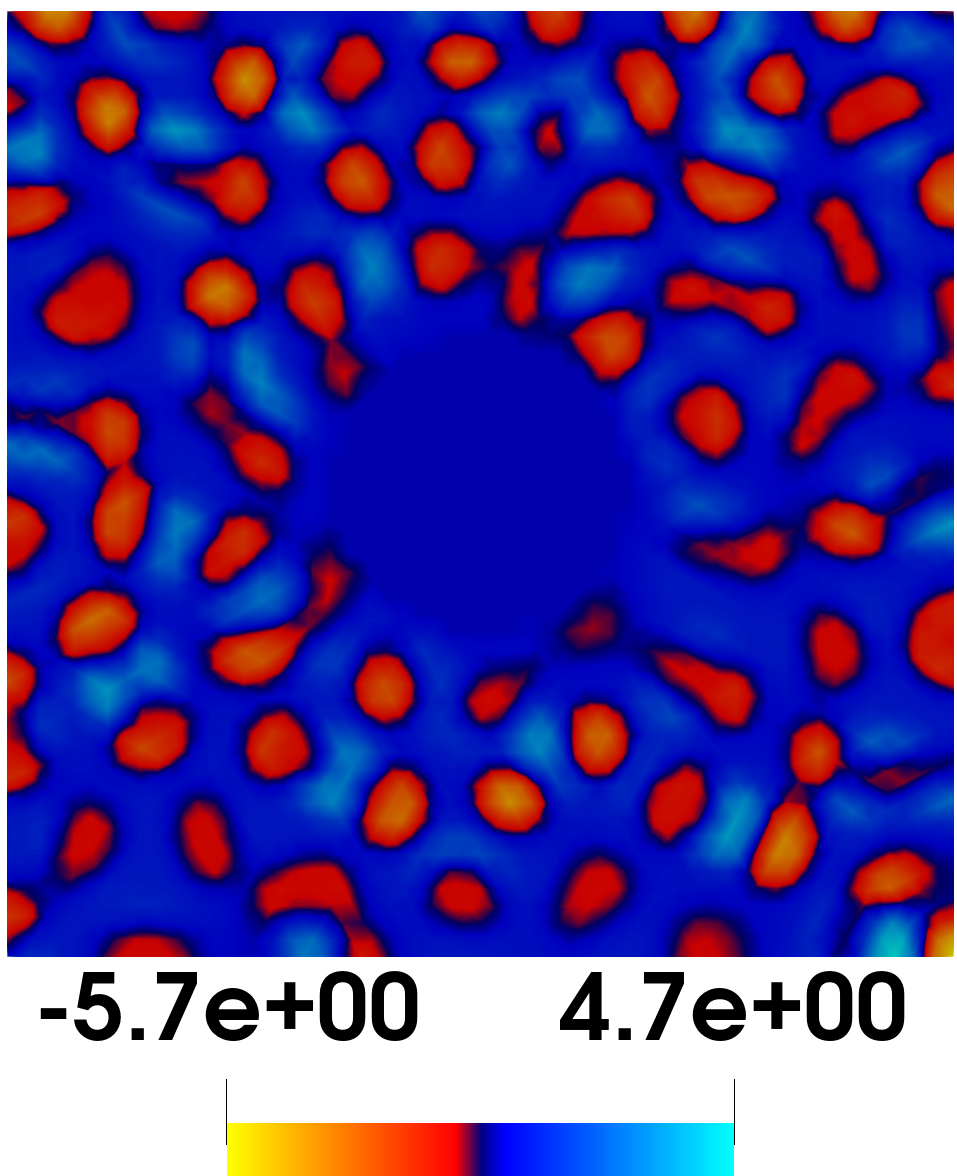}
\end{minipage}
}
\newline
\vskip20pt
\makebox[\textwidth][c]{%
\hskip4pt
\begin{minipage}{0.225\textwidth}
  \includegraphics[width=\textwidth]{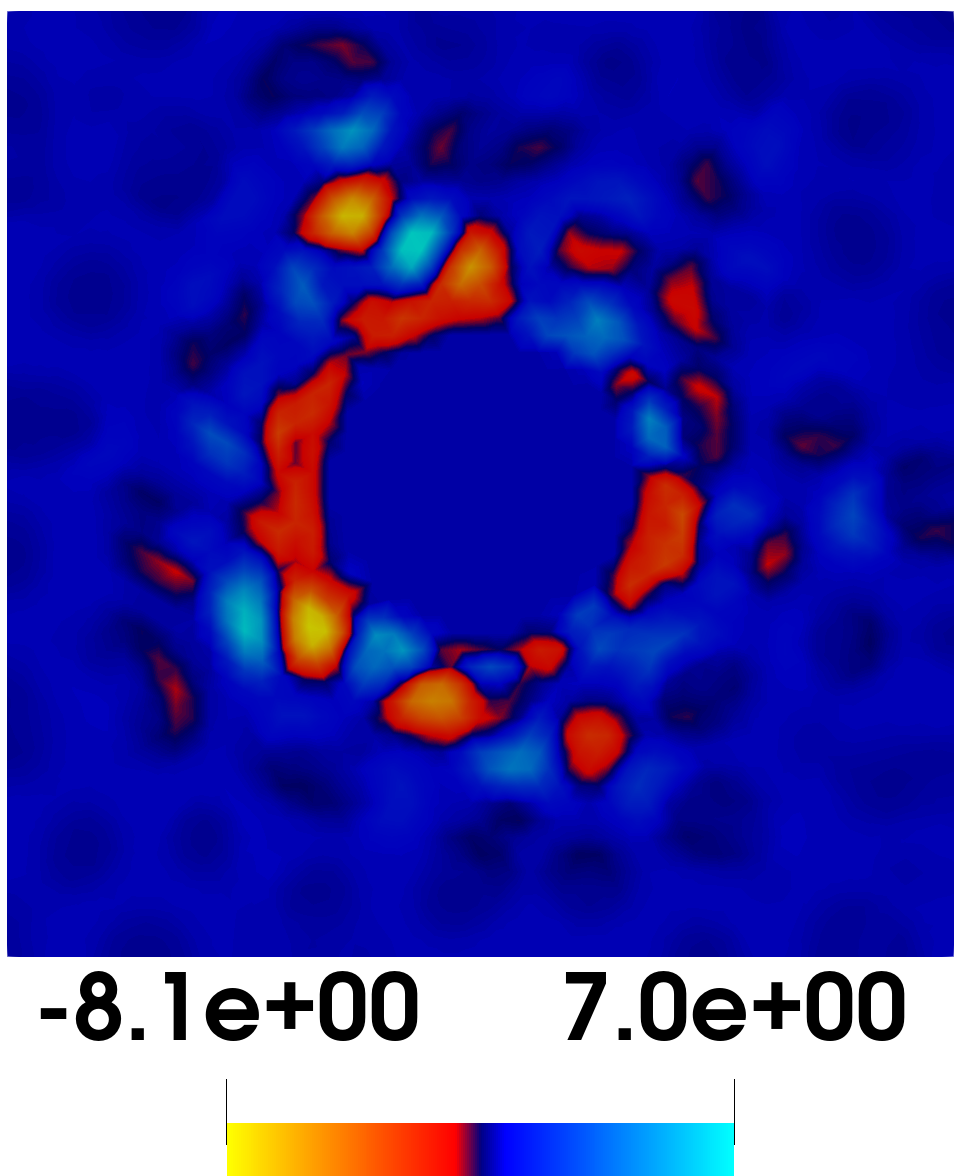}
\end{minipage}
\hskip18pt
\begin{minipage}{0.225\textwidth}
  \includegraphics[width=\textwidth]{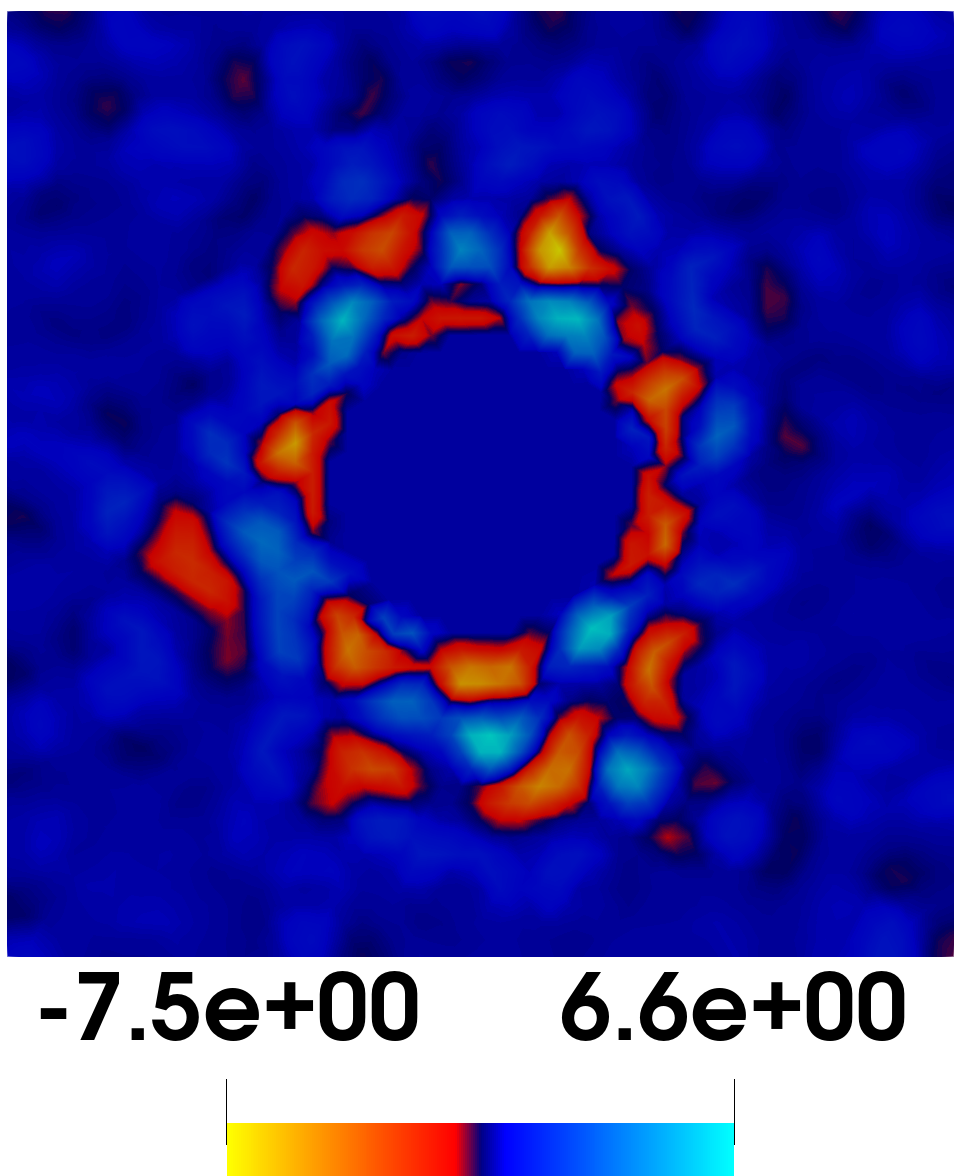}
\end{minipage}
\hskip18pt
\begin{minipage}{0.225\textwidth}
  \includegraphics[width=\textwidth]{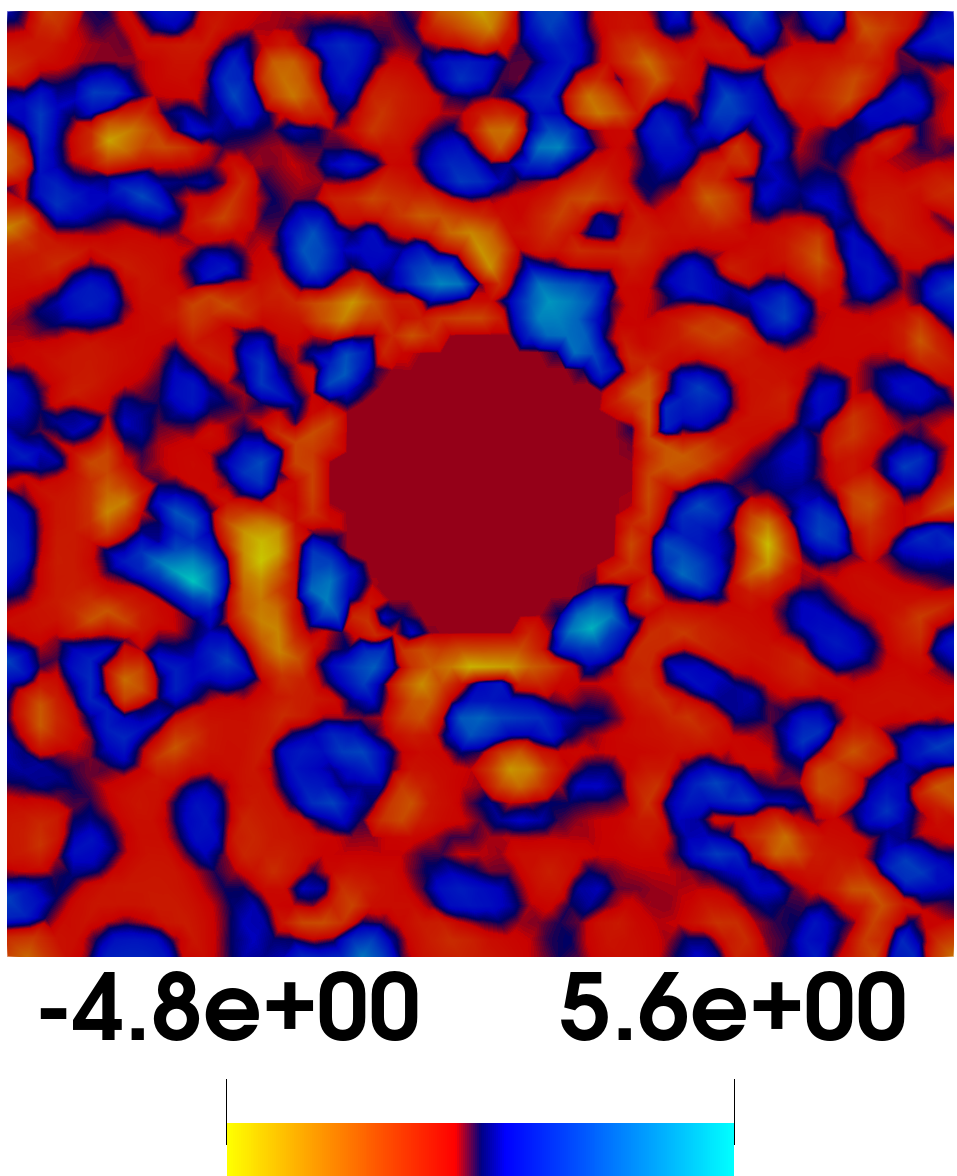}
\end{minipage}
}
\end{minipage}
\vskip20pt
\begin{minipage}{\textwidth}
\centering
\makebox[\textwidth][c]{
\begin{minipage}{0.225\textwidth}
  \includegraphics[width=\textwidth]{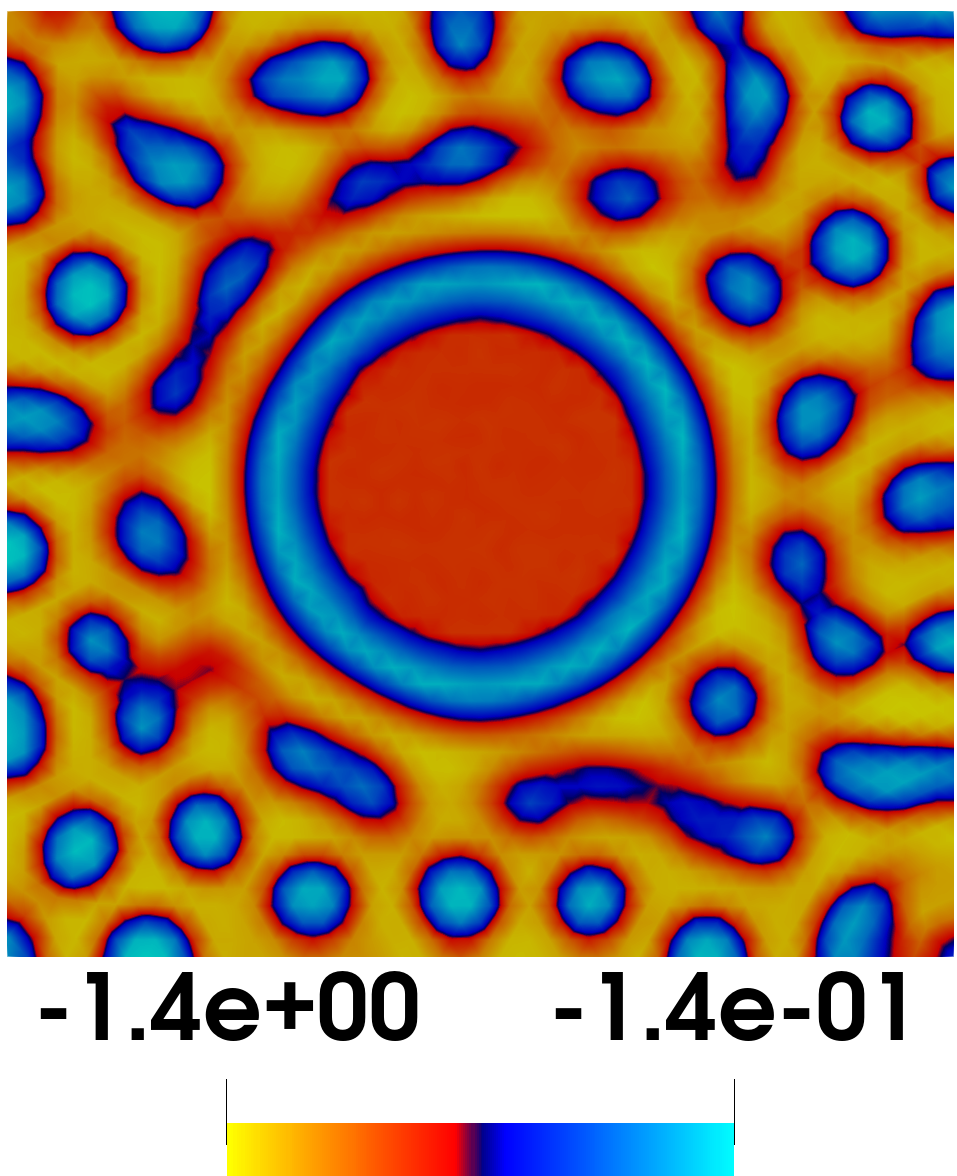} 
\end{minipage}
\hskip18pt
\begin{minipage}{0.225\textwidth}
  \includegraphics[width=\textwidth]{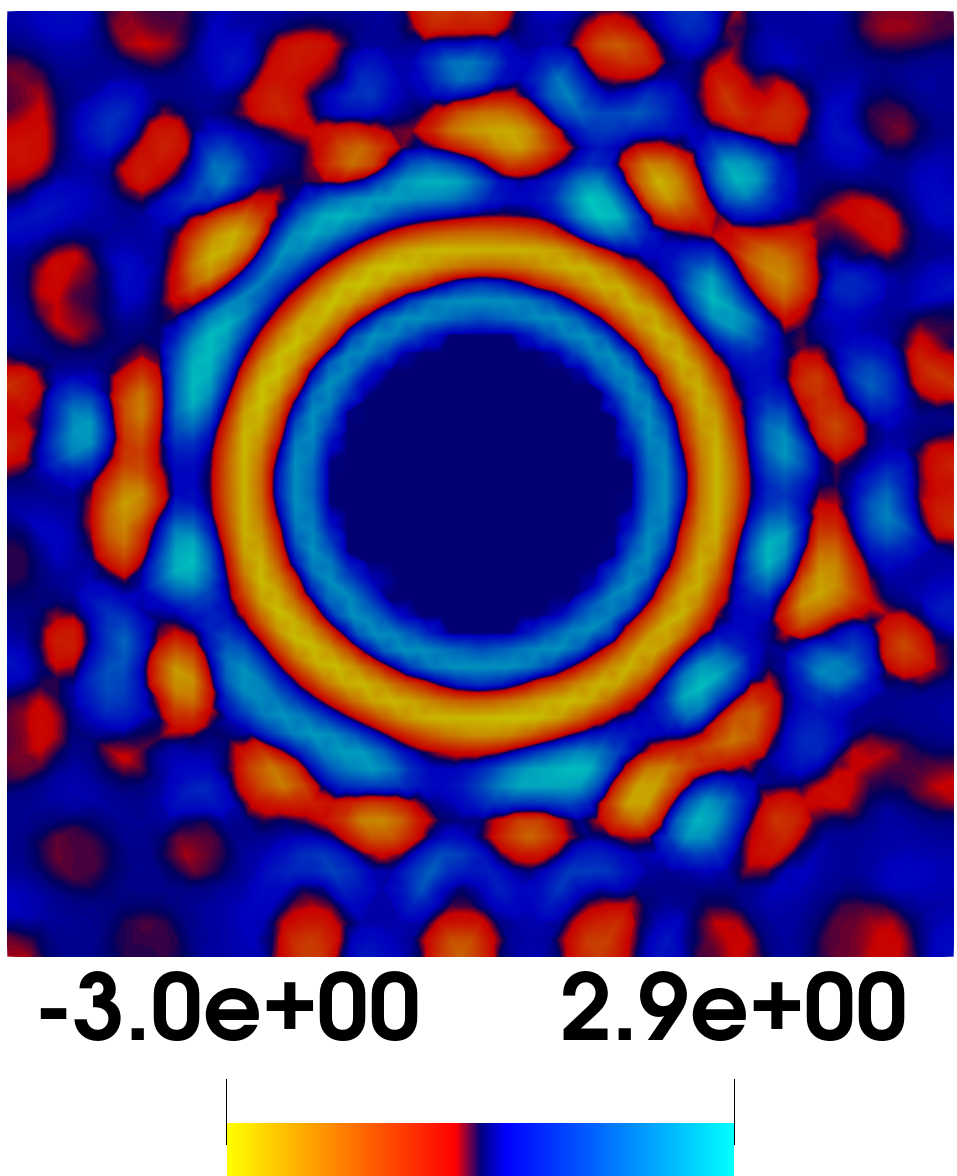}
\end{minipage}
\hskip18pt
\begin{minipage}{0.225\textwidth}
  \includegraphics[width=\textwidth]{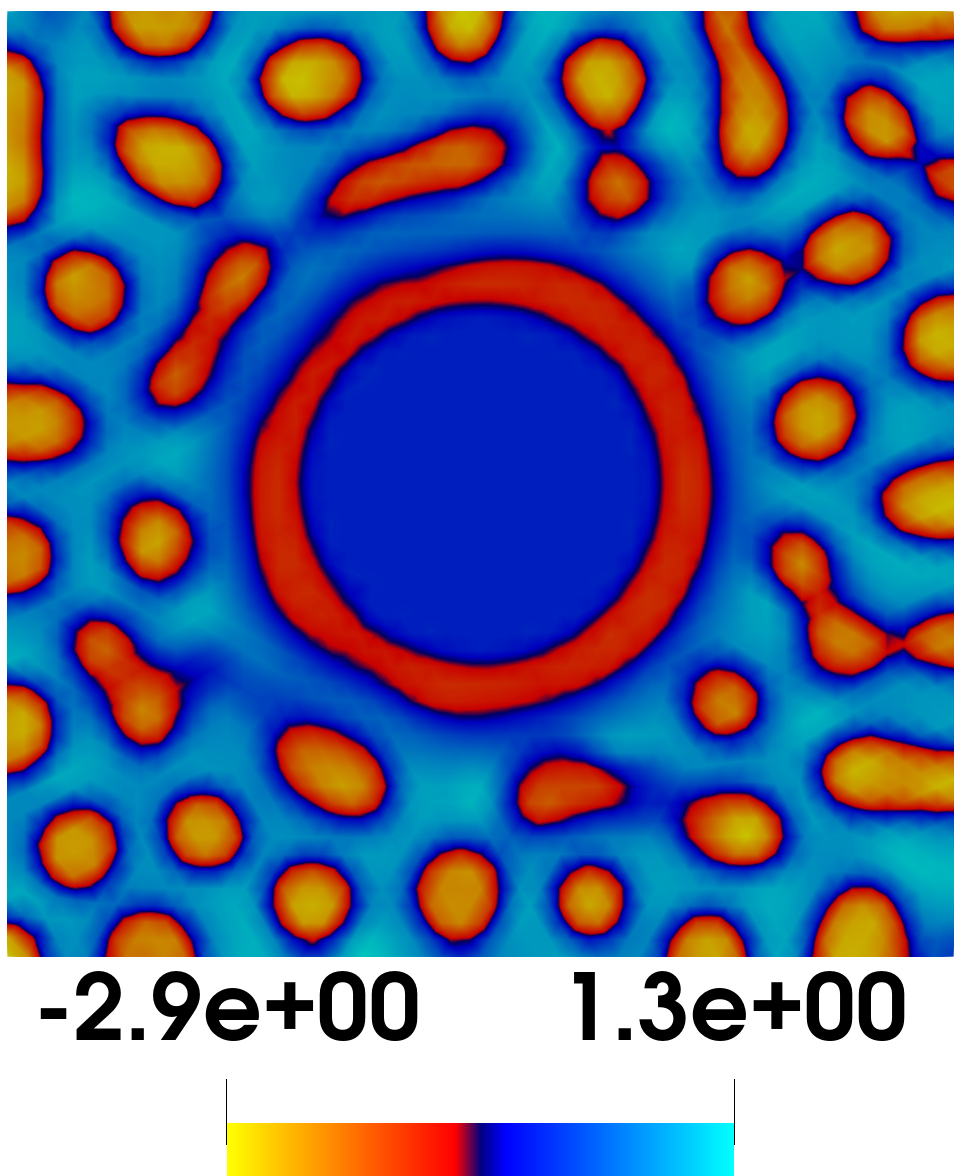}
\end{minipage}
}
\newline
\vskip20pt
\makebox[\textwidth][c]{%
\hskip4pt
\begin{minipage}{0.225\textwidth}
  \includegraphics[width=\textwidth]{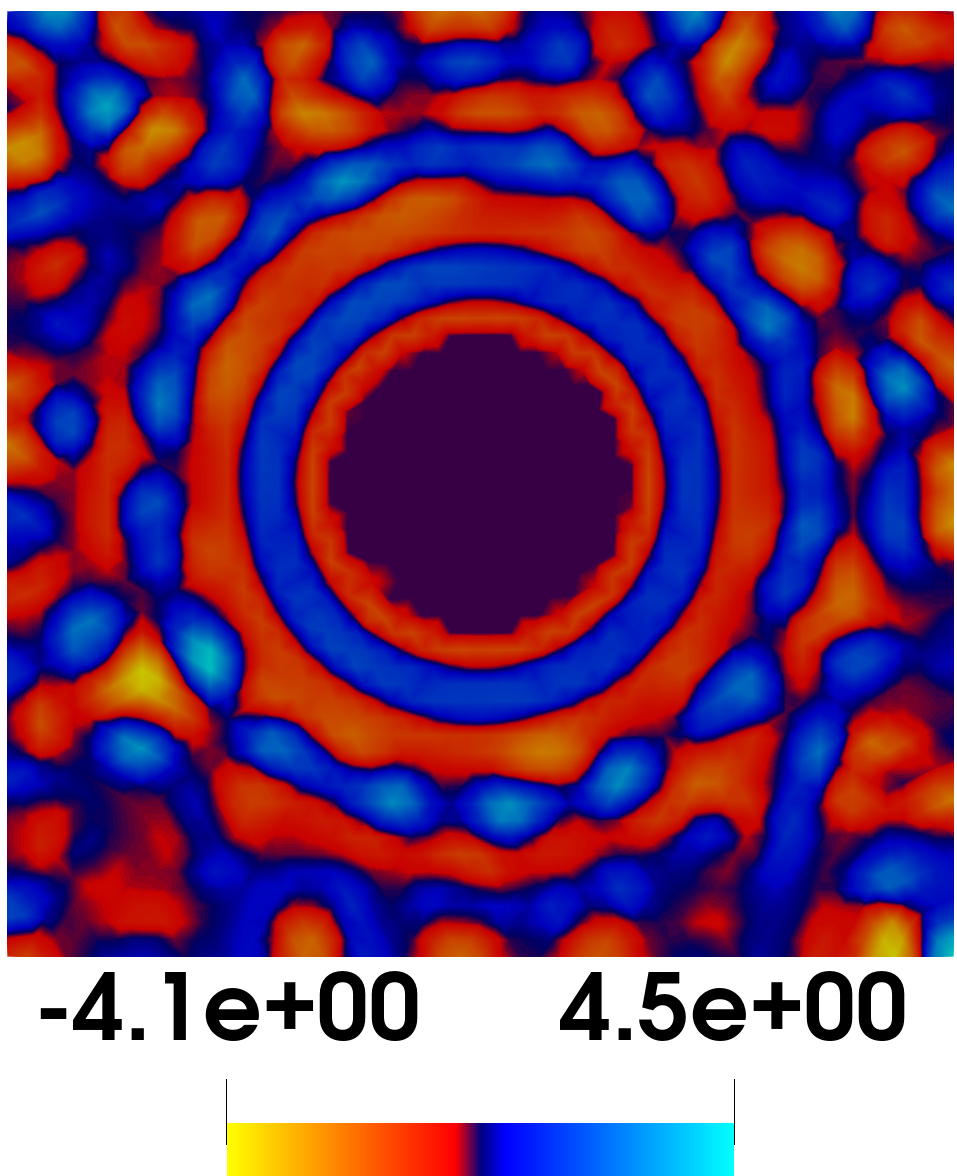}
\end{minipage}
\hskip18pt
\begin{minipage}{0.225\textwidth}
  \includegraphics[width=\textwidth]{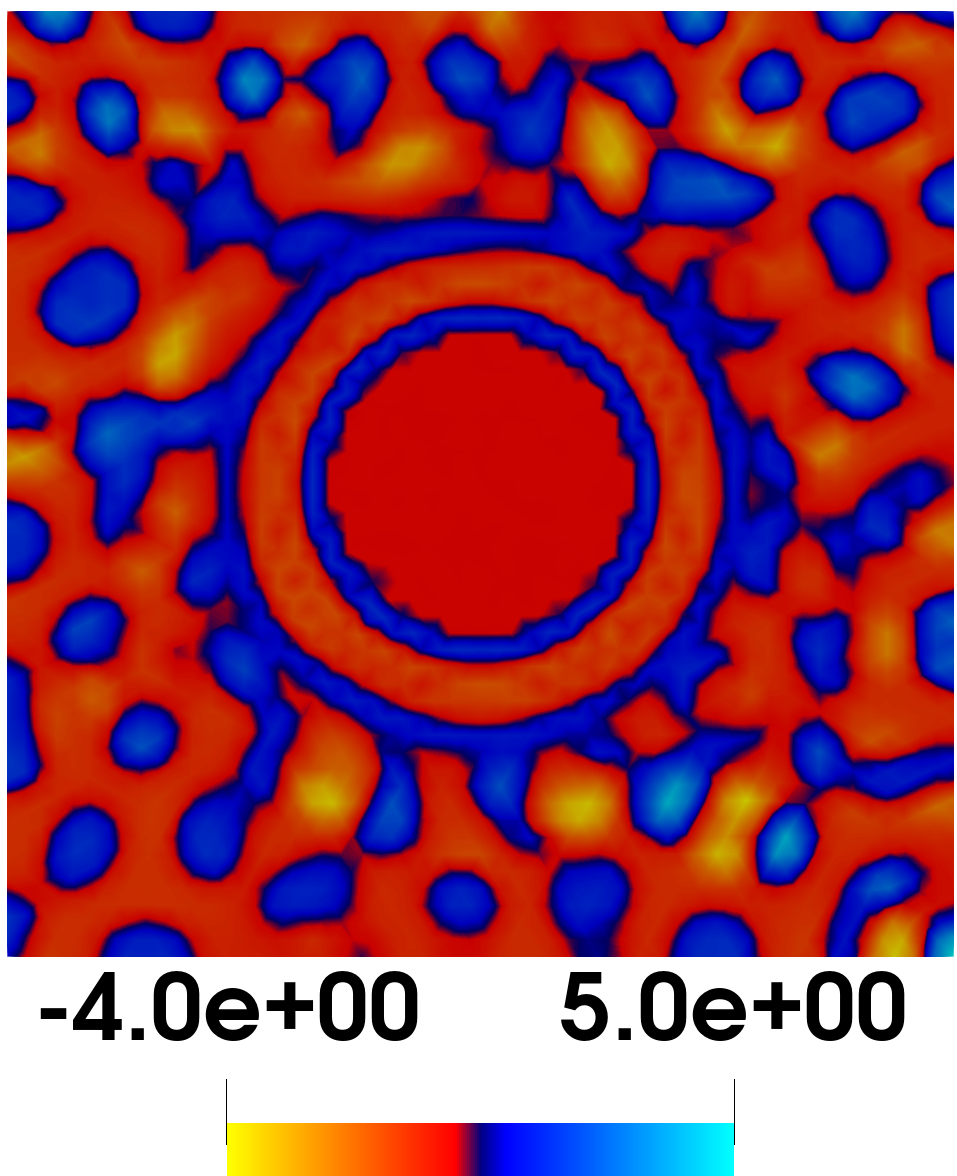}
\end{minipage}
\hskip18pt
\begin{minipage}{0.225\textwidth}
  \includegraphics[width=\textwidth]{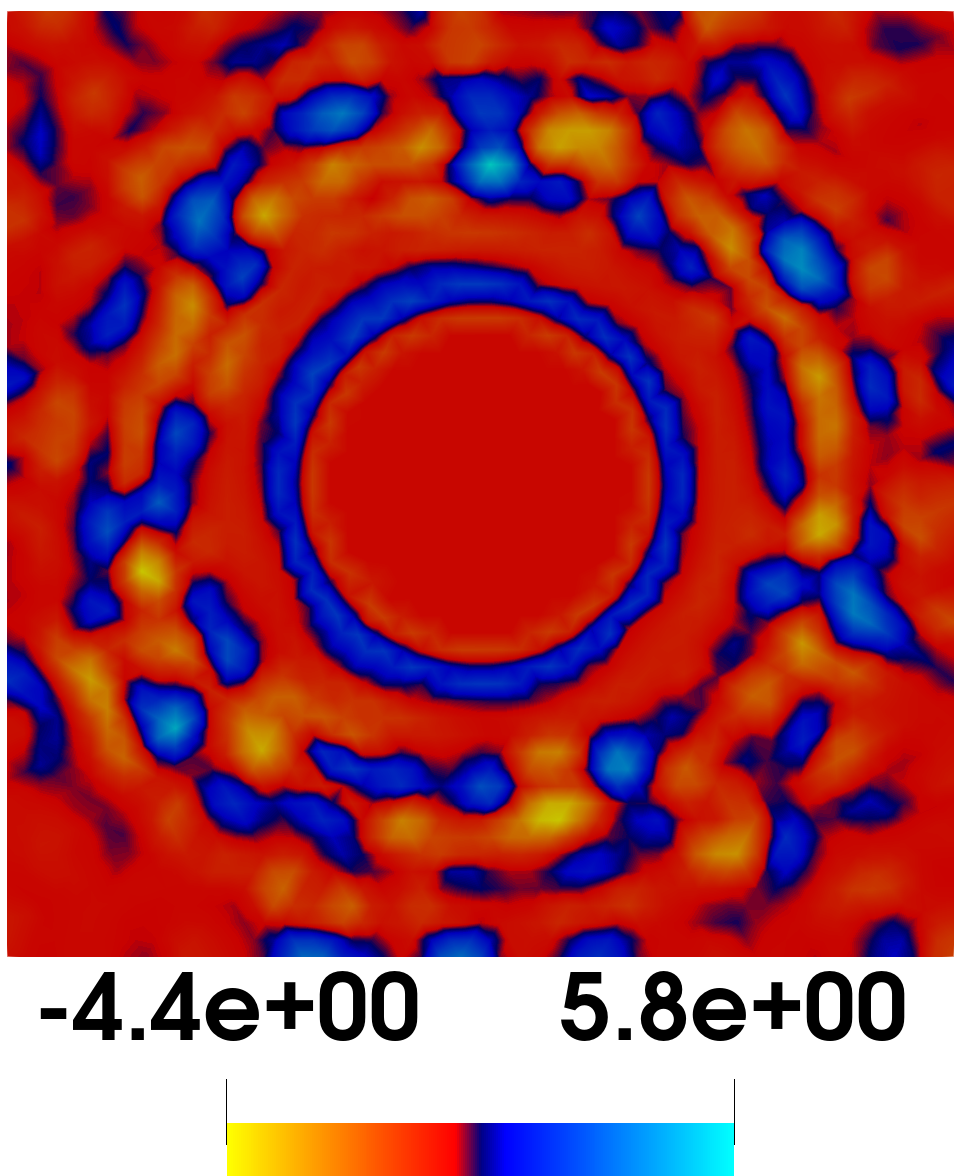}
\end{minipage}
}
\end{minipage}
  \caption{The first six basis {\magenta{functions}} {\magenta{plots refer to Neumann and the last six to Dirichlet parametrization experiments}}
  .}
  \label{Fig:Modes}
\end{figure} 
During the initial period of the evolution, $[0, T]=[0,100\tau_n]$, {\magenta{one}} can notice fast phase-field changes starting from pseudo-random initial data, both challenging for the ROM construction. Nevertheless, as time passes these changes {\magenta{weaken}},  see e.g. Figure \ref{FULL_3D} for a visualization of the evolution of an {\magenta{Neumann}} and a Dirichlet type embedded boundary.

For the reduced basis solution, the ROM has been trained onto $900$ parameter {\magenta{samples  chosen} randomly inside the parameter space. From this {\magenta{snapshot selection, and for the snapshot matrices (\ref{Su})-(\ref{Sw}), a ROM basis is derived}} with the}} POD procedure {\magenta{(\ref{eq:pod_energy})-(\ref{eq:pod_energy2}) and parameter  range}} ${\magenta{\mu_{\text{train}} \in }}[0.36, 0.48]$. The visualization of the first six basis {\magenta{functions}} can be seen in Figure \ref{Fig:Modes}. {\magenta{We underline that the  ROM is tested with parameters not contained in the training set}}. To test the accuracy of the ROM we compared {\magenta{the FOM, with the ROM solutions for}} 
$30$ additional samples {\magenta{in a parameter range $\mu_{\text{test}} \in {\magenta{[0.40, 0.44]}}$}} which were not used to create the ROM and were selected randomly within the aforementioned range. 
{{We notice that whether we increase the range of the parameter and also enlarge the test space, the method reacts well, although the errors especially for the auxiliary variable $w$ start to increase. Important is though, that the concentration relative error for which we are interested is kept low. 
For such numerical testing, we refer to the Section \ref{subsec:challenges} and the paragraph related to the part where a larger range of the parameter is considered as well as the training and testing parameters  $\mu$ are in intervals of much larger and/or much smaller values.}}

\begin{figure}
\centering 
\begin{minipage}{\textwidth}
\centering
\begin{minipage}{0.25\textwidth}
  \includegraphics[width=\textwidth]{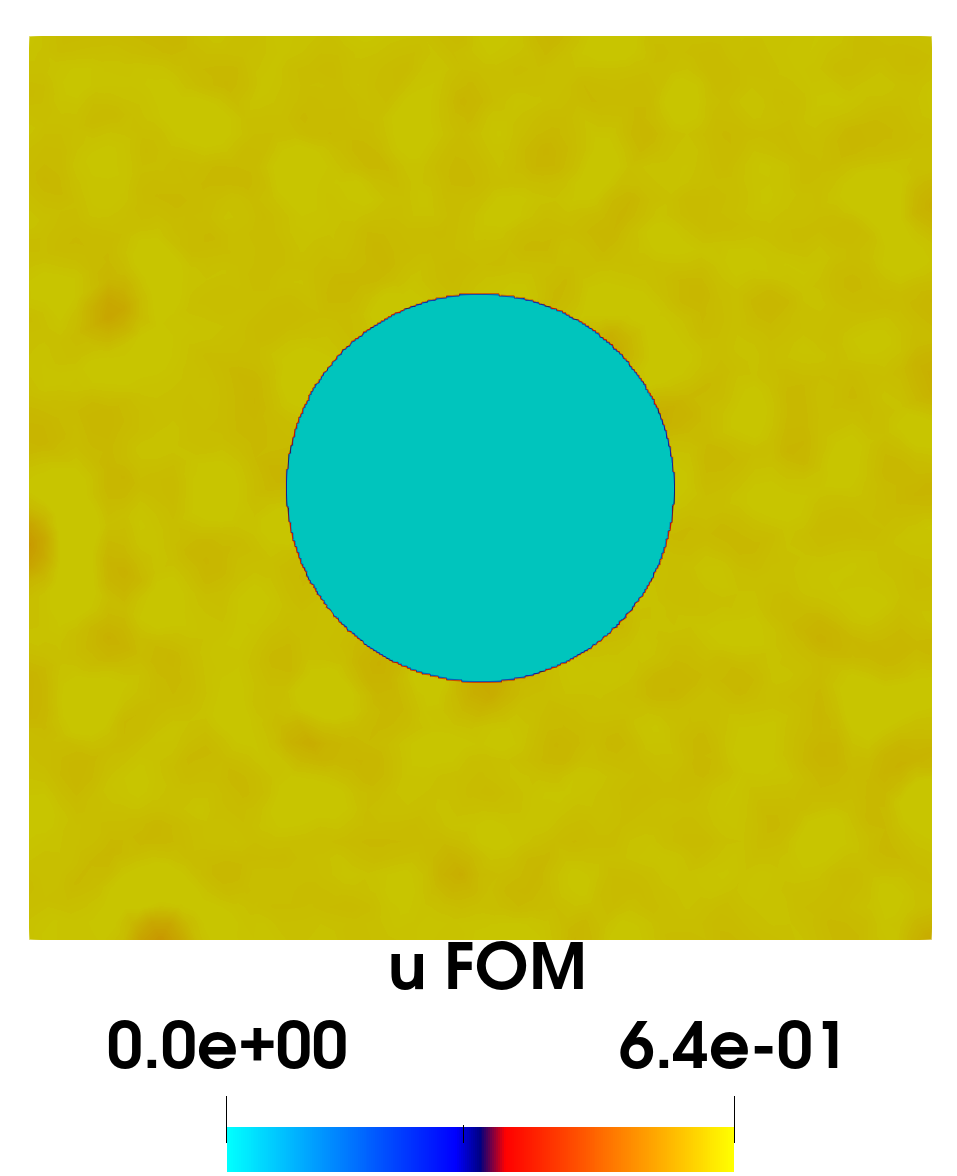}
\end{minipage}
 \qquad 
 \begin{minipage}{0.25\textwidth}
  \includegraphics[width=\textwidth]{{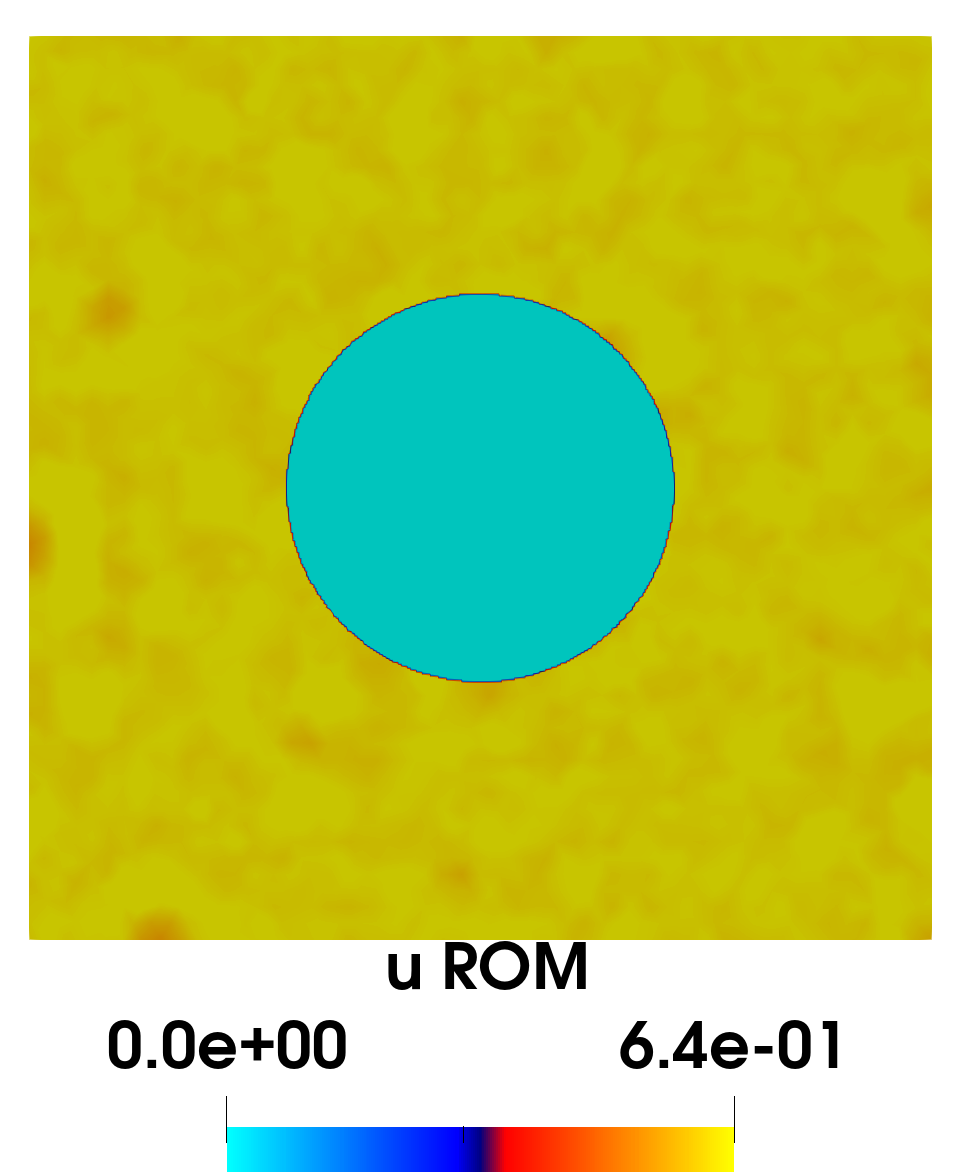}} 
\end{minipage}
 \qquad 
 \begin{minipage}{0.25\textwidth}
  \includegraphics[width=\textwidth]{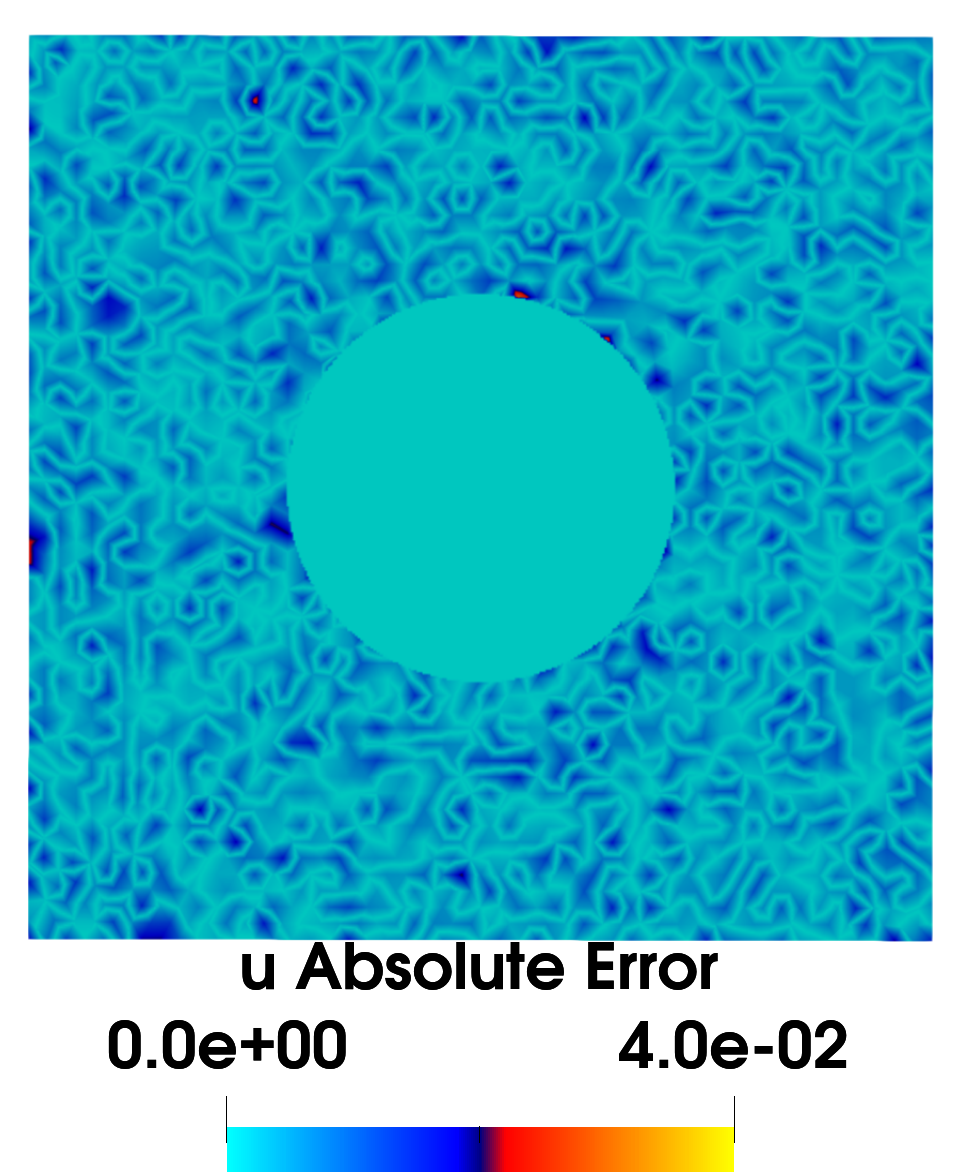}
\end{minipage}
\begin{minipage}{0.25\textwidth}
  \includegraphics[width=\textwidth]{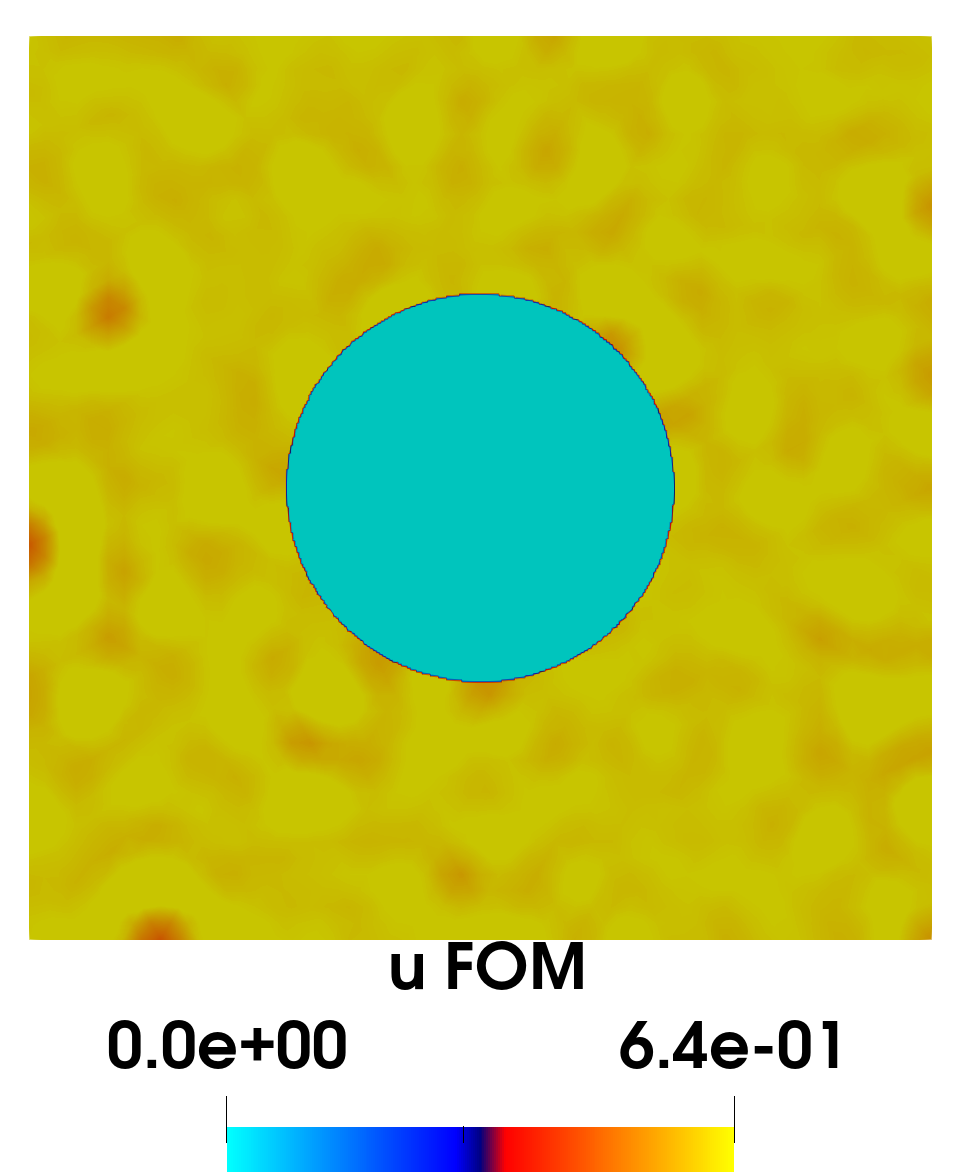}
\end{minipage}
\qquad
\begin{minipage}{0.25\textwidth}
  \includegraphics[width=\textwidth]{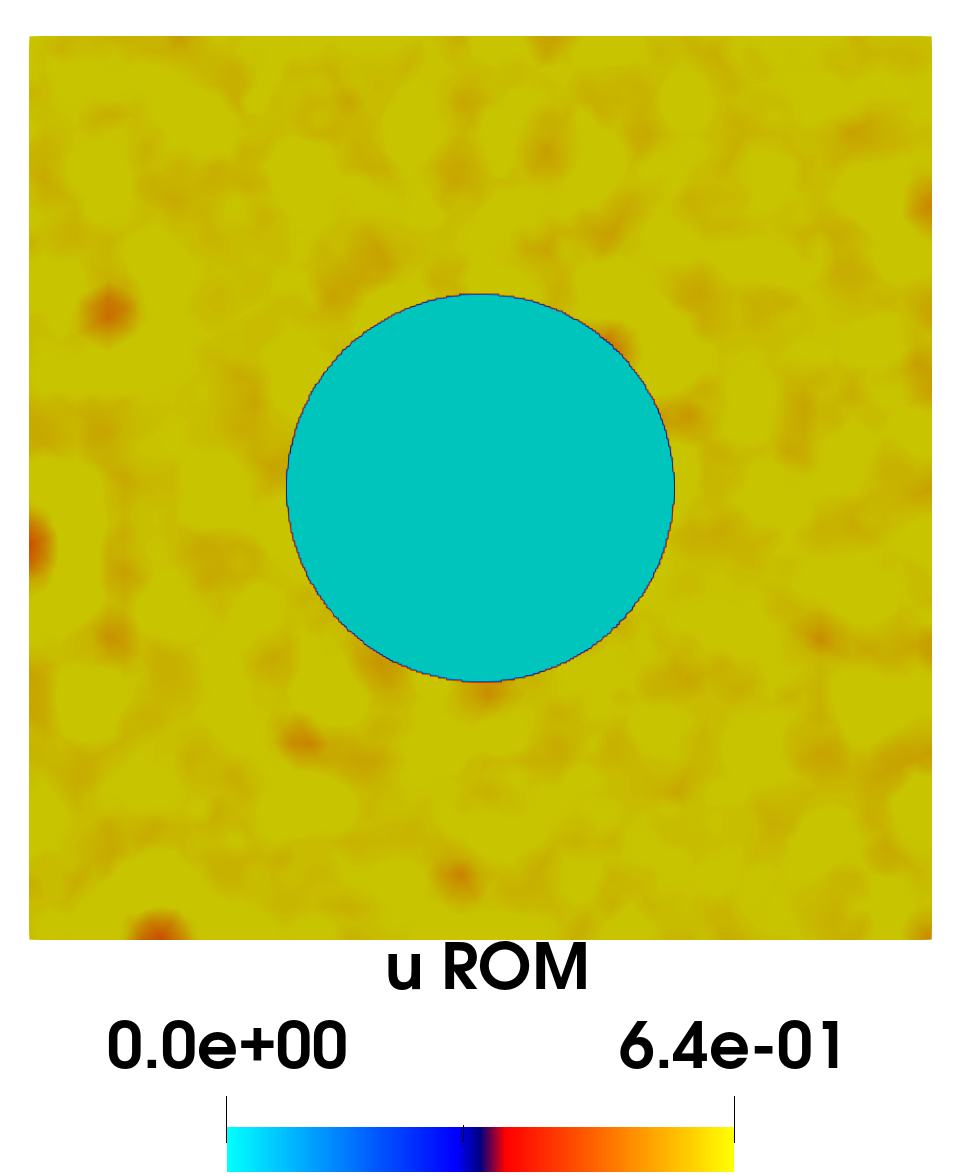} 
\end{minipage}
\qquad
\begin{minipage}{0.25\textwidth}
  \includegraphics[width=\textwidth]{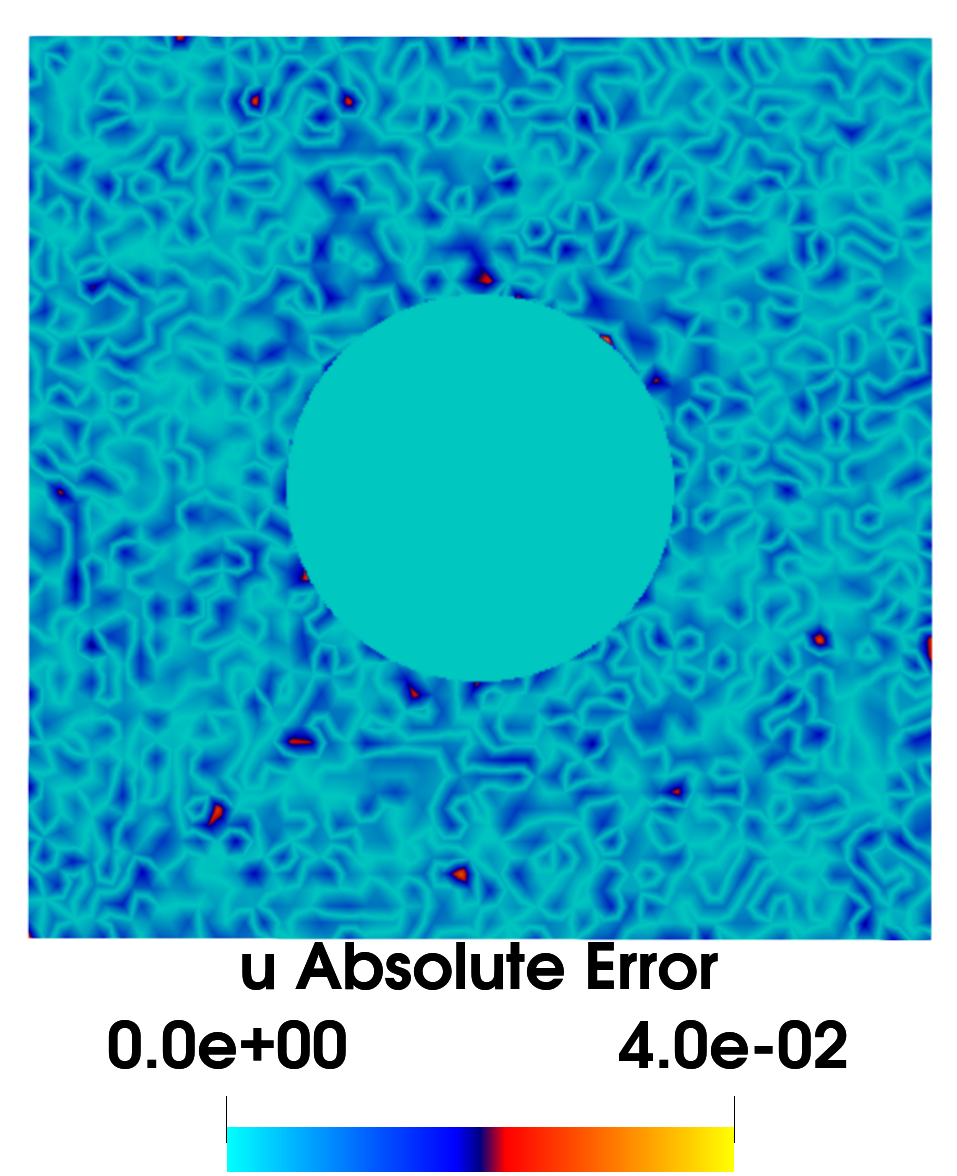}
\end{minipage}
\begin{minipage}{0.25\textwidth}
  \includegraphics[width=\textwidth]{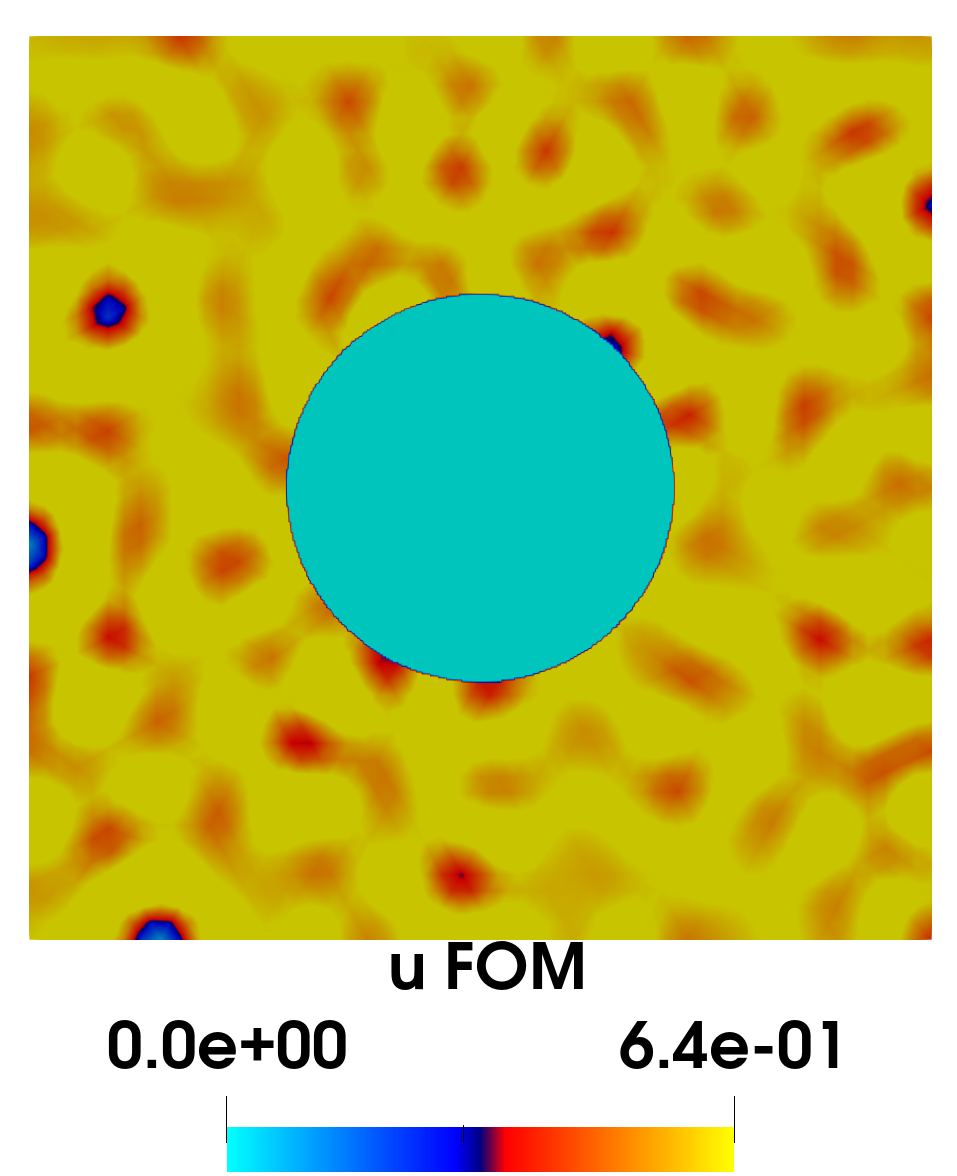}
\end{minipage}
\qquad
\begin{minipage}{0.25\textwidth}
  \includegraphics[width=\textwidth]{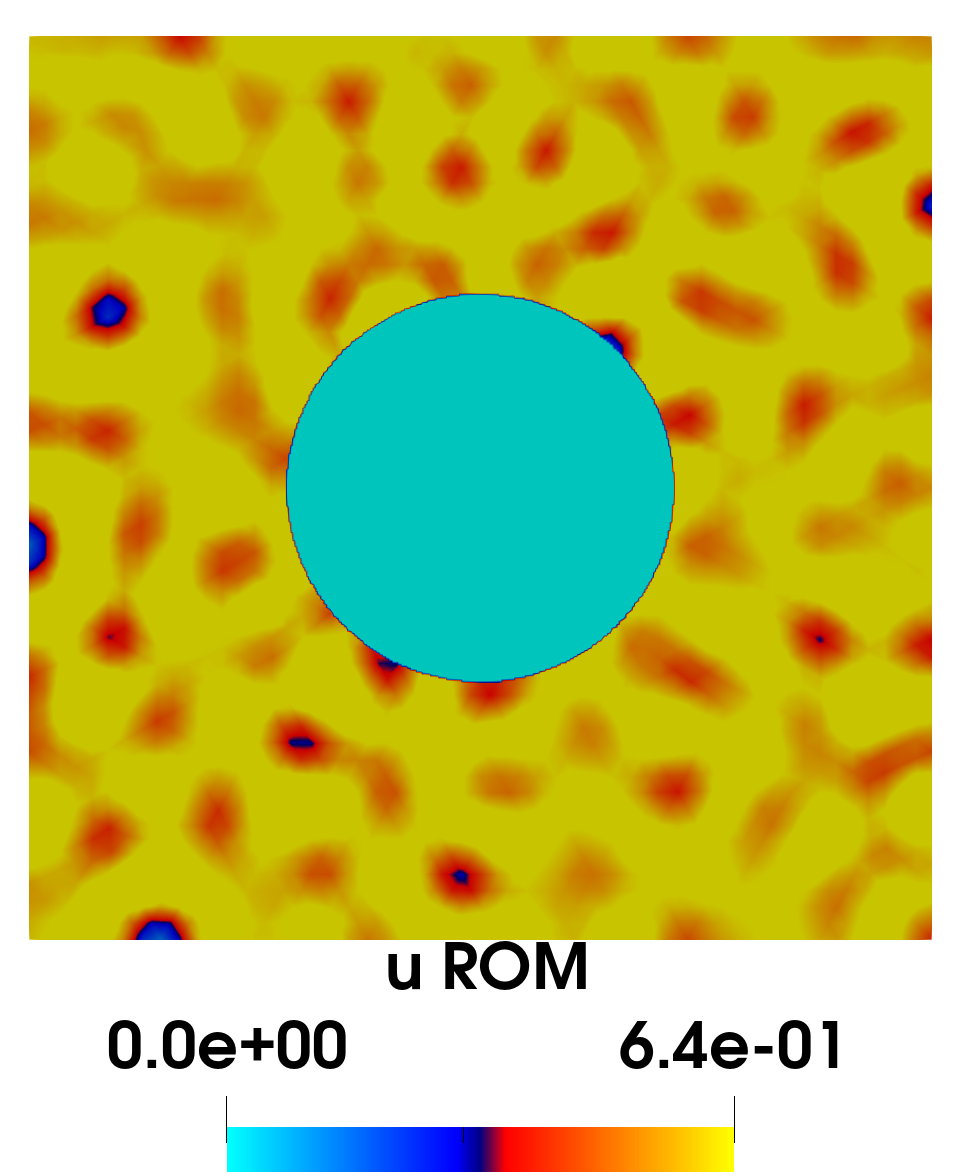}
\end{minipage}
\qquad
\begin{minipage}{0.25\textwidth}
  \includegraphics[width=\textwidth]{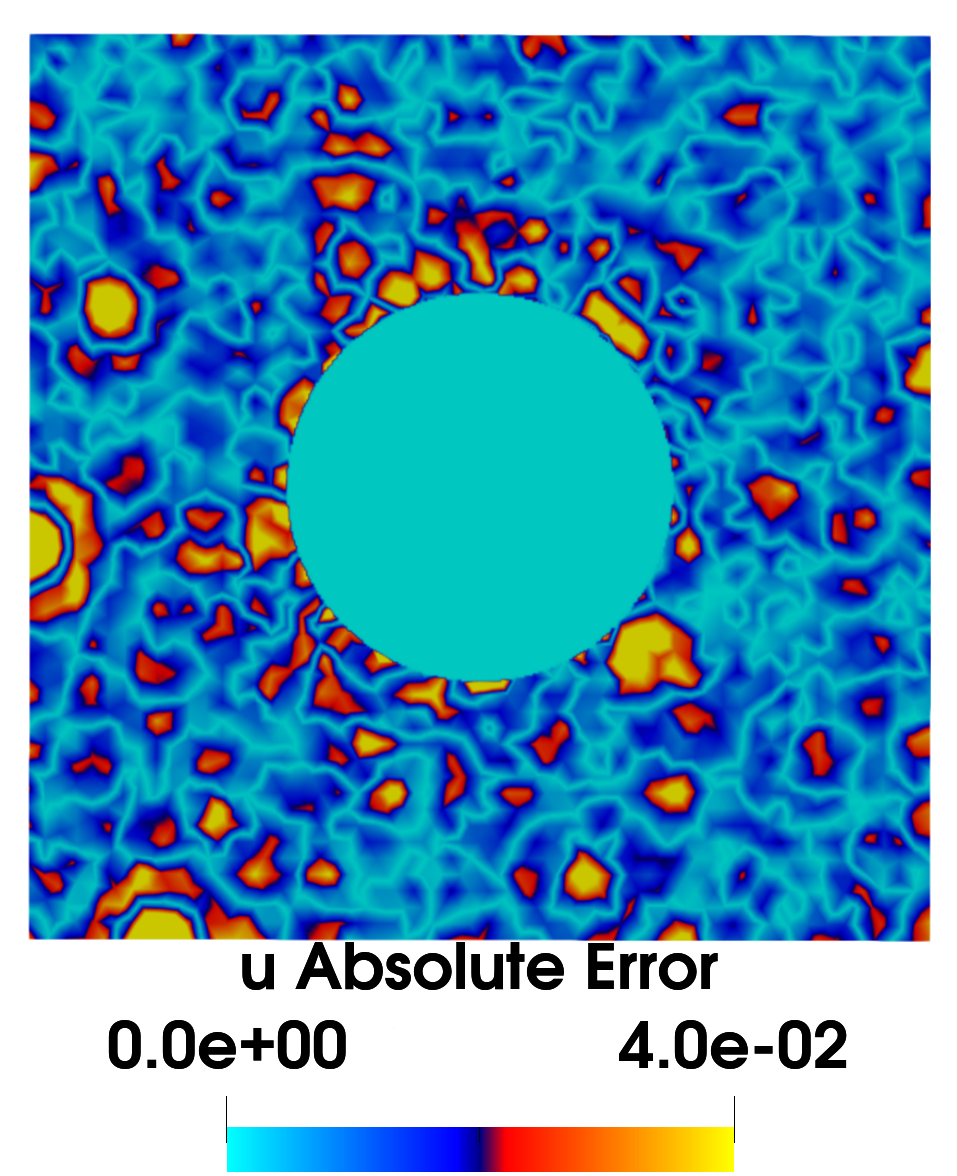}
\end{minipage}
\begin{minipage}{0.25\textwidth}
  \includegraphics[width=\textwidth]{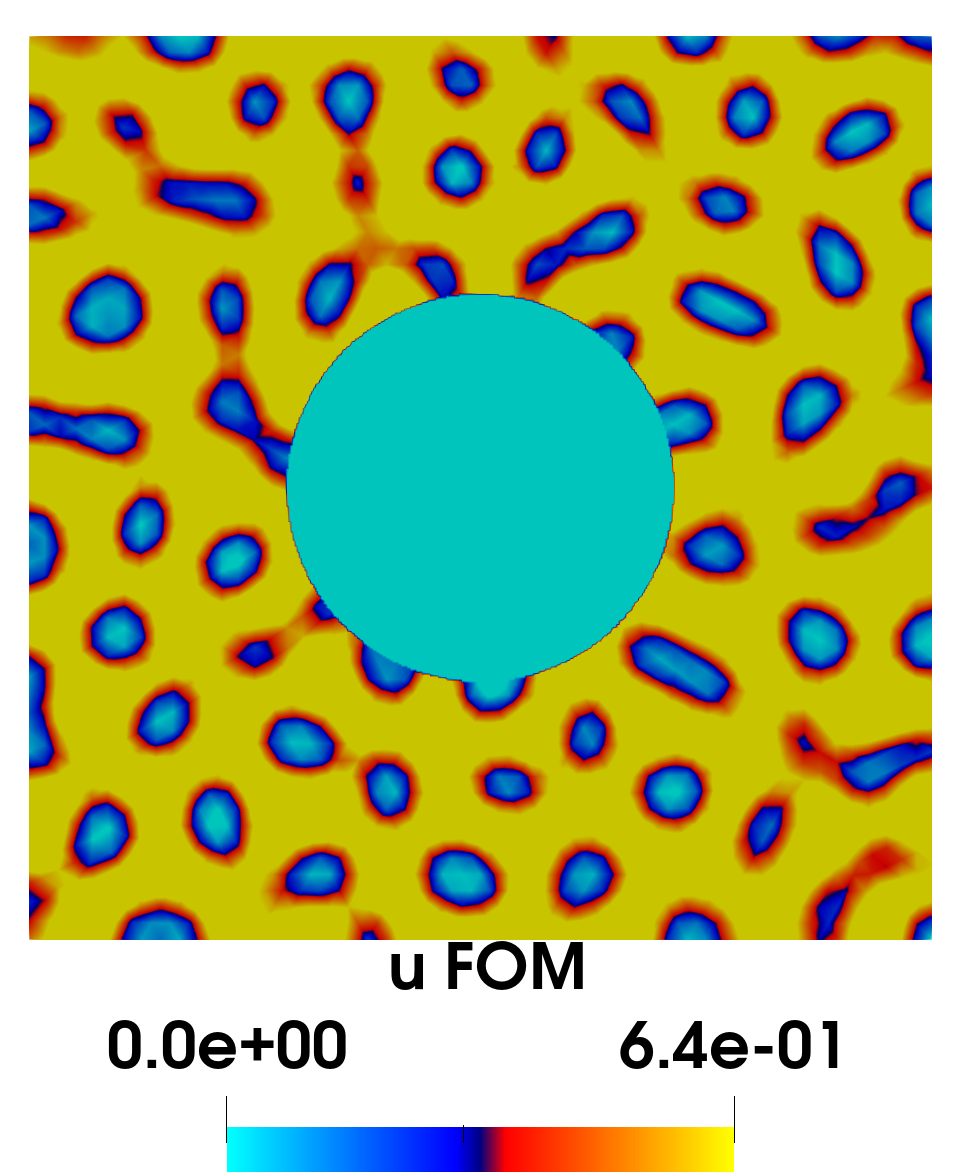}
\end{minipage}
\qquad
\begin{minipage}{0.25\textwidth}
  \includegraphics[width=\textwidth]{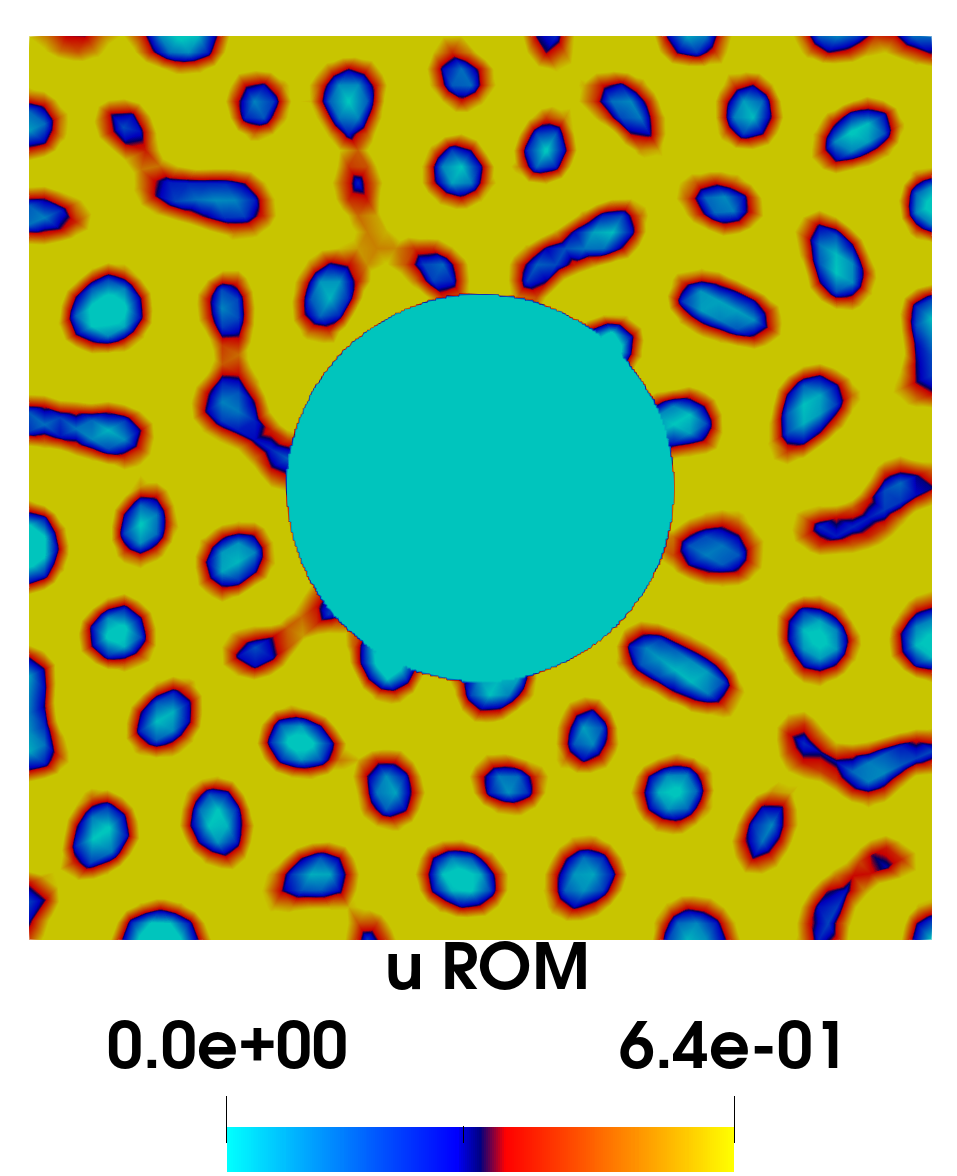}
\end{minipage}
\qquad
\begin{minipage}{0.25\textwidth}
  \includegraphics[width=\textwidth]{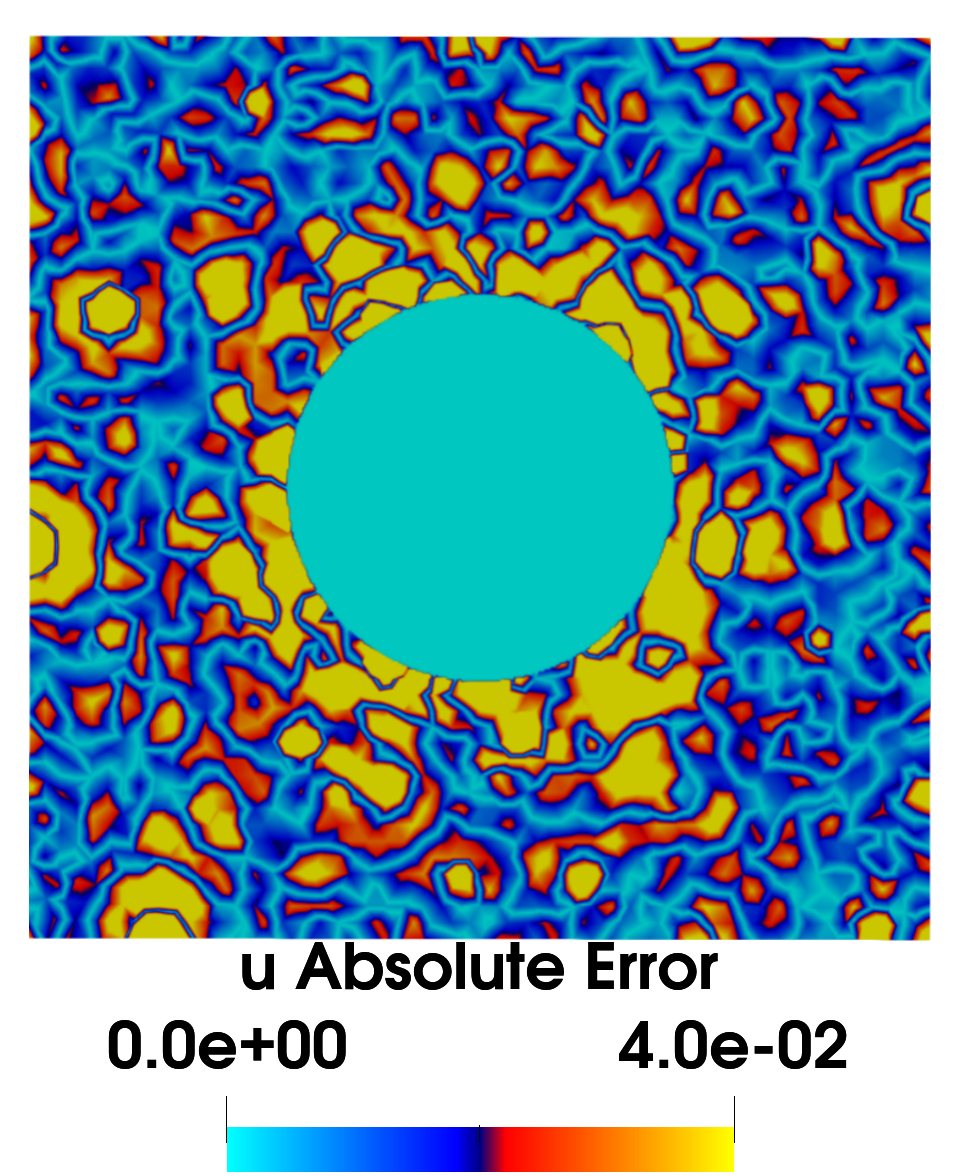}
  \end{minipage}
  \begin{minipage}{0.25\textwidth}
  \includegraphics[width=\textwidth]{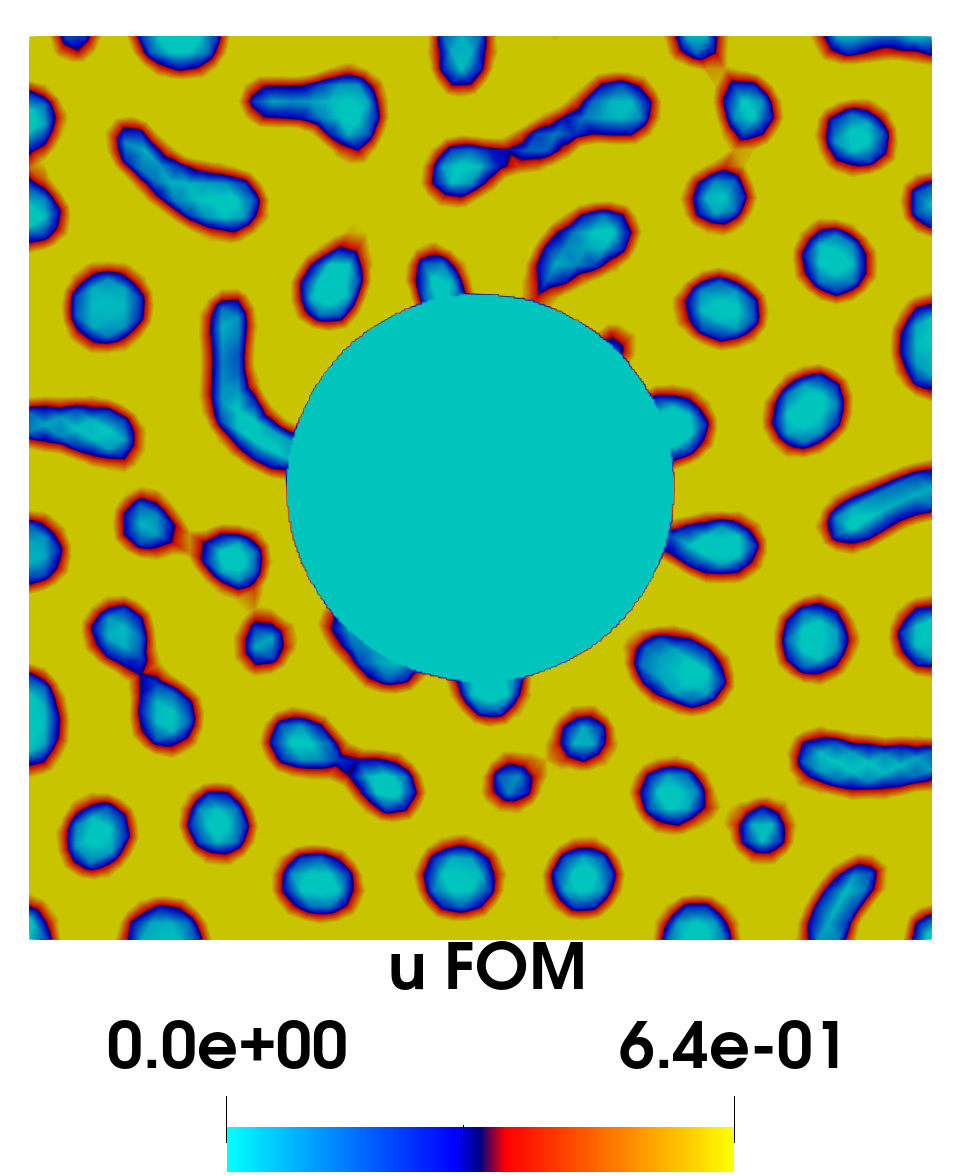}
\end{minipage}
\qquad
\begin{minipage}{0.25\textwidth}
  \includegraphics[width=\textwidth]{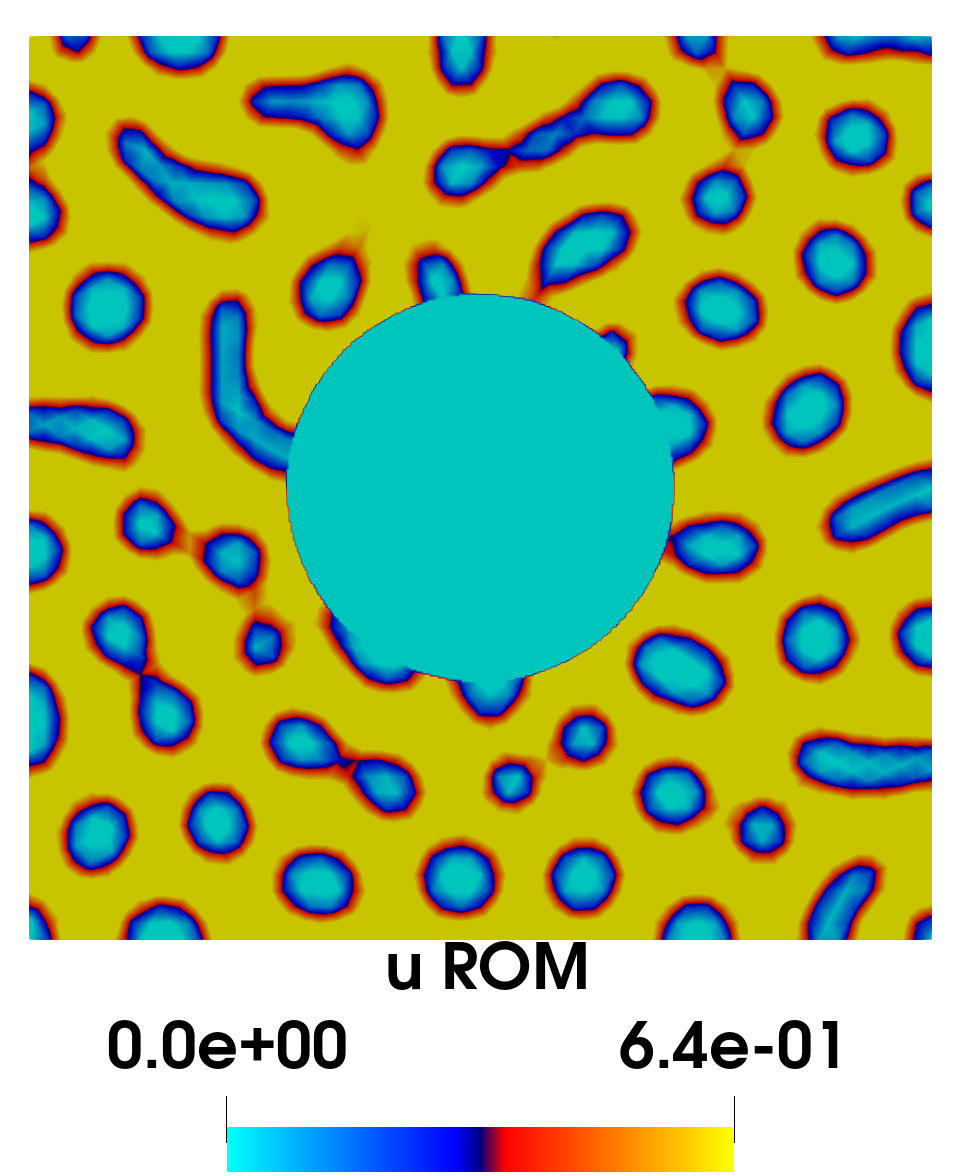}
\end{minipage}
\qquad
\begin{minipage}{0.25\textwidth}
  \includegraphics[width=\textwidth]{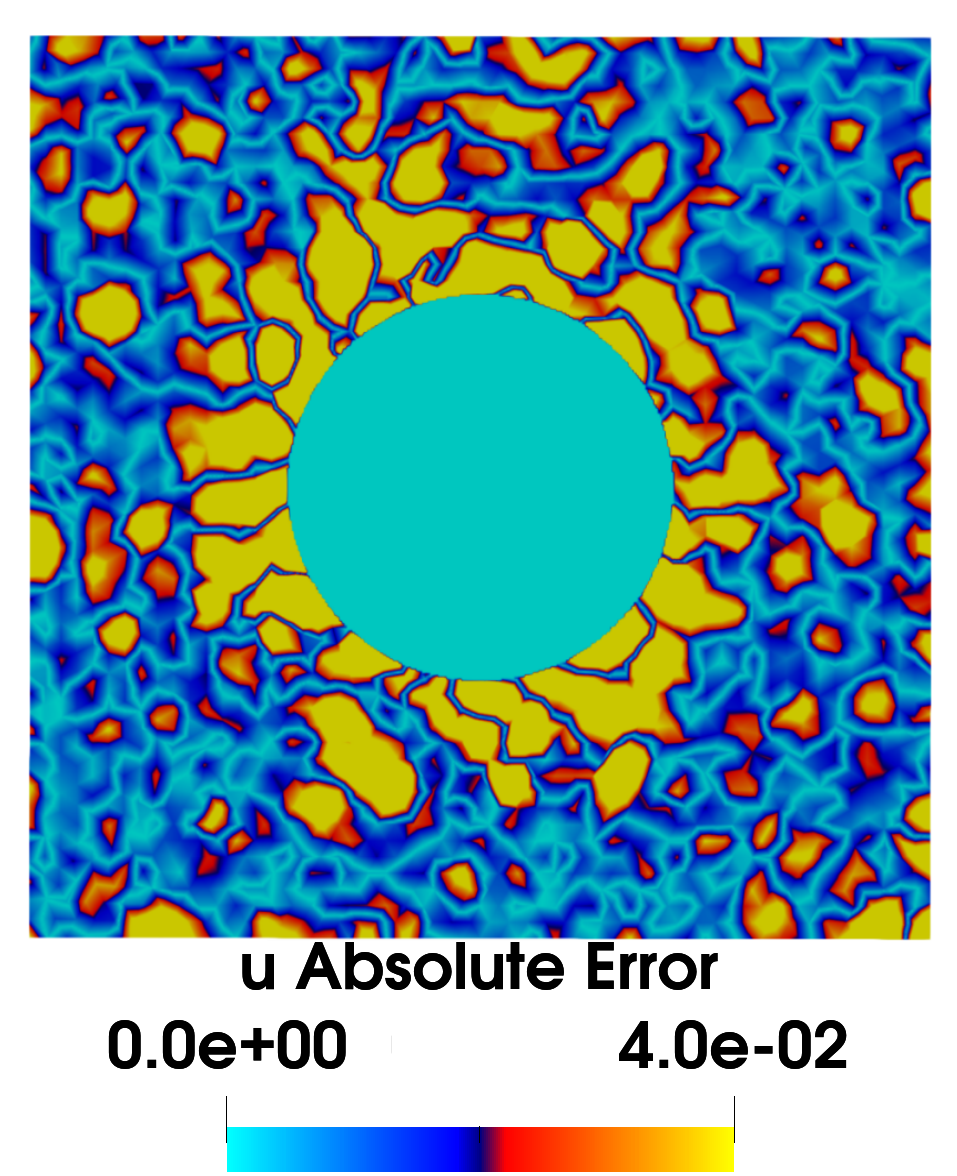}
\end{minipage} 
\end{minipage}
\caption{Results for the geometrical parametrized embedded circle  with diameter $\mu_{\text{test}} =0.42917$ and {\textit{{\magenta{Neumann}} boundary}}. 
In the first, second and third column we report the {\magenta{full order}} solution, the reduced order solution and the {\blue{absolute}} error plots for the concentration field. Each row corresponds to a different time $t=[10,20,40,60,100]\tau_n$.} \label{FULL_RED_ERROR_1P}
\end{figure}
\begin{figure} 
\centering
\begin{minipage}{\textwidth}
\centering
\begin{minipage}{0.25\textwidth}
  \includegraphics[width=\textwidth]{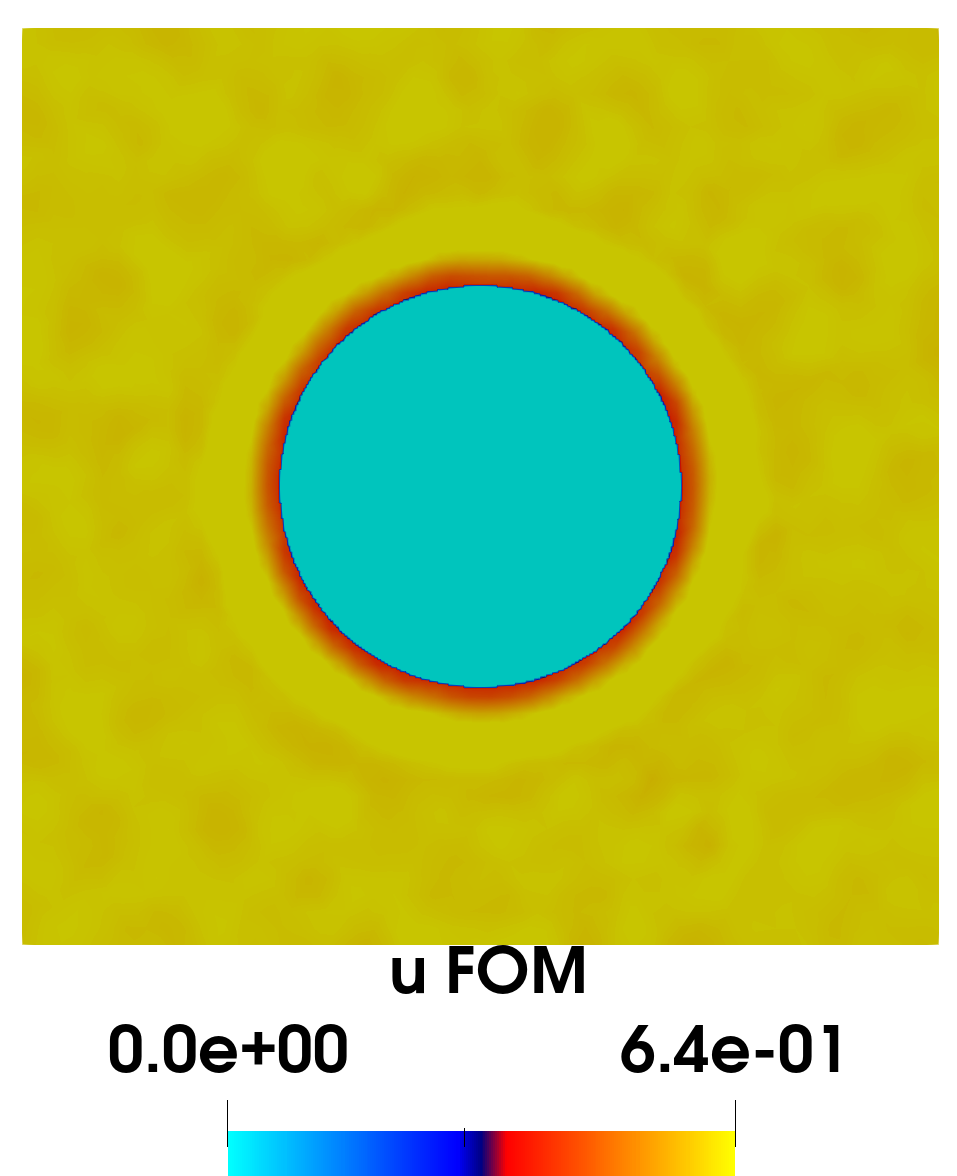}
\end{minipage}
 \qquad 
 \begin{minipage}{0.25\textwidth}
  \includegraphics[width=\textwidth]{{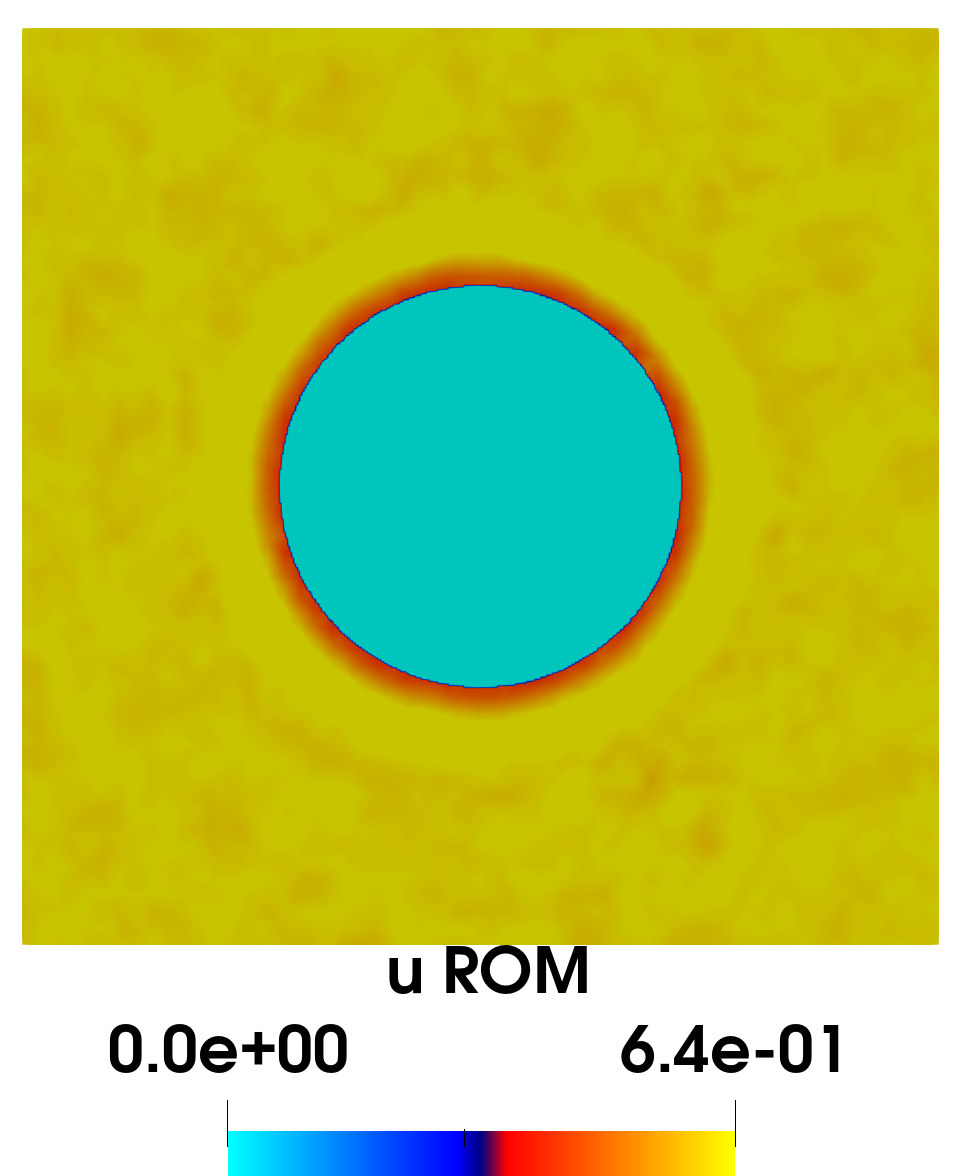}} 
\end{minipage}
 \qquad 
 \begin{minipage}{0.25\textwidth}
  \includegraphics[width=\textwidth]{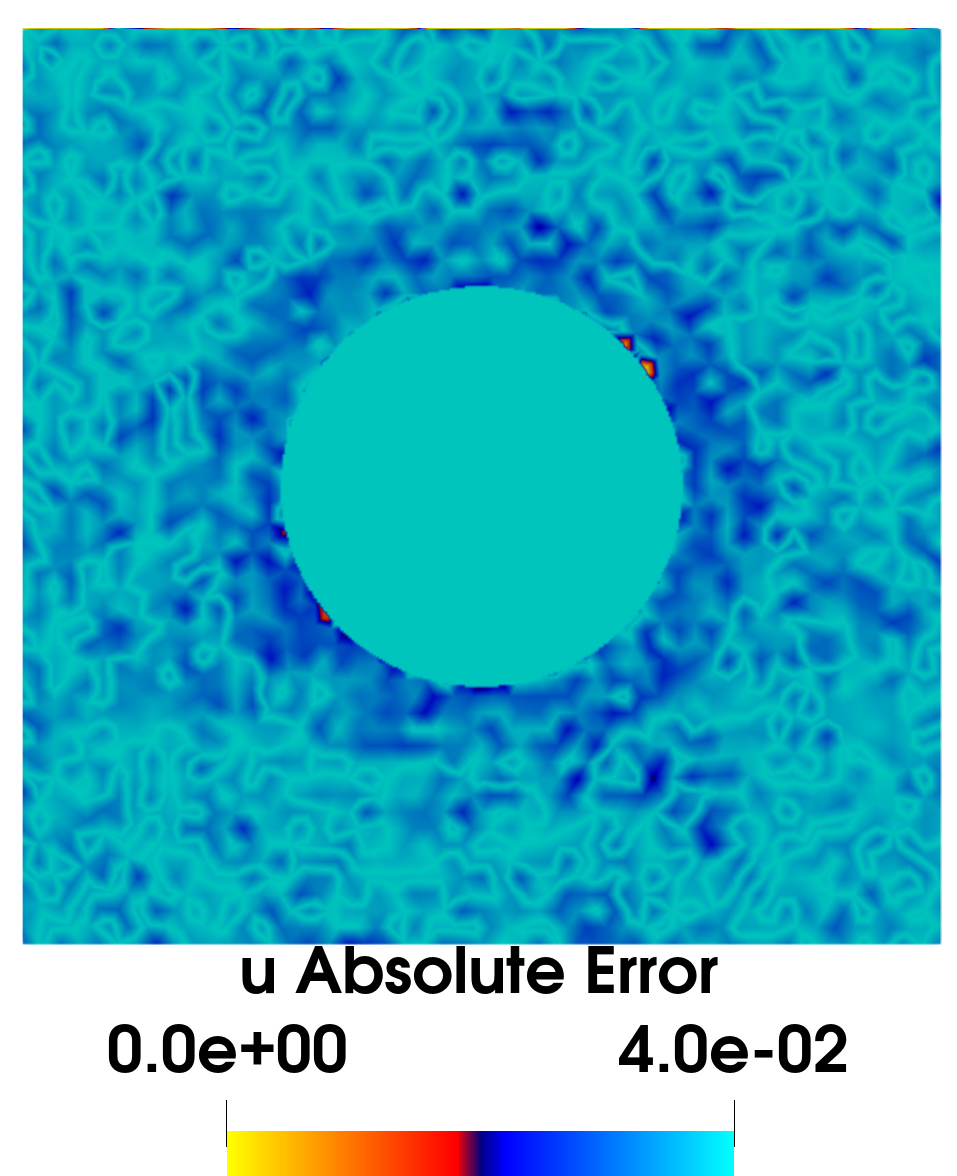}
\end{minipage}
\begin{minipage}{0.25\textwidth}
  \includegraphics[width=\textwidth]{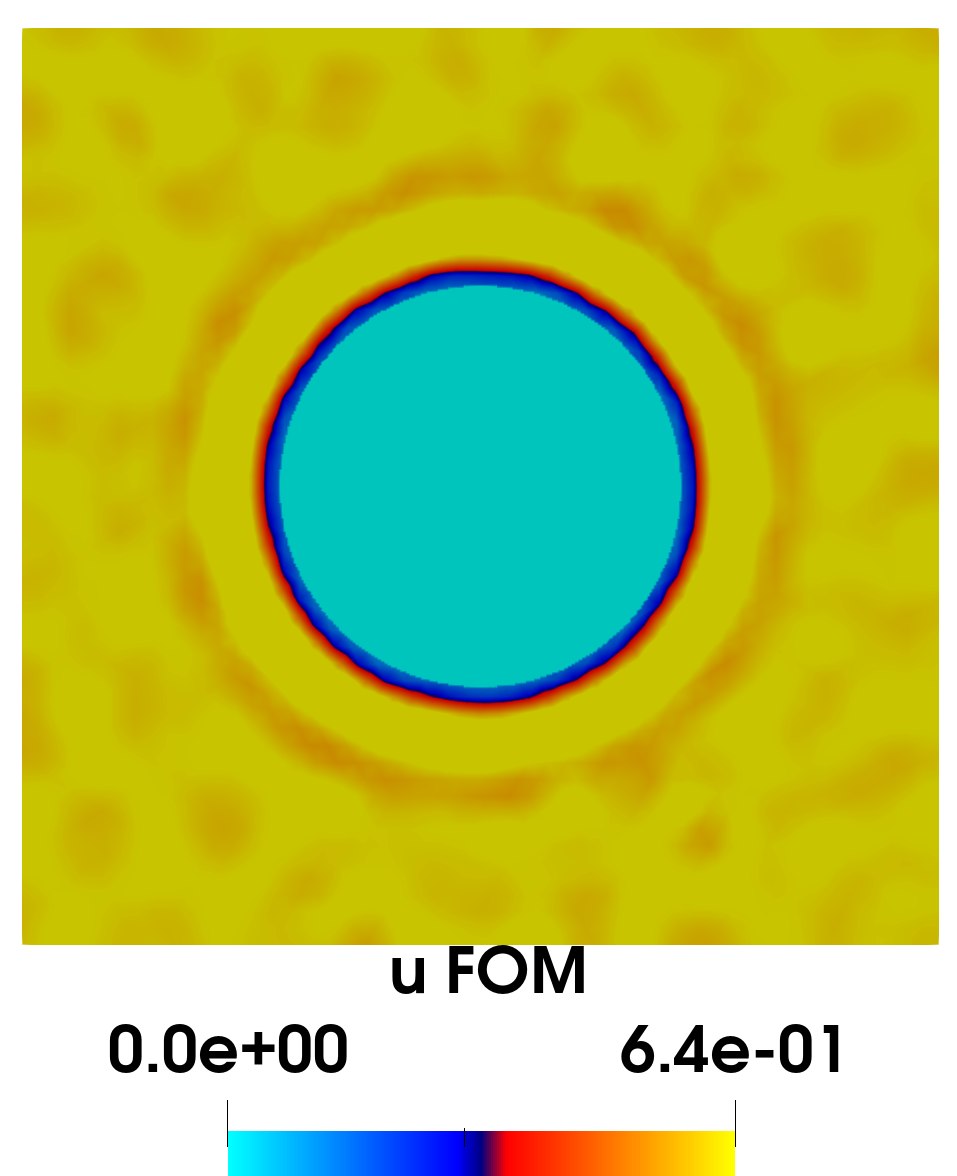}
\end{minipage}
\qquad
\begin{minipage}{0.25\textwidth}
  \includegraphics[width=\textwidth]{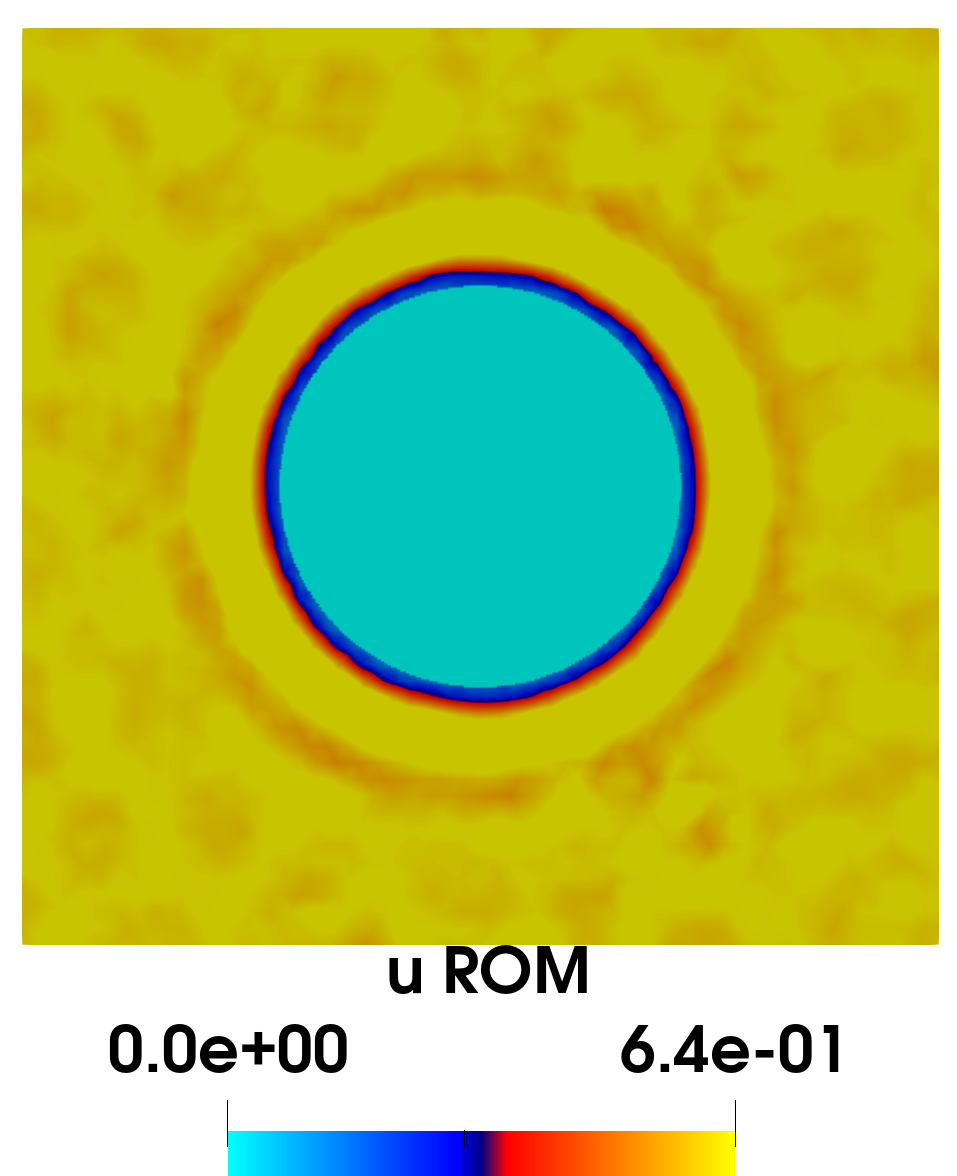} 
\end{minipage}
\qquad
\begin{minipage}{0.25\textwidth}
  \includegraphics[width=\textwidth]{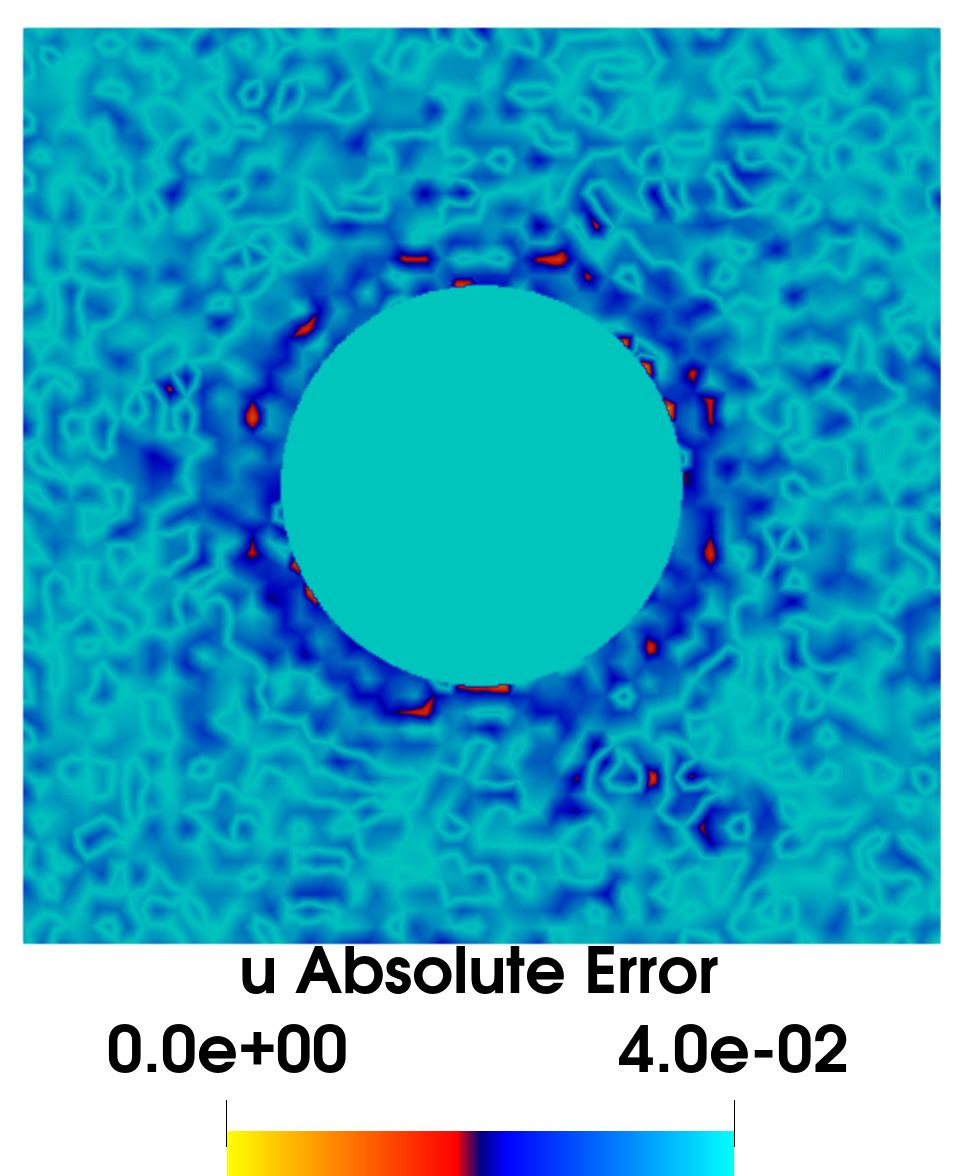}
\end{minipage}
\begin{minipage}{0.25\textwidth}
  \includegraphics[width=\textwidth]{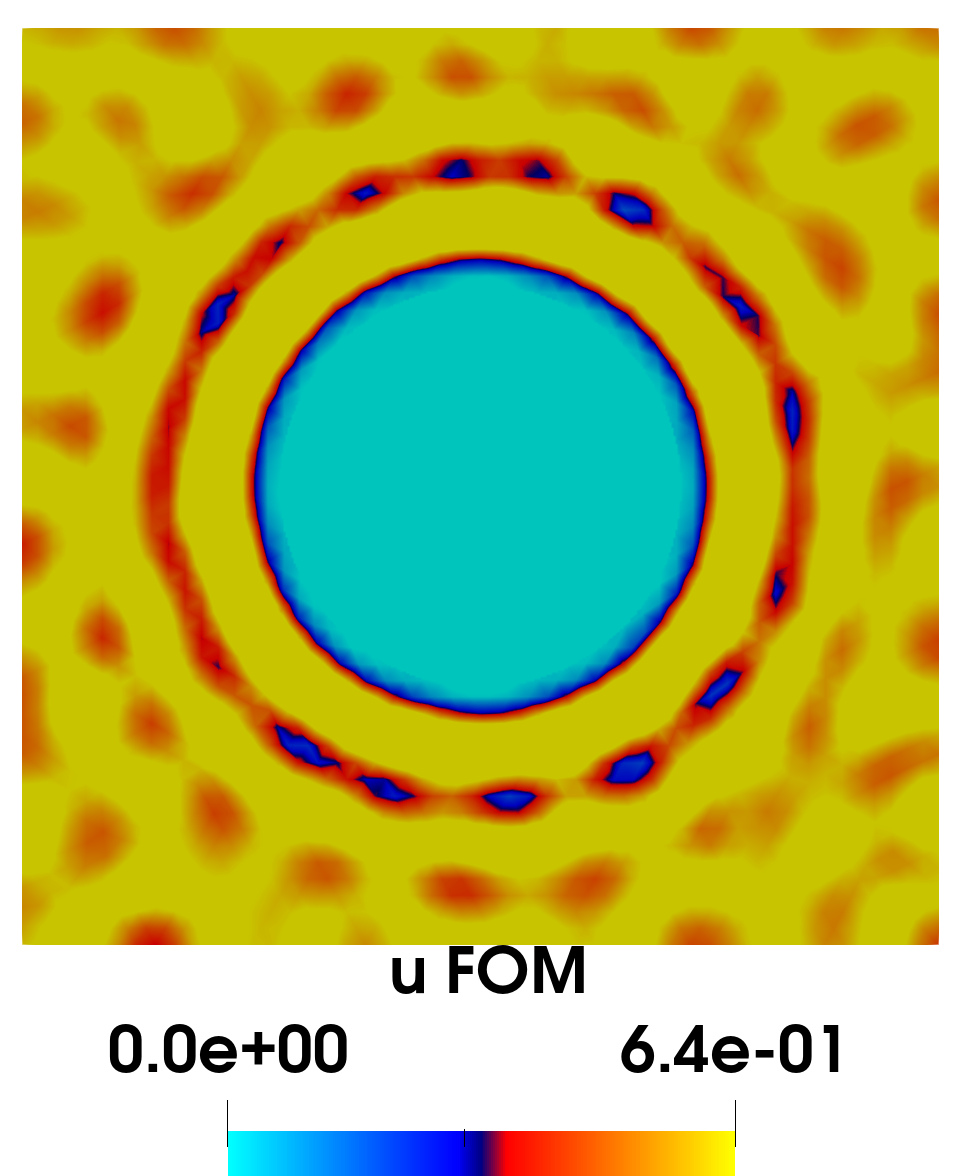}
\end{minipage}
\qquad
\begin{minipage}{0.25\textwidth}
  \includegraphics[width=\textwidth]{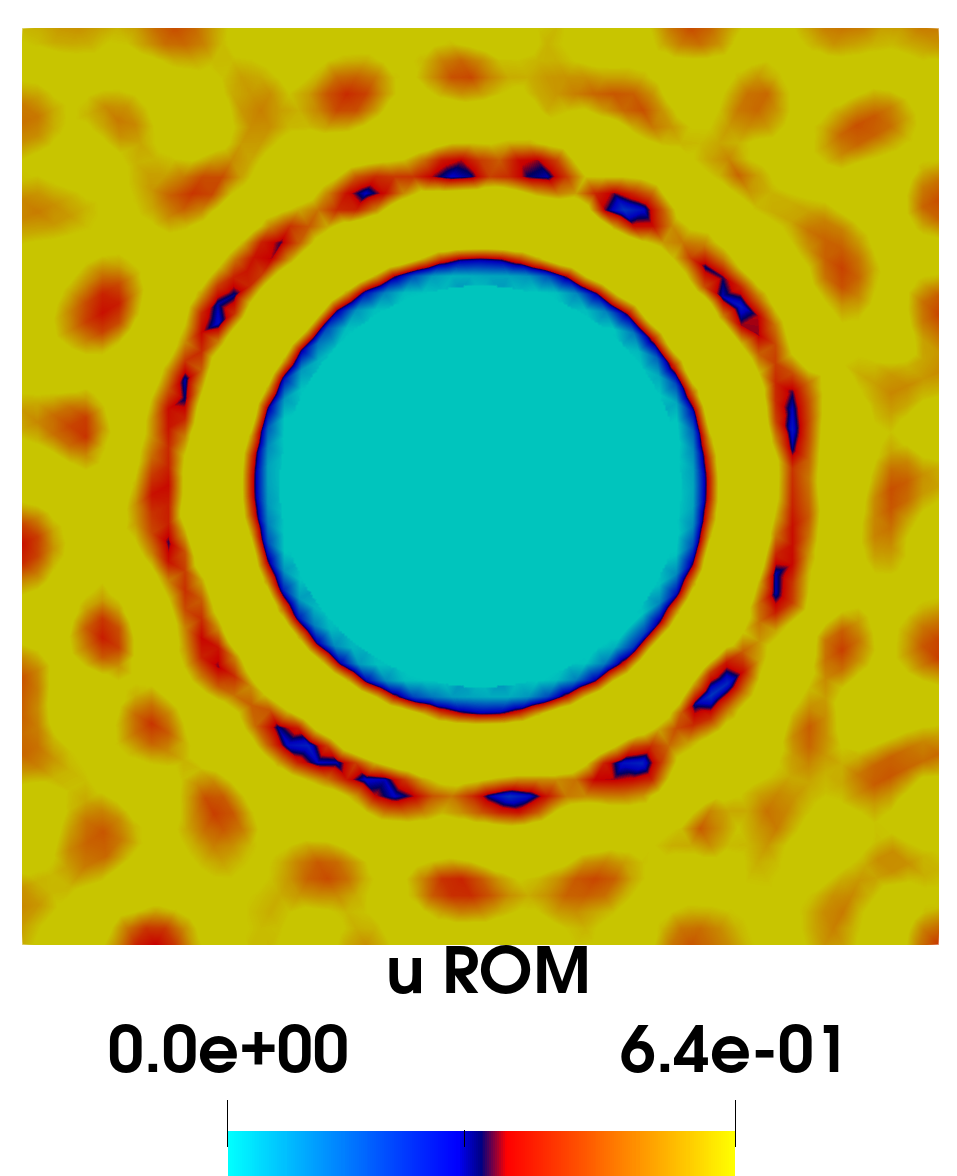}
\end{minipage}
\qquad
\begin{minipage}{0.25\textwidth}
  \includegraphics[width=\textwidth]{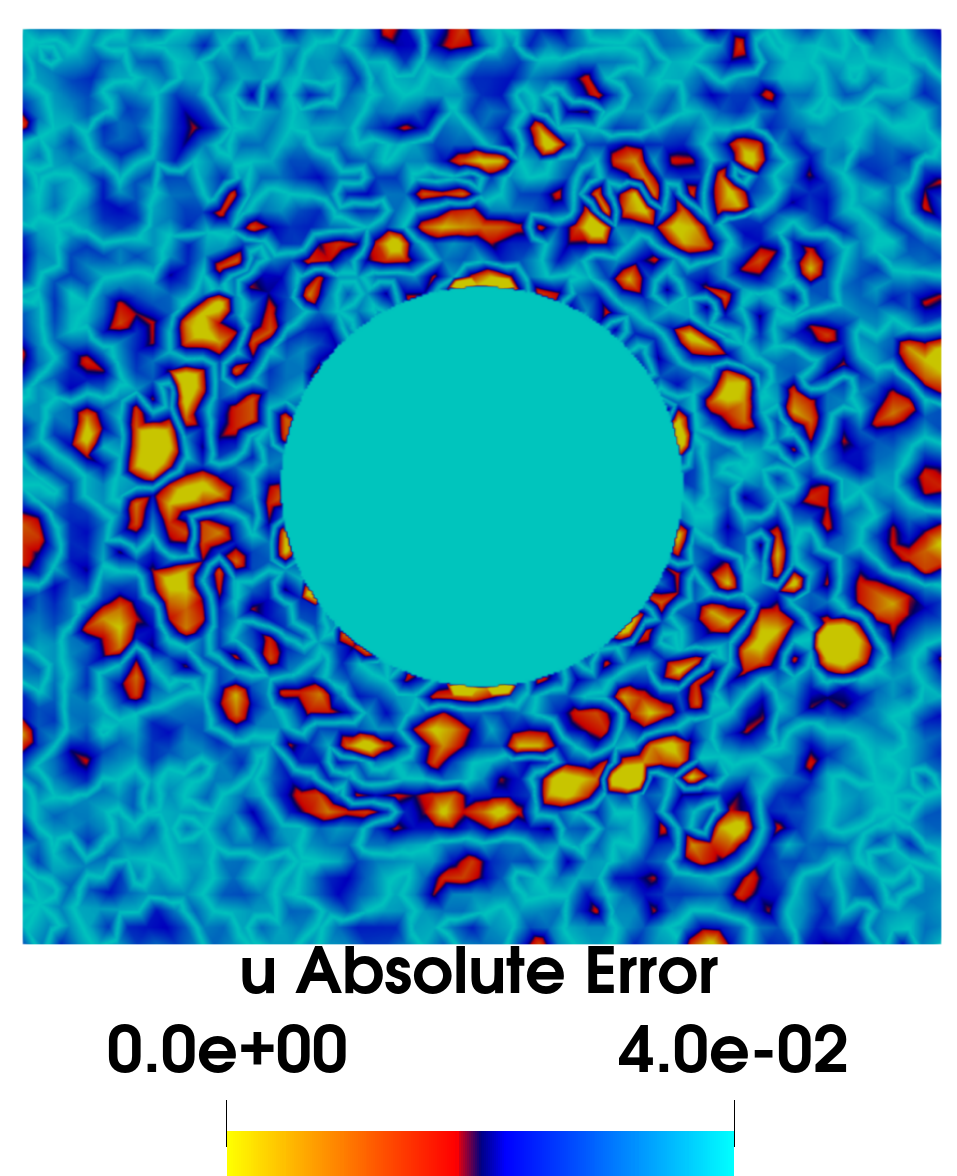}
\end{minipage}
\begin{minipage}{0.25\textwidth}
  \includegraphics[width=\textwidth]{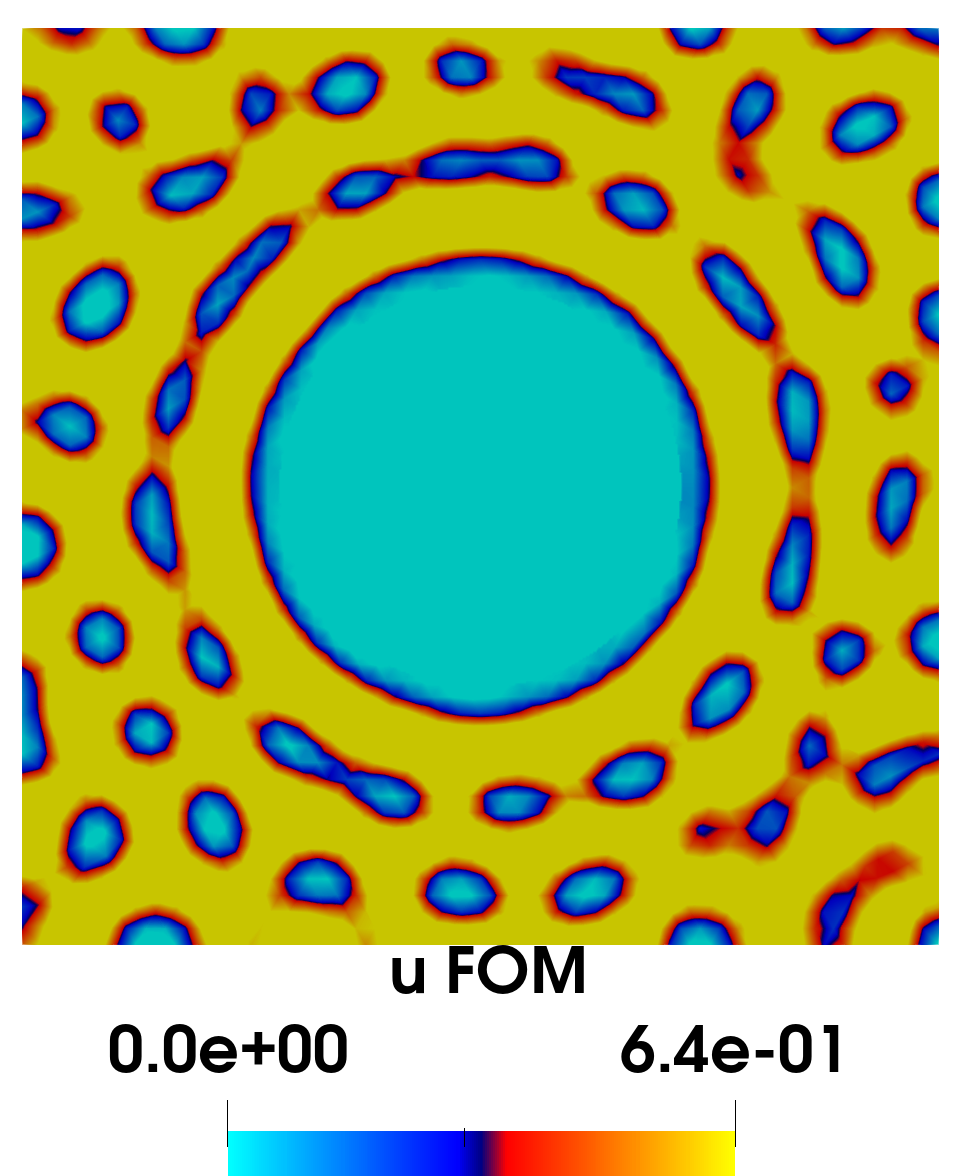}
\end{minipage}
\qquad
\begin{minipage}{0.25\textwidth}
  \includegraphics[width=\textwidth]{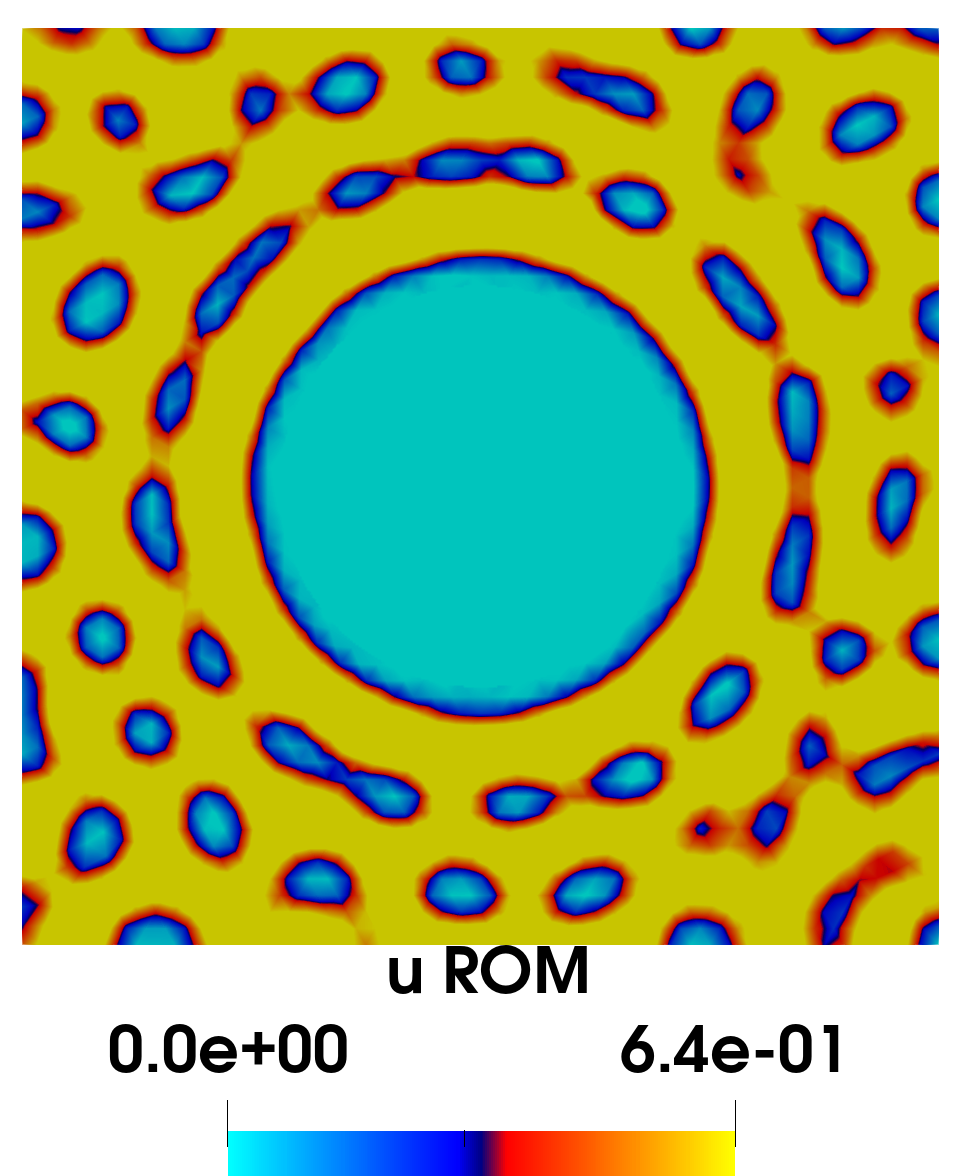}
\end{minipage}
\qquad
\begin{minipage}{0.25\textwidth}
  \includegraphics[width=\textwidth]{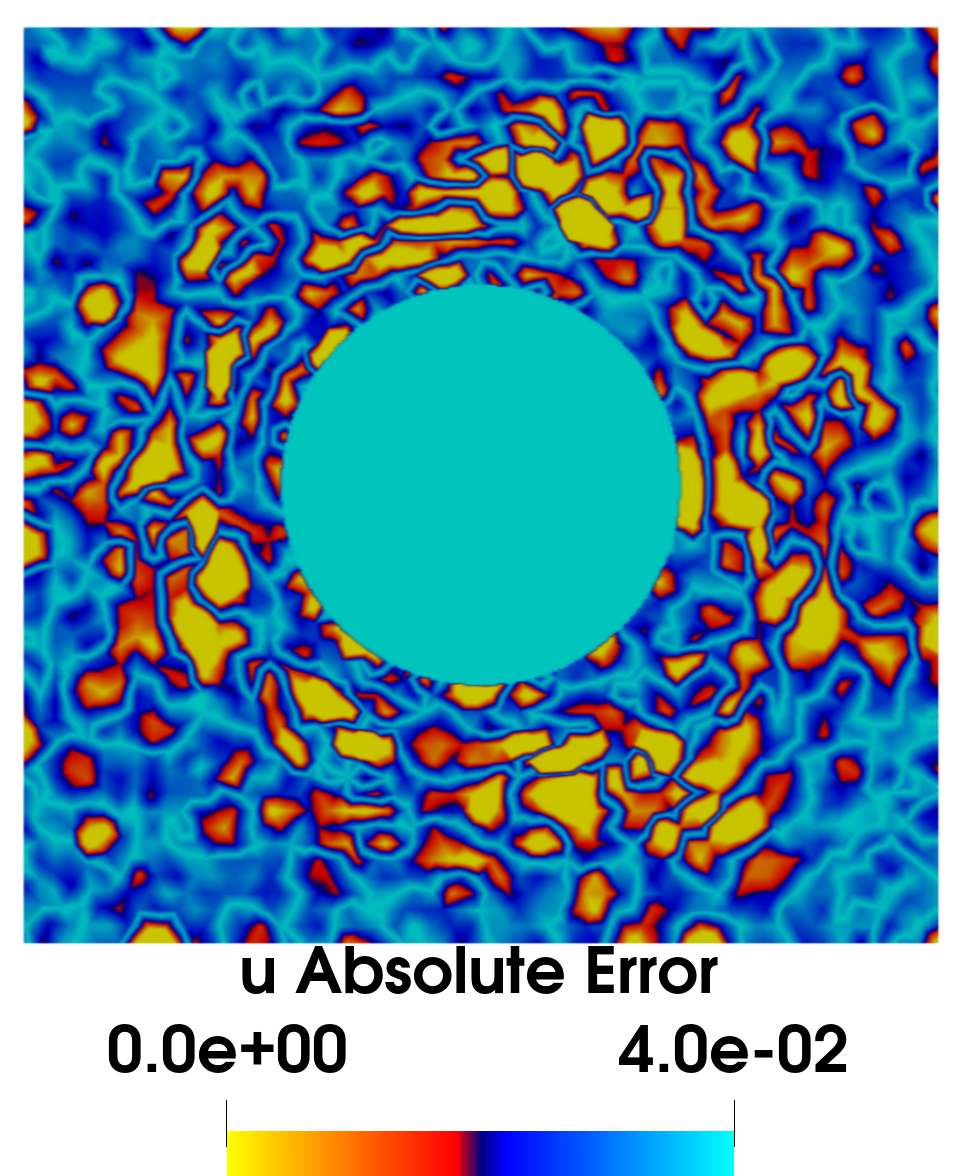}
  \end{minipage}
  \begin{minipage}{0.25\textwidth}
  \includegraphics[width=\textwidth]{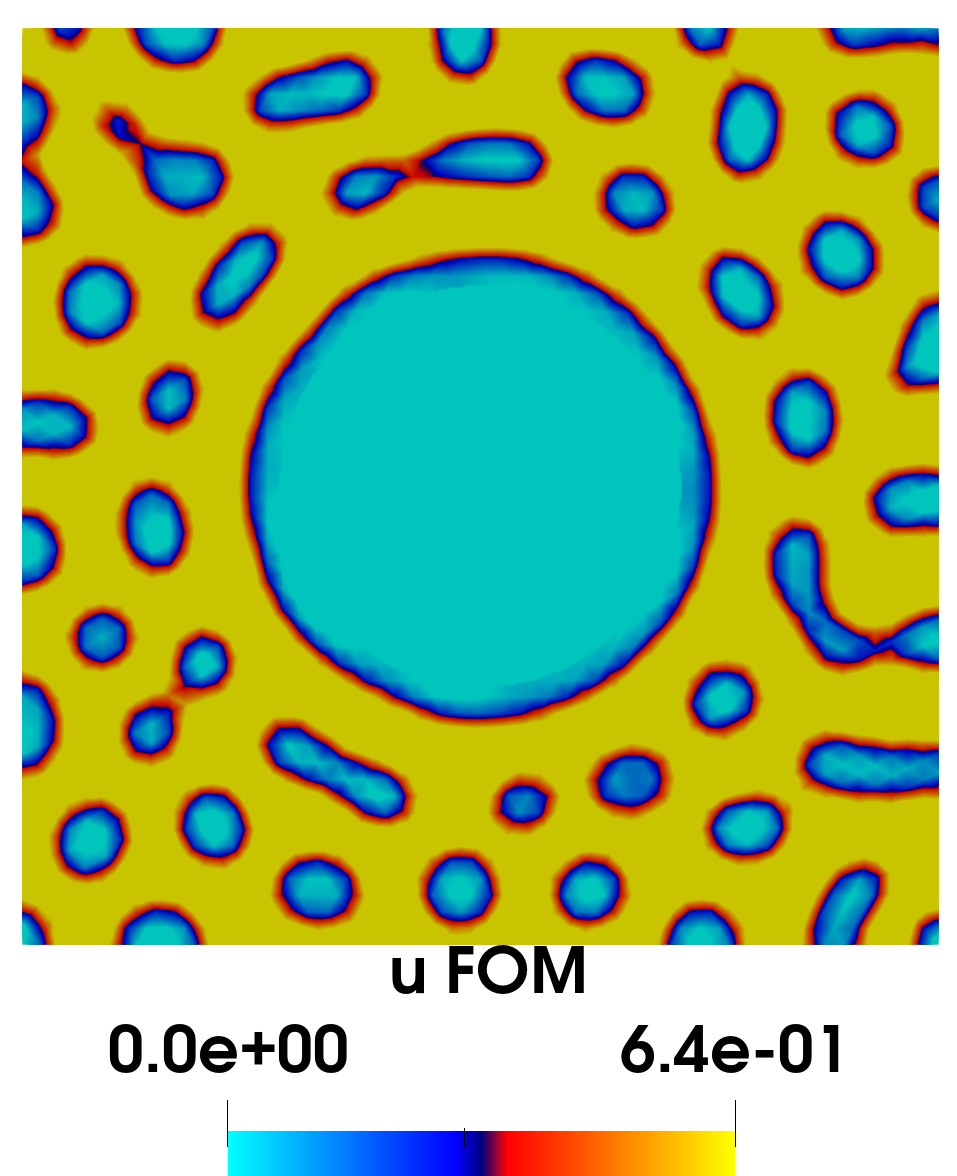}
\end{minipage}
\qquad
\begin{minipage}{0.25\textwidth}
  \includegraphics[width=\textwidth]{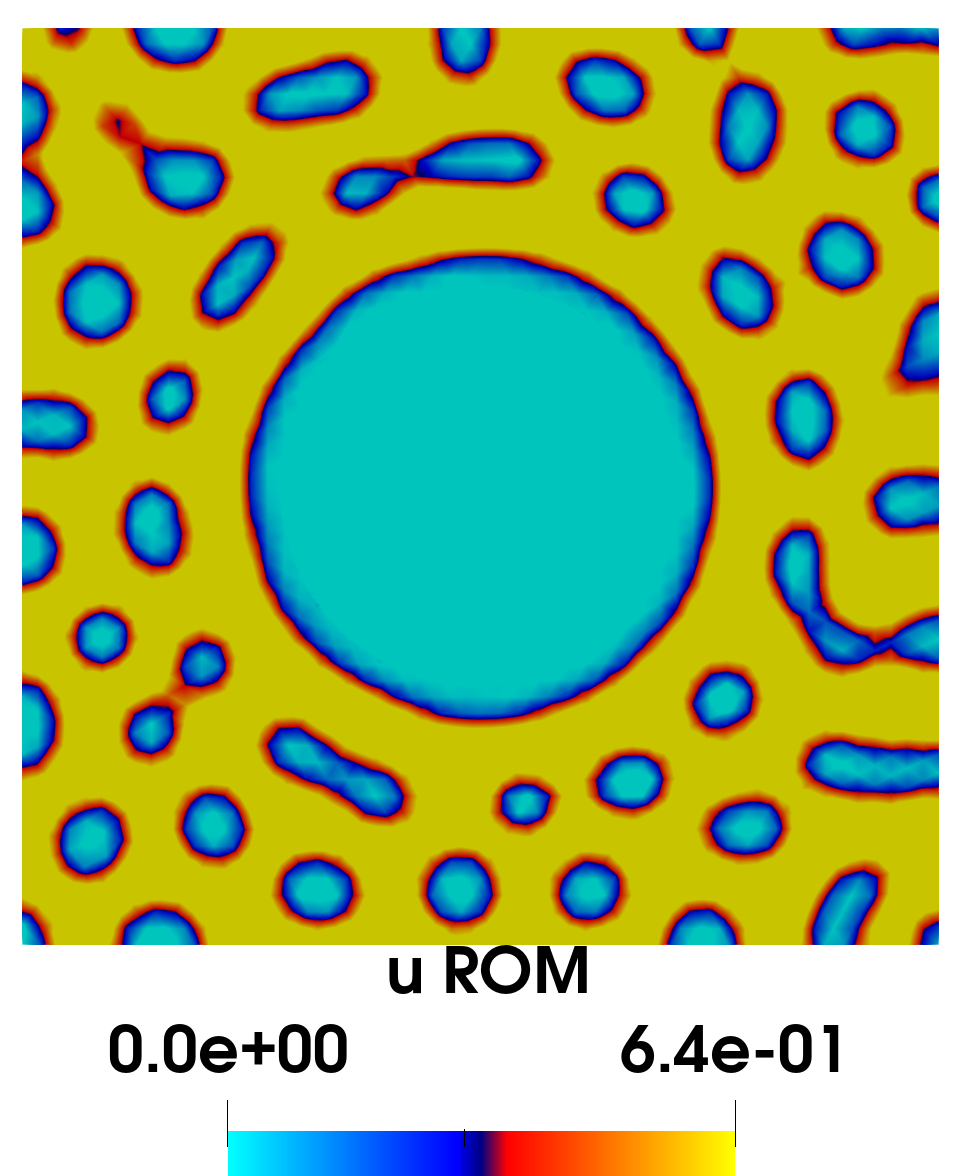}
\end{minipage}
\qquad
\begin{minipage}{0.25\textwidth}
  \includegraphics[width=\textwidth]{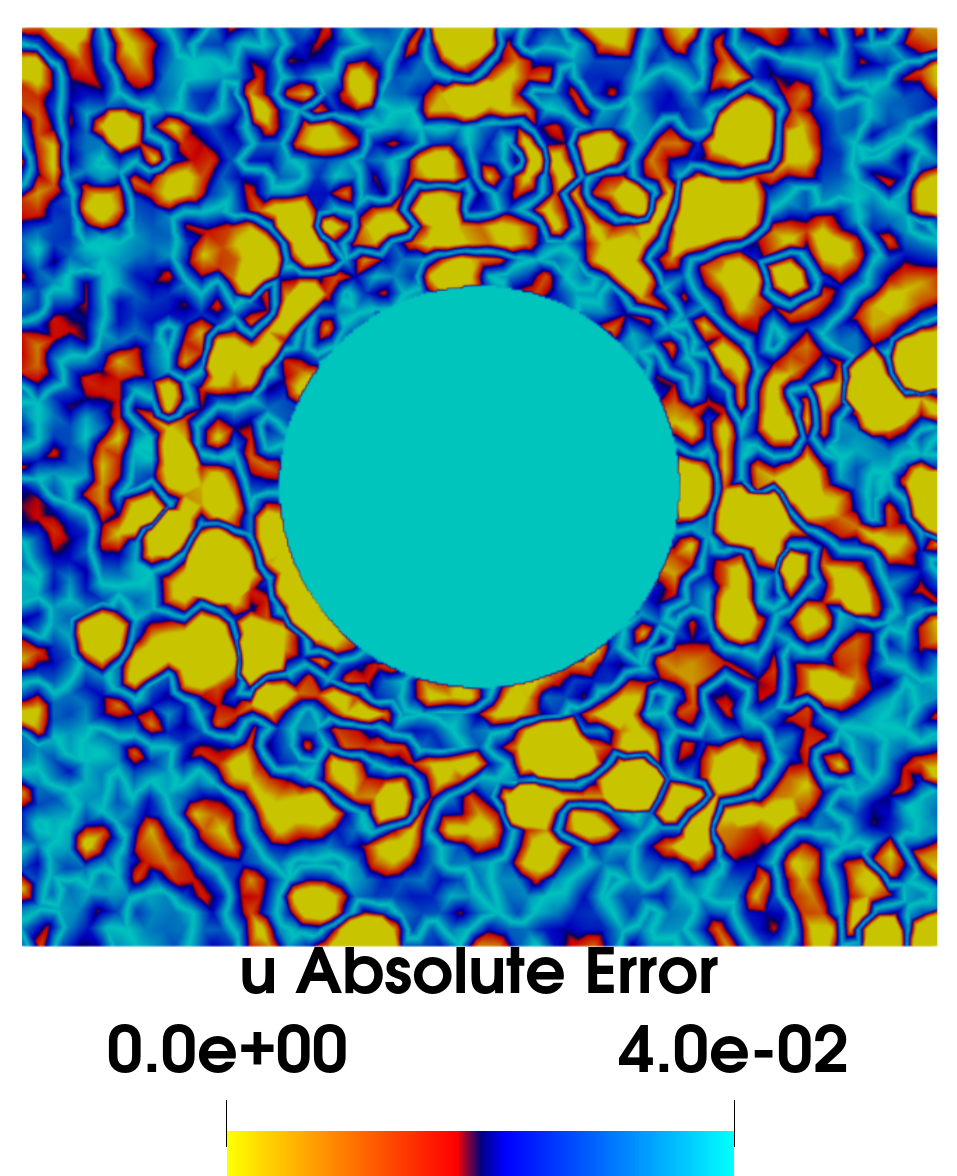}
\end{minipage} 
\end{minipage}
\caption{{\blue{Results for the geometrical parametrized embedded circle  with diameter $\mu_{\text{test}} = 0.43774
$ 
and {\textit{Dirichlet embedded boundary}}. In the first, second and third column we report the {\magenta{full order}} solution, the reduced order solution and the {\blue{absolute}} error plots for the concentration field. Each row corresponds to a different time $t=[10,20,40,60,100]\tau_n$.}}
} \label{Dirichlet_FULL_RED_ERROR_1P}
\end{figure}
\begin{table}\centering
  \begin{tabular}{|c||c||c|c|c|}
   \hline
    Snapshots for:& \multicolumn{4}{c|}{900 train parameters} \\ \hline \hline
    Modes     & \multicolumn{2}{c|}{Relative error for Neumann case}& \multicolumn{2}{c|}{Relative error for Dirichlet case} \\ \hline
    ({\magenta{$N_{POD}$}})       & Concentration $u$  & Potential $w$& Concentration $u$  & Potential $w$\\
    \hline
    1  & 0.18370	& 0.49895 & 0.22298 & 0.96624 \\
    5  & 0.15292	& 0.41387 & 0.20528 & 0.37814 \\
	10 & 0.06069	& 0.18397 & 0.11857 & 0.28665 \\
	15 & 0.04086	& 0.11956 & 0.08977 & 0.21658 \\
	20 & 0.04046	& 0.11311 & 0.07631 & 0.18301 \\
	25 & 0.03498	& 0.09183 & 0.04808 & 0.13259 \\
	30 & 0.02883	& 0.07852 & 0.03932 & 0.10802 \\
	35 & 0.02641	& 0.07085 & 0.03284 & 0.09350 \\
	40 & 0.02063	& 0.05981 & 0.02879 & 0.08377 \\
    45 & 0.01987	& 0.05681 & 0.02749 & 0.08072 \\
    \hline
  \end{tabular} 
      \caption{Geometrical parametrization and the {\magenta{mean relative error (for 30 tests)}} between the {\magenta{full order}} and the reduced basis solution for the concentration and potential component. Results are reported for different dimensions of the reduced basis spaces {\magenta{while the ROM has been trained onto 900 parameter samples and tested onto 30 samples 
      chosen randomly inside the parameter space.
}}}
%
%
\label{table:r_errors}
\end{table}
\begin{figure} \label{fig:comp_time}
\centering
\begin{minipage}{\textwidth}
\centering
\hskip1pt
(i)\includegraphics[width=0.42\textwidth]{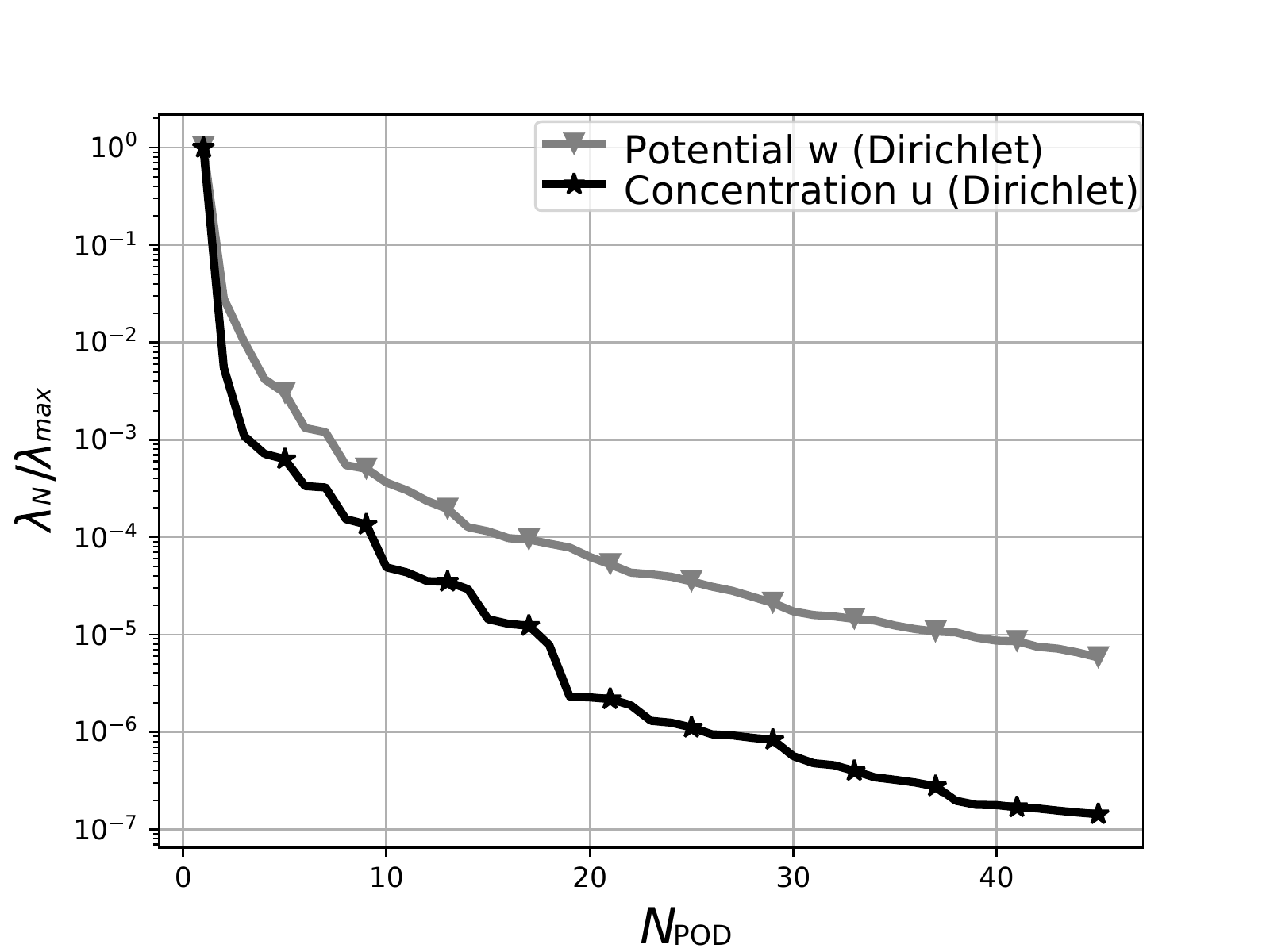}
~~~~~~
 \hskip-00pt \includegraphics[width=0.42\textwidth]{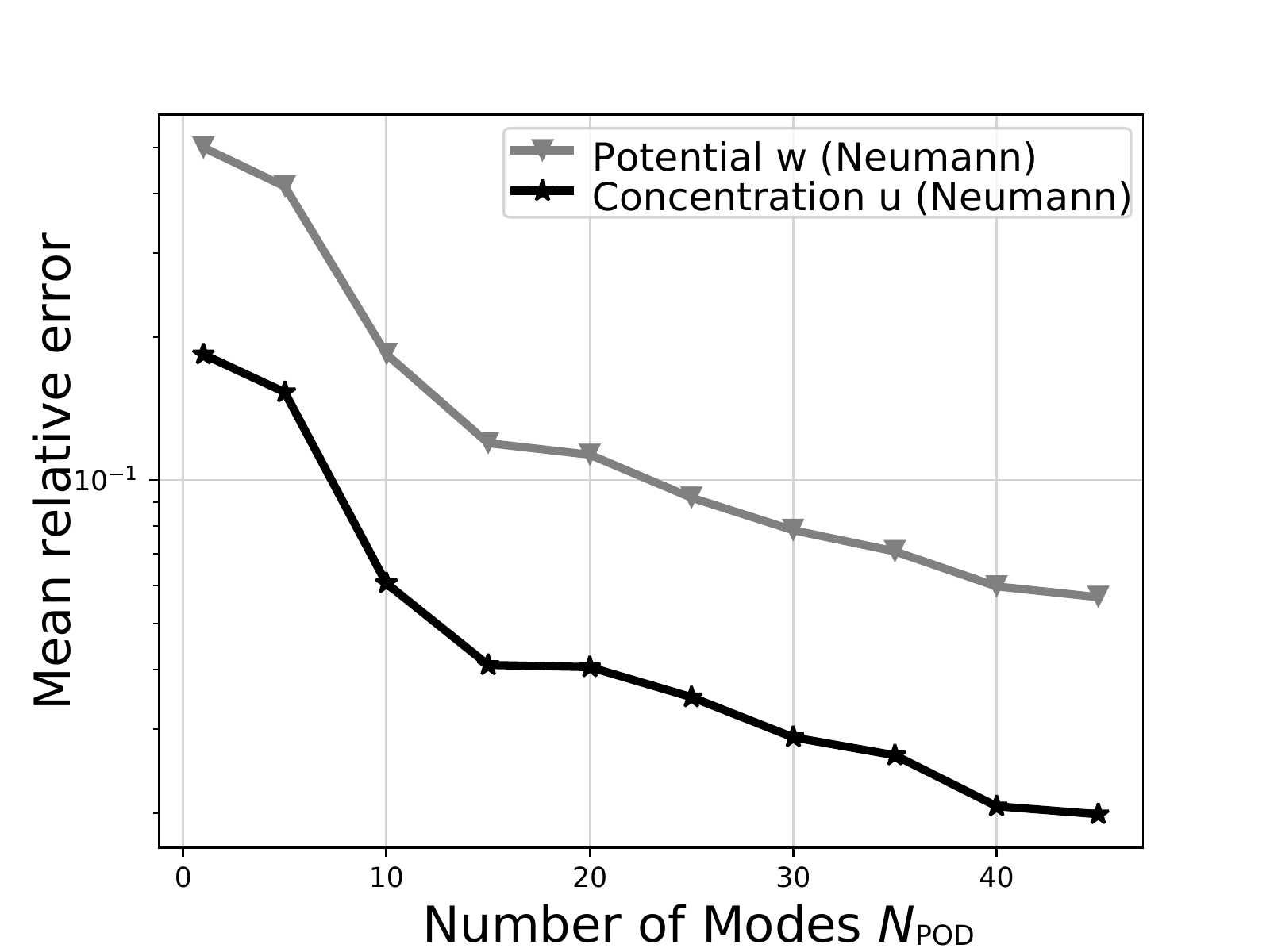}
(ii)
~~~~~~
\end{minipage}
\begin{minipage}{\textwidth}
\centering
(iii)\includegraphics[width=0.42\textwidth]{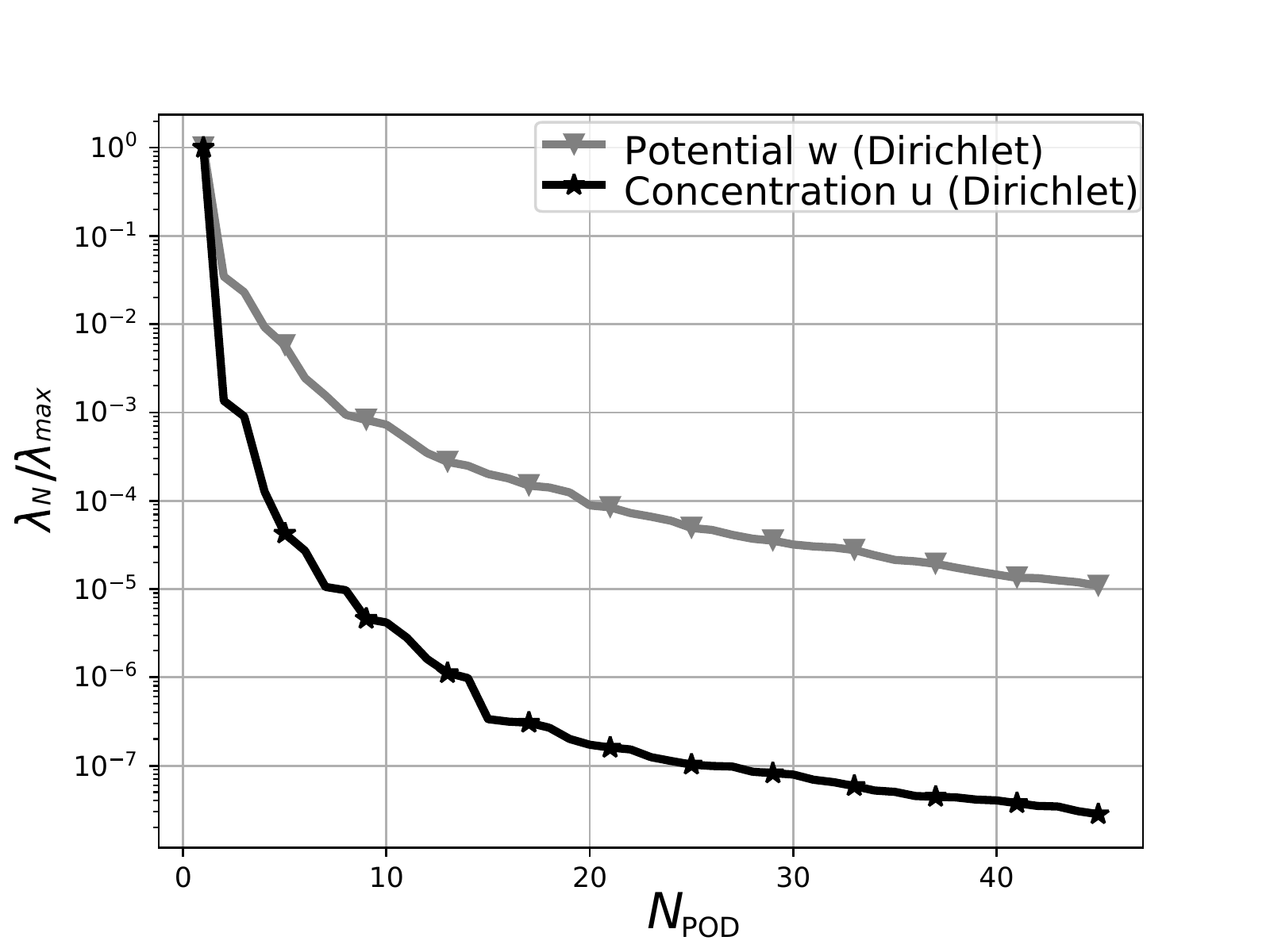}
~~~~~~
 \hskip-00pt \includegraphics[width=0.42\textwidth]{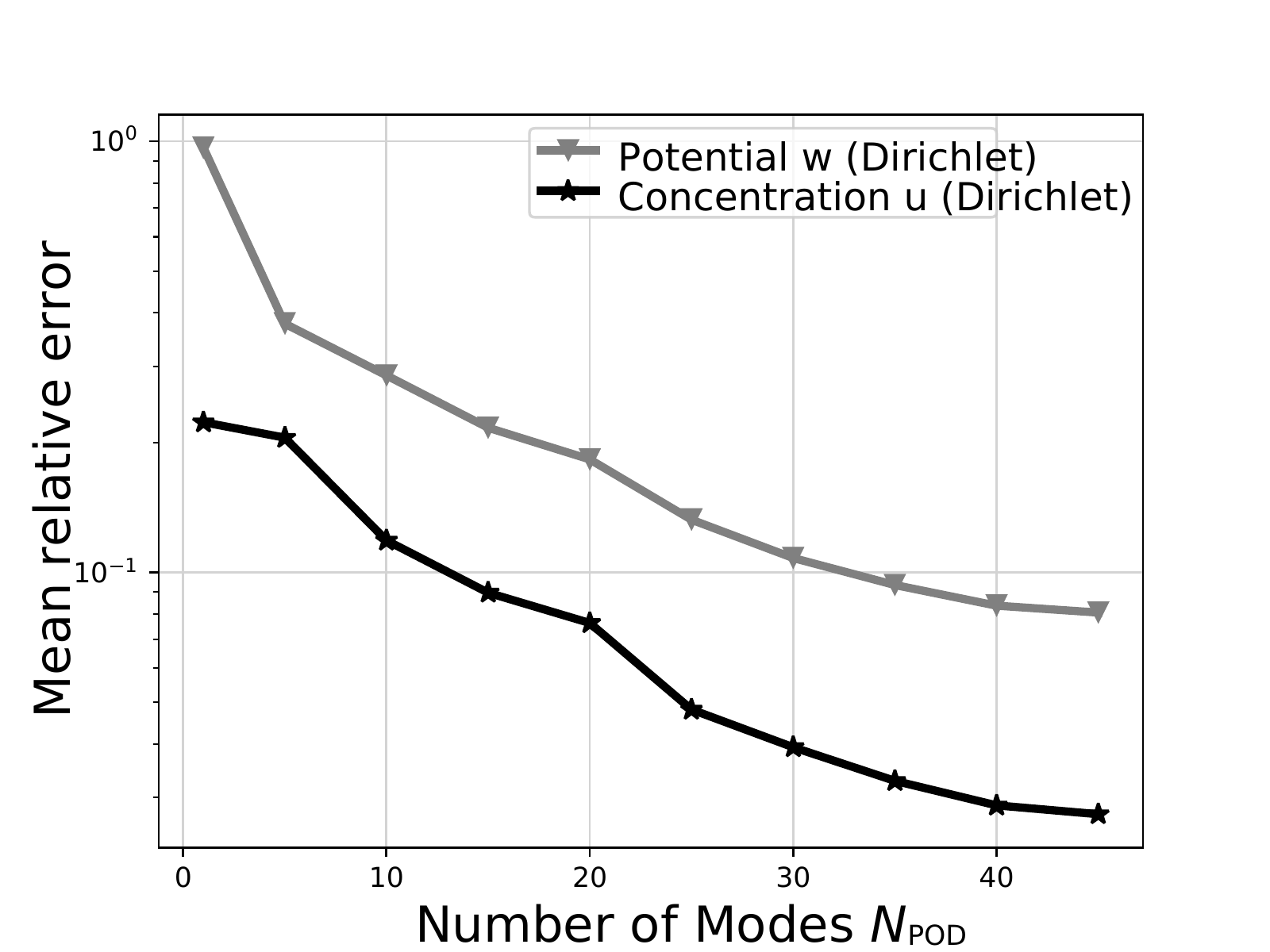}
(iv)
~~~~~~
\end{minipage}
  \caption{Geometrical parametrization and visualization of the results. The left plots depict the eigenvalues decay with respect to the numbers of the modes for the concentration and potential variable. On the right plots, we visualize {\magenta{the mean value}} of the  relative errors of the reduced-order problem for {\magenta{thirty}} random parameter values and for various number of modes.}
%
  %
\label{fig:r_errors}
\end{figure}
\begin{table}\centering
 \begin{tabular}{|c||c|c|}
   \hline
    Snapshots for:& \multicolumn{2}{c|}{900 train parameters}   \\ \hline \hline
    Modes     & {execution times} & Savings  \\
    (${\magenta{N_{POD}}}$)       & (t {\magenta{seconds}}) & $ (t_{\text{FOM}} - t_{\text{RB}})/t_{\text{FOM}}$ 
    \\
    \hline
     1 
     & 0.02808
     & 99.718\%
      \\
     10 
     & 0.05299 
     & 99.468\%
     \\
     20 
     & 0.08422 
     & 99.154\%
     \\
     30 
     & 0.10766 
     & 98,919\% 
     \\
     40 
     & 0.14128  
     & 98.581\%\\
     45 
     & 0.16291
     & 98.364\%
     \\
    \hline
  \end{tabular} 
  \caption{Execution time at the reduced order level {\blue{and per cent (\%) savings}}
. The computational time includes the projection of the full order matrices, the execution time of the online solver and the {\magenta{solution}} of the reduced problem. {{Times are for the solution of one random value of the input parameter. The time execution at full order method level (FOM) is {\magenta{approximate}} to $\approx 9.9623$ sec.}}}
\label{table:timers}
\end{table}
\begin{figure} \label{fig:comp_time} \centering
\begin{minipage}{\textwidth}
\centering
(i)\includegraphics[width=0.42\textwidth]{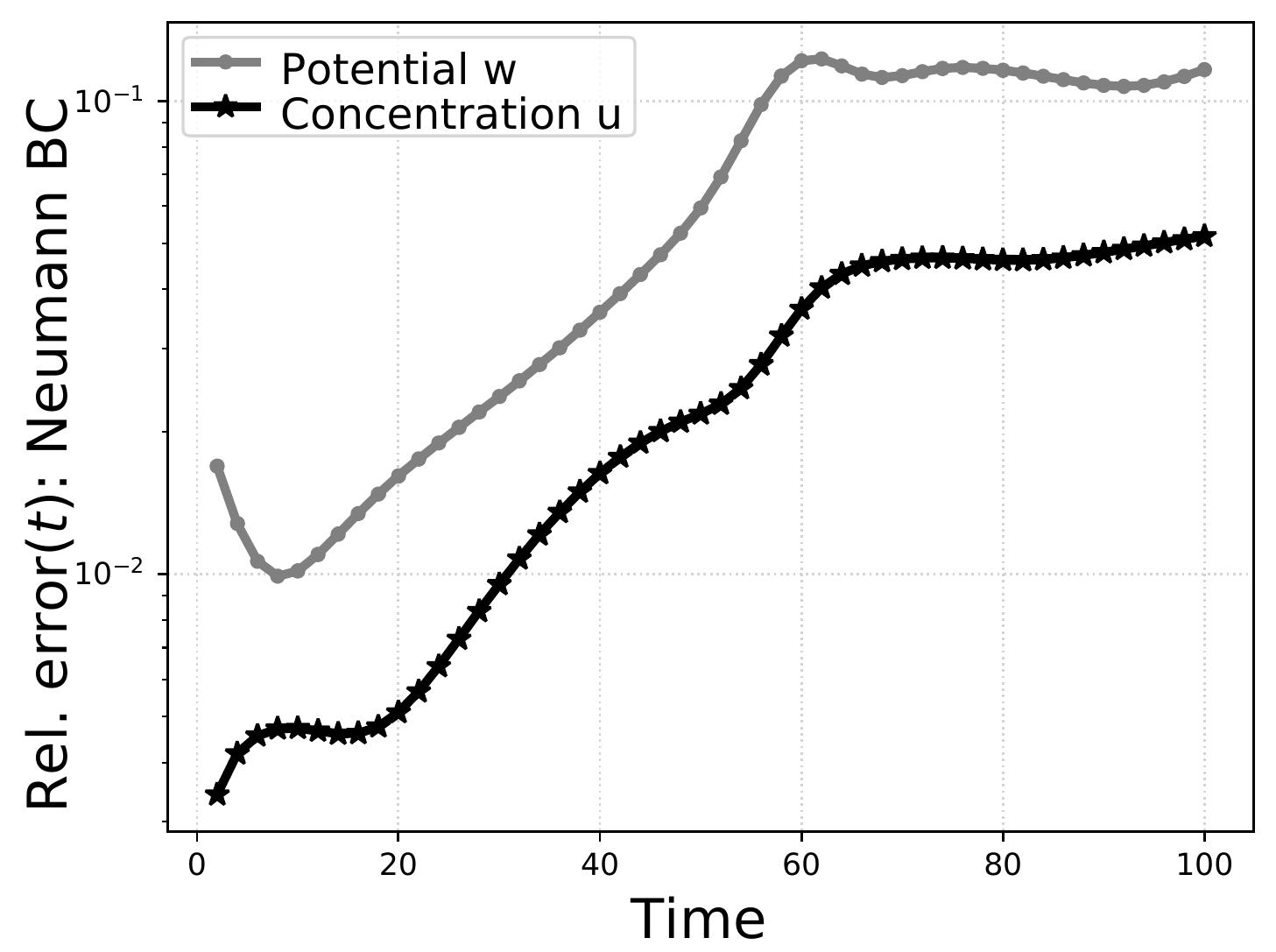}
~~~~~~
\includegraphics[width=0.42\textwidth]{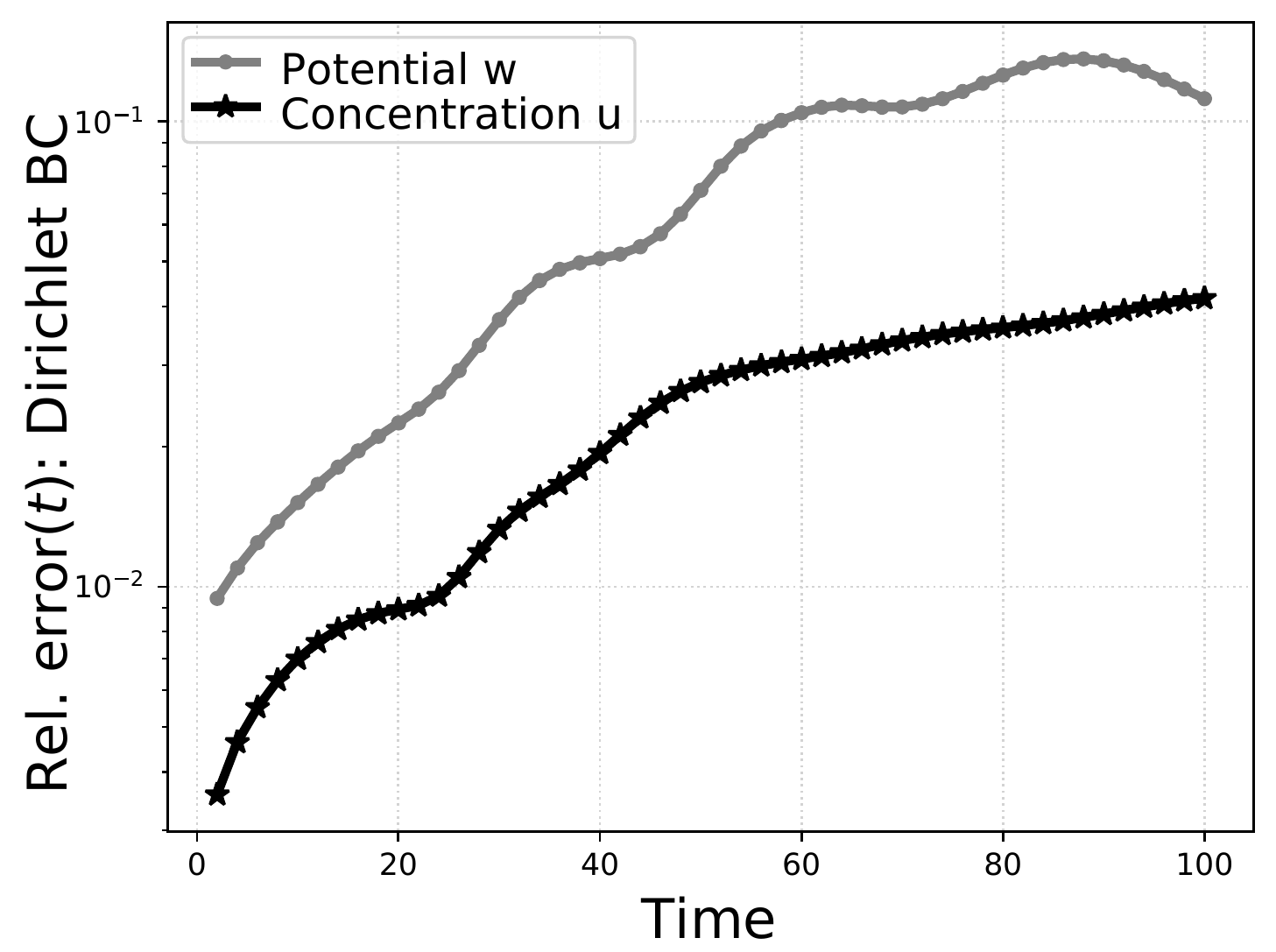}{\blue{(iii)}}
\\
(ii)\includegraphics[width=0.42\textwidth]{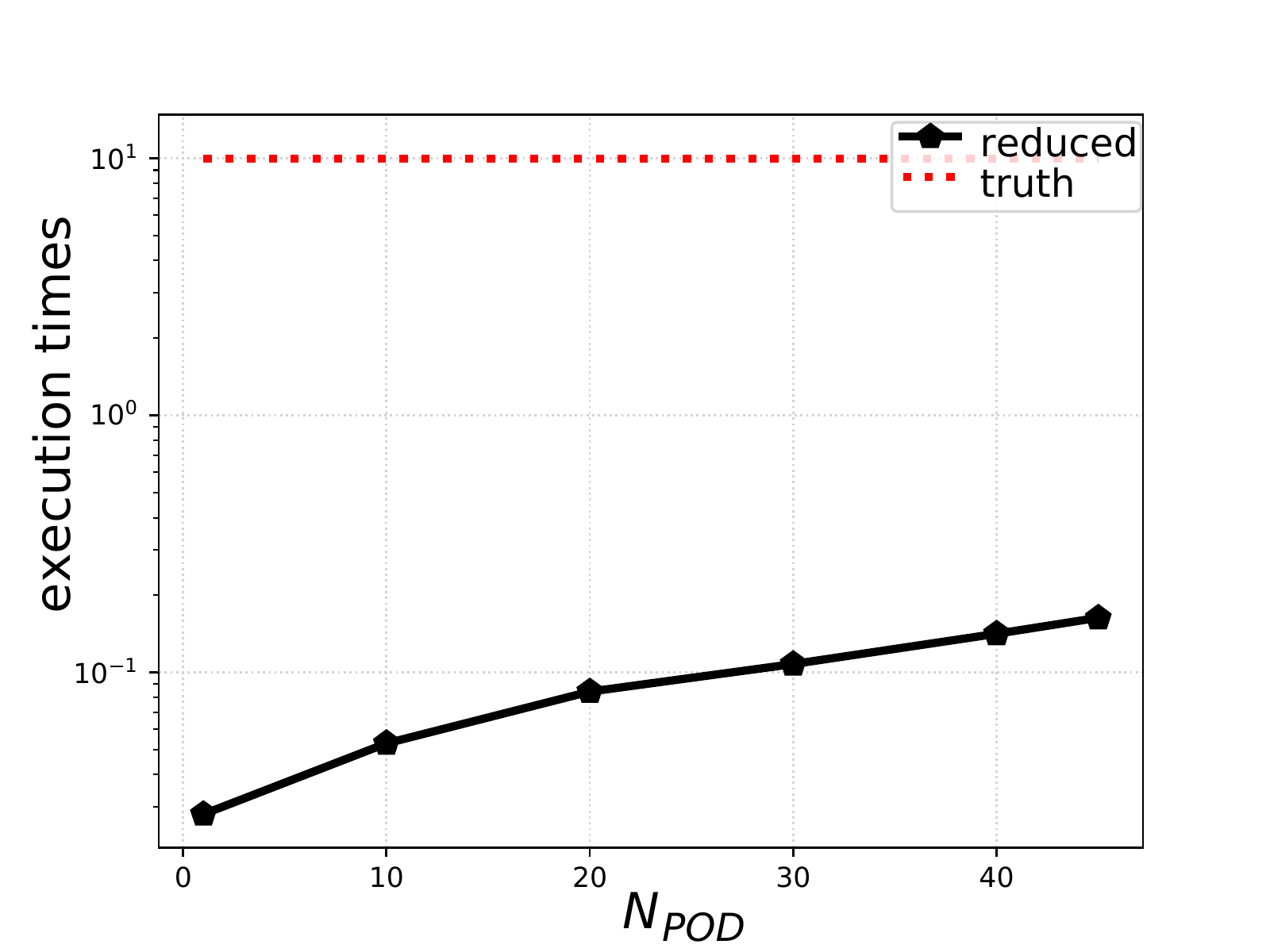}
~~~~~~
\end{minipage}
  \caption{
The {\magenta{(i)}} and {\magenta{(iii)}} plots depict the relative errors with respect {\magenta{to}} time for the concentration and potential variable and 45 modes. {\Xgreen{In the {{(ii)}} plot}}, we visualize the execution times of the reduced-order problem for one random parameter value and for various {\magenta{numbers}} of modes.}
\label{fig:timers}
\end{figure}
\begin{table} \centering
 \begin{tabular}{|c||c|c|c|c|}
   \hline
    Snapshots for:& \multicolumn{4}{c|}{900 train parameters} \\ \hline \hline
    time ($i\times \tau_n$)     & \multicolumn{2}{c|}{Relative error (Neumann)} & \multicolumn{2}{c|}{{\blue{Relative error (Dirichlet)}}}\\ \hline
    $i$       & Concentration u & Potential w & {\blue{Concentration u}} & {\blue{Potential w}}\\
    \hline
     1   & 0.00342 & 0.01692 & 0.00357  & 0.00944 \\
     10  & 0.00472 & 0.01016 & 0.00699  & 0.01517 \\
     20  & 0.00510 & 0.01613 & 0.00894  & 0.02250 \\
     30  & 0.00952 & 0.02374 & 0.01330  & 0.03751 \\
     40  & 0.01634 & 0.03580 & 0.01938  & 0.05078 \\
     50  & 0.02184 & 0.05944 & 0.02752  & 0.07117 \\
     60  & 0.03639 & 0.12172 & 0.03086  & 0.10452 \\
     70  & 0.04636 & 0.11315 & 0.03381  & 0.10750 \\
     80  & 0.04611 & 0.11615 & 0.03606  & 0.12585 \\
     90  & 0.04787 & 0.10791 & 0.03856  & 0.13514 \\                                           
     100 & 0.05181 & 0.11652 & 0.04169  & 0.11201 \\
    \hline
  \end{tabular} 
  \caption{The $L^2$ relative {\magenta{error}} is reported {\magenta{over time}} for the concentration and potential field for the Neumann and Dirichlet type of boundaries experiments. The ROMs have been obtained with 45 modes for the concentration {\magenta{and potential variables}} for both cases.}
\label{table:r_error_t}
\end{table}

As it is displayed in Figures \ref{FULL_RED_ERROR_1P} and \ref{Dirichlet_FULL_RED_ERROR_1P} {\blue{--Neumann and Dirichlet embedded boundary experiments respectively--
and with a minds eye comparison of the FOM and ROM solutions in first and second column, with a brief look they seem identical for every row associated with the time instances $t=[10,20,40,60,100]\tau_n$. 
In the third column and looking from a more detailed point of view, the absolute error, {\magenta{$|u-u_r|$}},  at each point of the geometry domain, {\bblue{using $45$ modes,}}  is visualized 
 for the {\bblue{zero Neumann --on all sub-boundaries--}} and the {\bblue{zero}} Dirichlet {\bblue{--only the embedded is Dirichlet while the remaining sub-boundaries are zero Neumann--} boundary conditions cases.}}  
In Table \ref{table:r_errors} and again for both types of boundaries, 
the {\magenta{mean relative errors $||u-u_r||_{L^2(\Omega)}/ ||u||_{L^2(\Omega)}$ and  $||w-w_r||_{L^2(\Omega)}/ ||w||_{L^2(\Omega)}$,}}
for {\magenta{the phase}} field for various number of basis {\magenta{functions}}  are reported and their graph can be seen in Figure \ref{fig:r_errors} as well as the normalized eigenvalues and their decay which have been used for the ROM. Thereafter, for a better understanding, the relative errors evolution 
with respect to time, for a  $45$ modes test, is demonstrated in Table \ref{table:r_error_t} and visualized in Figure \ref{fig:timers} (i), {\blue{(iii)}}. {\magenta{We clarify that we use the same number of modes for both variables $u$ and $w$. Moreover, the relative error increases over time due to additive --reduced basis approximation-- error in every time point. Nevertheless, and even if the error increases, the max error for the concentration --that we are interested in-- is stabilizing after approximately sixty time steps evolution in an error of the order of $10^{-2}$.  In}} Table \ref{table:timers}, and in Figure \ref{fig:timers} (ii) we report the execution times for several numbers of basis {\magenta{functions}} including the projection of the full order matrices, the execution time of the online solver and the determination of the reduced problem and we compare them with the time execution at the full order level, namely $9.9623$ sec.

{\magenta{{We remark that the full order discretized differential operators that appear in equation (\ref{eq:system_linear0}) are parameter dependent and therefore, also at the reduced order level in order to compute the reduced differential operator, we need to assemble the full order operators. Possible ways to avoid such potentially expensive operation, relying on an affine approximation of the full order differential operator, could be to use hyper reduction techniques, \cite{Grepl2007,HeRoSta16}. In this work, since the attention is mainly devoted to the methodological development of a reduced-order method in an embedded boundary setting, rather than in its efficiency, we do not rely on such hyper reduction techniques and we assemble the full order differential operators also during the online stage. Considering that the most demanding computational effort is spent during the solution of the full order problem rather than in the assembly of the differential operators, as reported in Section \ref{exp:Geometrical parametrization},  it is anyway possible to achieve a computational speed up, and the related results have been reported in the present section.}
Finally, we underline that the cut elements FEM stabilization robustness is additionally verified using the conservation of mass test and in particular the relative error between the reduced and truth mass as it is illustrated in Figure \ref{Energy_Mass}. We can easily notice that the ROM solution and in particular the approximated mass conserves competently, as time passes and increasing the number of modes we use.
}}

 \begin{figure} \centering
    \includegraphics[width=0.67
    \textwidth]{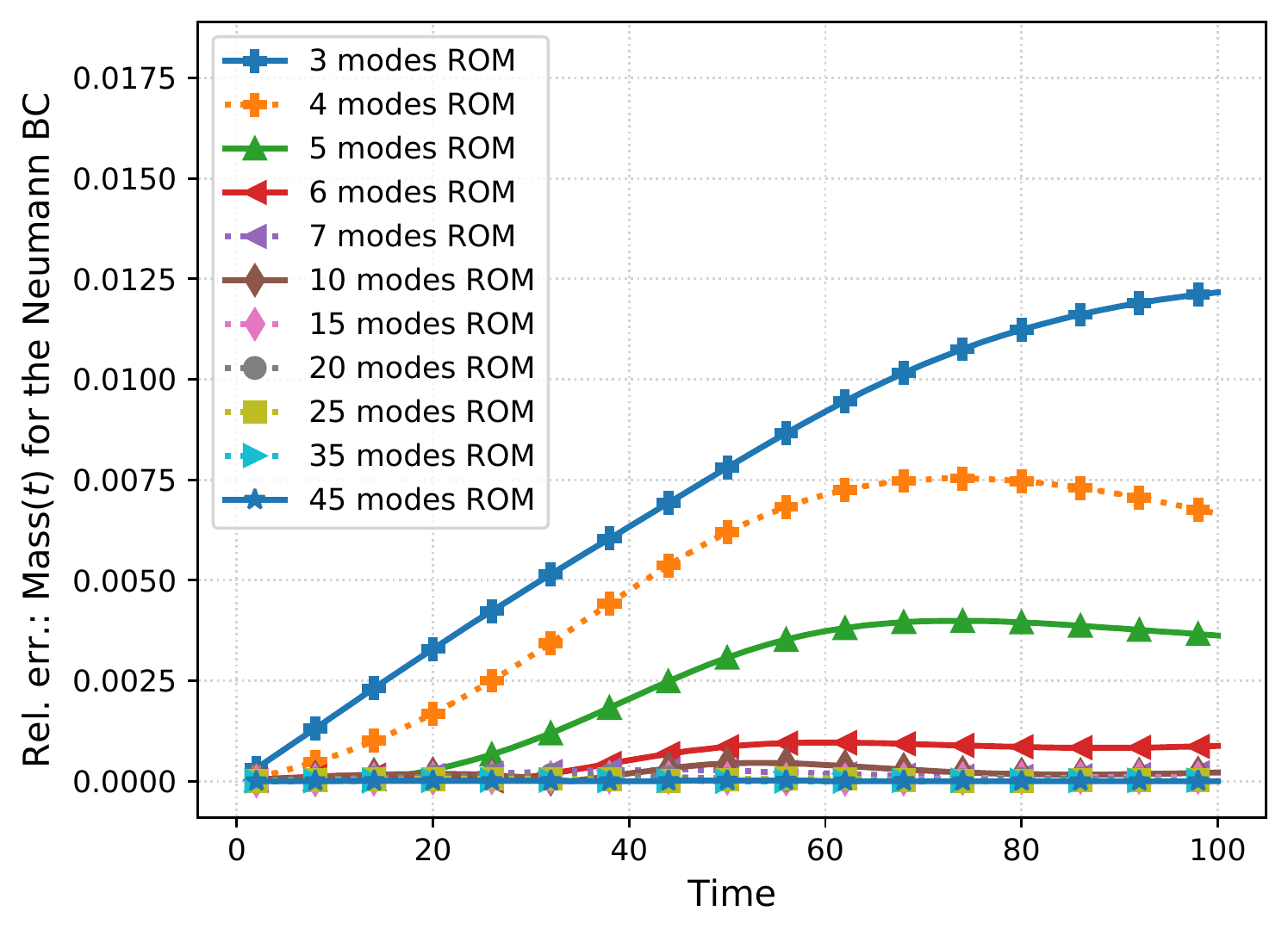}
    \\
    \includegraphics[width=0.67
    \textwidth]{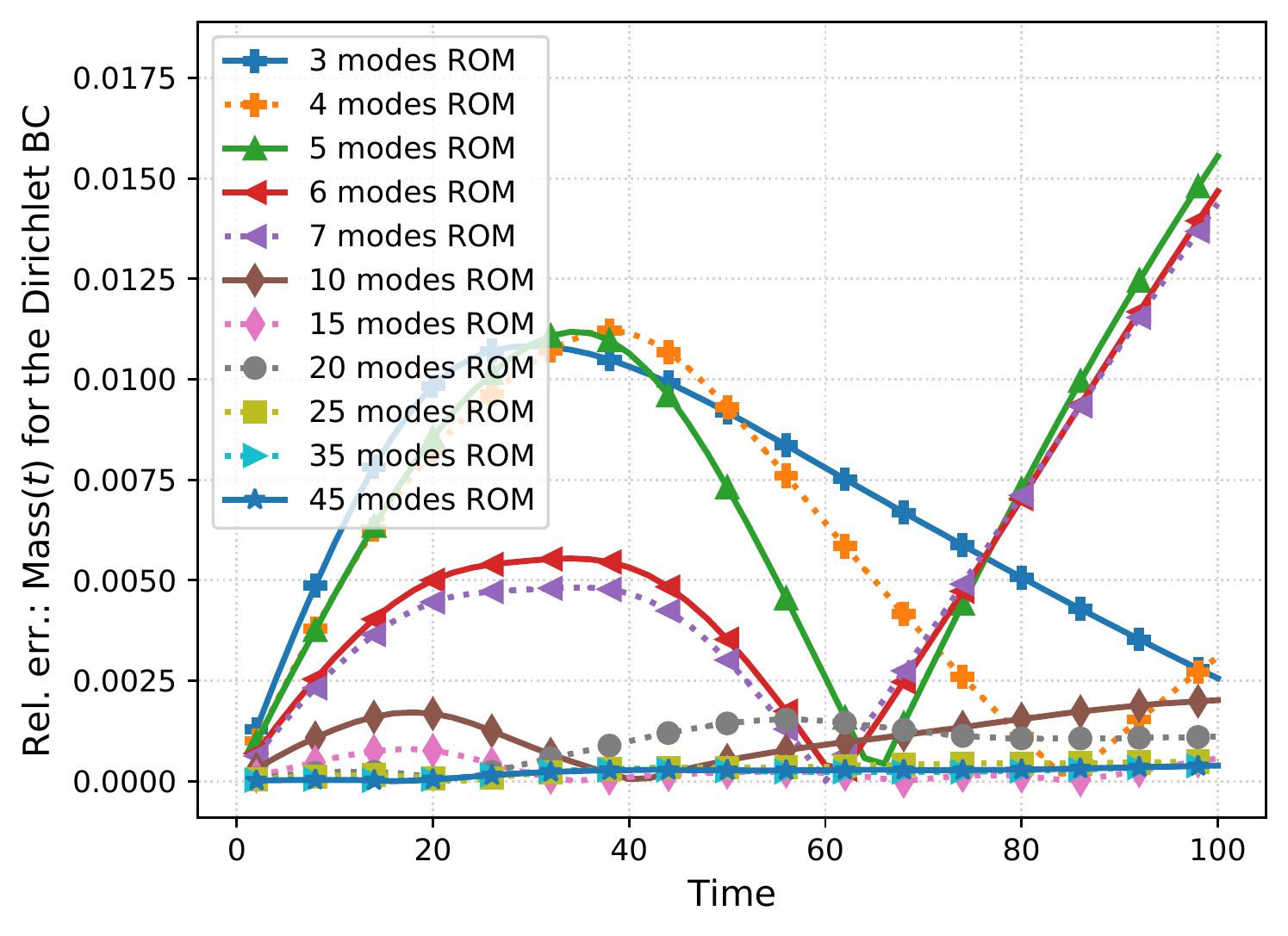}
\caption{
The conservation of mass for Neumann and Dirichlet embedded boundaries{\bblue{: evolution of the reduced basis approximation mass relative error  with respect to time, 
}} for various number of modes as considered in the numerical examples, for the parameter values {\magenta{$\mu=0.42917$ and $\mu=0.43774$}} respectively and
time instances $t=n\tau _n$, for $n = 1,...,100$. 
} 
  \label{Energy_Mass}
\end{figure}
%
\subsubsection{{\red{Challenges: larger scale geometric parametrization, time dependent geometry}}}\label{subsec:challenges}
{\red{In this paragraph, we examine our method numerically in cases of larger parameter range, as well as, evolutionary in time geometries. Two tests are provided, both emphasizing the more challenging Dirichlet embedded boundary condition case.
}}
\paragraph{{\red{Larger scale geometric parametrization.}}}
%
{\red{
We consider two experiments in which we show that if $\mu$ is in an interval with much larger or much smaller values, and in a larger range, our method reacts well.
In Table \ref{table:r_errors_extended_range} we report the relative errors for an extended train and test parameter set in two cases with (I) train set $[0.36,0.54]$, test set $[0.38,0.52]$ and (II) train set $[0.20, 0.48]$, test set $[0.22, 0.46]$. The latter results show that  even if we test the more difficult case of Dirichlet embedded boundary, the results are comparable with an order of $10^{-2}$ with the ones with much smaller train and test parameter intervals as in Table \ref{table:r_errors}.
\begin{table}\centering
  \begin{tabular}{|c||c||c|c|c|}
   \hline
    {\red{Snapshots for:}}& \multicolumn{4}{c|}{{\red{900 train parameters}}    } 
   \\ 
    \hline \hline
    {\red{Modes}} & \multicolumn{2}{c|}{{\red{Relative error (Dirichlet) (I) test}}}& \multicolumn{2}{c|}{{\red{Relative error (Dirichlet) (II) test}}} \\ \hline
    ({\red{$N_{POD}$}})       & {\red{Concentration $u$}}  & {\red{Potential $w$}}& {\red{Concentration $u$}}  & {\red{Potential $w$}}\\
    \hline
    {\red{1}}  & {\red{0.12816}}	& {\red{0.73905}} & {\red{0.16457}} & {\red{0.49708}} \\
    {\red{5}}  & {\red{0.12613}}	& {\red{0.73053}} & {\red{0.15532}} & {\red{0.48670}} \\
	{\red{10}} & {\red{0.15761}}	& {\red{0.31983}} & {\red{0.13032}} & {\red{0.31029}} \\
	{\red{15}} & {\red{0.10306}}	& {\red{0.23707}} & {\red{0.10653}} & {\red{0.24047}} \\
	{\red{20}} & {\red{0.08952}}	& {\red{0.21202}} & {\red{0.07884}} & {\red{0.19690}} \\
	{\red{25}} & {\red{0.07678}}	& {\red{0.18699}} & {\red{0.05915}} & {\red{0.16011}} \\
	{\red{30}} & {\red{0.05602}}	& {\red{0.14369}} & {\red{0.04719}} & {\red{0.12814}} \\
	{\red{35}} & {\red{0.04501}}	& {\red{0.12358}} & {\red{0.04175}} & {\red{0.11083}} \\
	{\red{40}} & {\red{0.03794}}	& {\red{0.10884}} & {\red{0.03454}} & {\red{0.09740}} \\
    {\red{45}} & {\red{0.03167}}	& {\red{0.09539}} & {\red{0.03254}} & {\red{0.09172}} \\
    \hline
  \end{tabular} 
      \caption{{\red{Geometrical parametrization and mean relative error between the full order solution and the reduced basis solution for the concentration and potential component tested onto 30 samples 
 chosen randomly inside the test parameter spaces and trained on 900 snapshots with train and test sets respectively: (I) train set $[0.36,0.54]$, test set $[0.38,0.52]$ and (II) train set $[0.20, 0.48]$, test set $[0.22, 0.46]$. 
      }}}
      \label{table:r_errors_extended_range}
\end{table}
%
Basic aspect so that we can achieve good relative error results is to keep the test geometry parameter $\mu$ inside an interval which is smaller of that of the training parameter. If we choose equal training and test parameter range sets the results are worse. Moreover, the smaller test  set we use, the fewer basis modes we need to use.  {\bblue{Further investigation has shown that there is a need of a relatively large number of $900$ snapshots, although they are calculated only once in the beginning of the ROM procedure -offline stage- to produce the reduced basis. We highlight at this point that the experiments as introduced in Section \ref{exp:Geometrical parametrization} using smaller parameter intervals can give similar results even if we use $600$ snapshots.}}}}
\paragraph{{\bblue{Time dependent geometry.}}}
%
{\bblue{
In this section we briefly show the results for an evolutionary in time embedded geometry. We tested the case in which the circle is moving periodically onto the $y$ axis with fixed diameter size $\delta=0.42$. So, the parameter $\mu(\cdot)$ now defines the motion which happens with respect to the center of the circle.
%
Our levelset function which is used for the embedded geometry description is
$
(x -\theta_1(t))^2 + (y  - \theta_2(t))^2 \leq \delta^2/4,
$
where $\theta(t) = (\theta_1(t), \theta_2(t))$ denotes the position of the center of the circle in the domain.
%
The motion of the circle is assumed to be known and in particular 
 we assume the periodic motion of the circle, i.e
$
\label{cylinder oscillation equation}
\theta(t) = \theta(0) + x_0\text{sin}(140\pi\mu(t))
 {\bm{j}},
$ 
where $x_0$ denotes the amplitude of the oscillation of the circle, $\mu(t)$ is a function of time and ${\bm{j}}$ is the unit vector in the vertical direction. Thus, the circle moves only vertically in our simulations. 
We choose
values for the constants 
$\theta(0)=(0, 0.1)$, 
$x_0=0.0039$,  
$\mu(t)=(\mu _\text{max} - \mu_\text{min})\frac{t}{T} + \mu_ \text{min}$, 
$\mu_\text{max}=0.15$, 
 $\mu_\text{min}=0.1$ and final time $T=100\tau_n$. 
%
All the rest of the data are as in Section \ref{exp:Geometrical parametrization} and the Dirichlet embedded geometry case. 
We choose $N_{\text{train}} = N_t = 100$ time instances, we run a POD on the set of the collected snapshots, and we obtain our basis functions with which we are going to compute the reduced solutions. 
%
For each $\theta^i$, $i=1, \dots, N_{\text{t}}$, we compute both the reduced solution 
and the corresponding full order solution. 
We compute the $L^2$ mean relative error 
for the concentration 
and for the potential 
 by taking the time average of the relative errors at each time $t_k$ and we obtain the mean approximation error as seen in Table \ref{table:r_errors_evol_geom}. 
%
%
%
%
%
%
\begin{table}
\centering
  \begin{tabular}{|c||c||c|}
   \hline
    {\bblue{Snapshots for:}}& \multicolumn{2}{c|}{{\bblue{100 time instances}}} \\ \hline \hline
    {\bblue{Modes}}     & \multicolumn{2}{c|}{{\bblue{Relative error for Dirichlet case}}} \\ \hline
    ({\bblue{{{$N_{POD}$}}}})       & {\bblue{Concentration $u$}}  & {\bblue{Potential $w$}} 
    \\
    \hline 
    {\bblue{1}} & {\bblue{0.15300}} & {\bblue{0.27402}} \\
    {\bblue{2}} & {\bblue{0.16976}} & {\bblue{0.29750}} \\
	{\bblue{3}} & {\bblue{0.12213}} & {\bblue{0.26491}} \\
	{\bblue{4}} & {\bblue{0.08475}} & {\bblue{0.16569}} \\
	{\bblue{5}} & {\bblue{0.06607}} & {\bblue{0.17079}} \\
	{\bblue{6}} & {\bblue{0.06032}} & {\bblue{0.16207}} \\
	{\bblue{7}} & {\bblue{0.03409}} & {\bblue{0.10952}} \\
	{\bblue{8}} & {\bblue{0.02403}}	& {\bblue{0.06251}} \\
	{\bblue{9}} & {\bblue{0.02093}}	& {\bblue{0.05422}} \\
    {\bblue{10}}& {\bblue{0.02040}}	& {\bblue{0.05319}} \\
    {\bblue{11}}& {\bblue{0.01506}} & {\bblue{0.03503}} \\
    \hline
  \end{tabular} 
      \caption{{\bblue{Geometrical parametrization and the average relative error between the full order solution and the reduced basis solution for the concentration and potential component --average for the errors at each time instance. Results are reported for different dimensions of the reduced basis spaces. The ROM has been trained onto 100 parameter samples coming from each time step.
      }}}
      \label{table:r_errors_evol_geom}
\end{table}
%
%
With this test case we can appreciate the advantage that we obtain by employing the CutFEM reduced order model that we proposed. Indeed, at every new time-step $t^i$ corresponds a different domain configuration, depending on the parameter $\theta^i$; without a CutFEM formulation, we would have to remesh at every time-step, making the offline phase of a reduced order model prohibitive from the computational cost point of view. Moreover, results show that we are able to obtain good results at the reduced order level even without employing a snapshot transportation during the offline phase. Although, the case of a larger motion scale and a more advanced development will be the topic of a future paper.
}}
\subsubsection{{\magenta{Outcome}}}
Numerical experiments consider {\magenta{Neumann}} boundary conditions as well as  Dirichlet ones. The tests clearly indicate that an efficient orthogonal decomposition projection-based reduced-order model can be derived over a  {\magenta{full order}} cut finite element method solver for the challenging (for both the full and the reduced level, see also the uncomfort  basis {\magenta{functions}} in Figure \ref{Fig:Modes}) nonlinear Cahn-Hilliard system. {\magenta{Efficient in the sense that the combination of embedded methods and reduced order models allowed us to obtain fast evaluation of parametrized problems, avoiding remeshing as well as the reference domain formulation, often used in the reduced order modeling for boundary fitted finite element formulations while}} 
 we rely only on an appropriate smooth enough fixed background grid. Last {\blue{Section's}} tests have clearly shown that {\magenta{sufficient good approximations}} can be obtained for this phase-field system at a reduced level.
\section{Concluding remarks and future developments}\label{sec:conclusions}
We conclude this work by noting that the above approach and the combination of unfitted mesh finite element methods with an embedded POD basis and IMEX type discretization, by using linear polynomials, imply good reduced basis approximation properties for {\blue{a}}  Cahn-Hilliard fourth-order diffusion nonlinear PDE system {\magenta{for which we have applied the splitting approach, leading to a coupled system of second order in space}}. Considering the reduced-order approximation, the background mesh approach appears beneficial even for nonsmooth pseudo-random initial data. 
As expected, and referring to previous related authors' works, \cite{KaratzasStabileNouveauScovazziRozza2018,KaratzasStabileAtallahScovazziRozza2018}, an increased error is noticed onto the Nitsche embedded boundary interface in the reduced level, {\magenta{which will be studied further in future work}}.  {\blue{We underline}} the significant execution time reduction considering the projection of the full order matrices, through the reduced execution time of the online solver and the solution of the reduced problem, {\blue{by}} capturing efficiently the full order solution information in a reduced level {\blue{solution.  As}} a perspective, we mention the {\blue{construction}} of higher-order IMEX {\blue{methods, more}} efficient methodologies for the affine decomposition of the discretization differential {\blue{operator,}} and the investigation of the applicability of well-known hyper reduction techniques, such as the empirical interpolation method, \cite{etna_vol32_pp145-161,KaKaTra21},  in the context of the aforementioned embedded reduced-order basis method.  Of future interest {\magenta{are}} also parabolic nonlinear partial differential equations {\blue{set}} in a more general framework and to test snapshots transportation techniques as presented in \cite{KaBaRO18}.
\section*{Acknowledgments}
This work is supported by the European Research Council Executive Agency by means of the H2020 ERC Consolidator Grant project AROMA-CFD ``Advanced Reduced Order Methods with  Applications in  Computational  Fluid  Dynamics'' - GA  681447, (PI: Prof. G. Rozza), {\blue{FARE-X-AROMA-CFD project by MIUR}}, INdAM-GNCS 2018 and {\blue{2019}} and by project FSE - European Social Fund - HEaD {\blue{``Higher}} Education and Development" SISSA operazione 1, Regione Autonoma Friuli - Venezia Giulia,  the Hellenic Foundation for Research and Innovation (HFRI) and  the  General  Secretariat  for  Research  and  Technology (GSRT),  
under  grant agreement No[1115], 
the ``First Call for H.F.R.I. Research Projects to support Faculty members and Researchers and the procurement of high-cost research equipment'' grant 3270,
and National Infrastructures for Research and Technology S.A. (GRNET S.A.) in the National HPC facility - ARIS - under project ID pa190902. {{The authors would like also to thank Dr Andrea Mola for useful instructions regarding the cutfem stabilization, and Dr Francesco Ballarin for fruitful discussions for the reduced order part.}}
\bibliographystyle{amsplain_gio}
\bibliography{bibfile_sissa}

\providecommand{\bysame}{\leavevmode\hbox to3em{\hrulefill}\thinspace}
\providecommand{\MR}{\relax\ifhmode\unskip\space\fi MR }
\providecommand{\MRhref}[2]{%
  \href{http://www.ams.org/mathscinet-getitem?mr=#1}{#2}
}
\providecommand{\href}[2]{#2}
\begin{thebibliography}{10}

\bibitem{ngsxfem}
\emph{{\magenta{{ngsxfem -- {A}dd-{O}n to {NG}Solve for unfitted finite element
  discretizations}}}}, {\magenta{https://github.com/ngsxfem/ngsxfem}}.

\bibitem{rbnics}
\emph{{\magenta{{RBniCS} - Reduced order modelling in {FEniCS}}}},
  {\magenta{https://www.rbnicsproject.org}}, {\magenta{2015}}.

\bibitem{ngsolve}
\emph{{\magenta{{NGSolve} - High performance multiphysics finite element
  software}}}, {\magenta{https://github.com/NGSolve/ngsolve}},
  {\magenta{2018}}.

\bibitem{Agosti17}
A.~Agosti, P.~F. Antonietti, P.~Ciarletta, M.~Grasselli, and M.~Verani,
  \emph{{A Cahn-Hilliard--type equation with application to tumor growth
  dynamics}}, Mathematical Methods in the Applied Sciences \textbf{40} (2017),
  no.~18, 7598--7626.

\bibitem{Alikakos2019}
N.~Alikakos, G.~Fusco, and P.~Smyrnelis, \emph{{Elliptic Systems of Phase
  Transition Type}}, 91, Monograph in the series Progress in Nonlinear
  Differential Equations and Their Applications, Birkhauser, 2018.

\bibitem{Alpak2016}
F.~O. Alpak, B.~Riviere, and F.~Frank, \emph{A phase-field method for the
  direct simulation of two-phase flows in pore-scale media using a
  non-equilibrium wetting boundary condition}, Computational Geosciences
  \textbf{20} (2016), no.~5, 881--908.

\bibitem{Anderson98}
D.~M. Anderson, G.~B. McFadden, and A.~A. Wheeler, \emph{Diffuse-interface
  methods in fluid mechanics}, Annual Review of Fluid Mechanics \textbf{30}
  (1998), no.~1, 139--165.

\bibitem{Antonopoulou2013EXISTENCEAR}
D.~C. Antonopoulou, D.~Farazakis, and G.~Karali, \emph{Malliavin calculus for
  the stochastic {C}ahn-{H}illiard/{A}llen-{C}ahn equation with unbounded noise
  diffusion}, {\magenta{Journal of Differential Equations}} \textbf{265}
  (2018), no.~7, 3168--3211.

\bibitem{ANTONOPOULOU20162383}
D.~C. Antonopoulou, G.~Karali, and A.~Millet, \emph{{{Existence and regularity
  of solution for a stochastic Cahn--Hilliard/Allen--Cahn equation with
  unbounded noise diffusion}}}, Journal of Differential Equations \textbf{260}
  (2016), no.~3, {\magenta{2383--2417}}.

\bibitem{BaFa2014}
M.~Balajewicz and C.~Farhat, \emph{Reduction of nonlinear embedded boundary
  models for problems with evolving interfaces}, Journal of Computational
  Physics \textbf{274} (2014), 489--504.

\bibitem{ballarin2015supremizer}
F.~Ballarin, A.~Manzoni, A.~Quarteroni, and G.~Rozza, \emph{{Supremizer
  stabilization of POD-Galerkin approximation of parametrized steady
  incompressible Navier--Stokes equations}}, International Journal for
  Numerical Methods in Engineering \textbf{102} (2015), no.~5, 1136--1161.

\bibitem{BeOhPaRoUr17}
P.~Benner, M.~Ohlberger, A.~Patera, G.~Rozza, and K.~Urban, \emph{{Model
  Reduction of Parametrized Systems}}, MS\&A series, vol.~17, Springer, 2017.

\bibitem{Bertozzi07}
A.~L. Bertozzi, S.~Esedoglu, and A.~Gillette, \emph{{Inpainting of Binary
  Images Using the Cahn--Hilliard Equation}}, Trans. Img. Proc. \textbf{16}
  (2007), no.~1, 285--291.

\bibitem{Blank2012}
L.~Blank, H.~Garcke, L.~Sarbu, T.~Srisupattarawanit, V.~Styles, and A.~Voigt,
  \emph{Phase-field approaches to structural topology optimization},
  pp.~245--256, Springer Basel, Basel, 2012.

\bibitem{BoPhD16}
J.~Bosch, \emph{{Fast Iterative Solvers for Cahn-Hilliard Problems}}, Ph.D.
  thesis, Otto-von-Guericke {\magenta{Universit\"{a}t}}, Magdeburg, 2016.

\bibitem{BURMAN20101217}
E.~Burman, \emph{Ghost penalty}, Comptes Rendus Mathematique \textbf{348}
  (2010), no.~21, {\magenta{1217--1220}}.

\bibitem{BuHa11}
E.~Burman and P.~Hansbo, \emph{{Fictitious domain finite element methods using
  cut elements: II. A stabilized Nitsche method}}, Applied Numerical
  Mathematics \textbf{52} (2011), no.~6, 2837--2862.

\bibitem{BuHa14}
\bysame, \emph{{Fictitious domain methods using cut elements: III. A stabilized
  Nitsche method for Stokes' problem}}, ESAIM: M2AN \textbf{48} (2014),
  no.~5-8, 859--874.

\bibitem{CaHi58}
J.~W. Cahn and J.~E. Hilliard, \emph{{Free Energy of a Nonuniform System. I.
  Interfacial Free Energy}}, The Journal of Chemical Physics \textbf{28}
  (1958), no.~2, 258--267.

\bibitem{Cha16}
F.~Chave, D.~Di~Pietro, F.~Marche, and F.~Pigeonneau, \emph{{A Hybrid
  High-Order Method for the Cahn--Hilliard problem in Mixed Form}}, SIAM
  Journal on Numerical Analysis \textbf{54} (2016), no.~3, 1873--1898.

\bibitem{Cherfils17}
L.~Cherfils, H.~Fakih, and A.~Miranville, \emph{{A Complex Version of the
  Cahn--Hilliard Equation for Grayscale Image Inpainting}}, Multiscale Modeling
  \& Simulation \textbf{15} (2017), no.~1, 575--605.

\bibitem{ChinestaEnc2017}
F.~Chinesta, A.~Huerta, G.~Rozza, and K.~Willcox, ch.~Model Reduction Methods,
  Encyclopedia of Computational Mechanics, Second Edition, pp.~1--36, John
  Wiley \& Sons, 2017.

\bibitem{Chinesta2011}
F.~Chinesta, P.~Ladeveze, and E.~Cueto, \emph{{A Short Review on Model Order
  Reduction Based on Proper Generalized Decomposition}}, Archives of
  Computational Methods in Engineering \textbf{18} (2011), no.~4, 395.

\bibitem{Choksi09}
R.~Choksi, M.~Peletier, and J.~Williams, \emph{{On the Phase Diagram for
  Microphase Separation of Diblock Copolymers: An Approach via a Nonlocal
  Cahn--Hilliard Functional}}, SIAM Journal on Applied Mathematics \textbf{69}
  (2009), no.~6, 1712--1738.

\bibitem{doi:10.1137/130943108}
K.~Chrysafinos and E.~N. Karatzas, \emph{{Error Estimates for Discontinuous
  Galerkin Time-Stepping Schemes for Robin Boundary Control Problems
  Constrained to Parabolic PDEs}}, SIAM Journal on Numerical Analysis
  \textbf{52} (2014), no.~6, 2837--2862.

\bibitem{ChKa2015}
K.~Chrysafinos and E.~N. Karatzas, \emph{{Symmetric error estimates for
  discontinuous Galerkin time-stepping schemes for optimal control problems
  constrained to evolutionary Stokes equations}}, Computational Optimization
  and Applications \textbf{60} (2015), no.~3, 719--751.

\bibitem{CLAUS2019185}
S.~Claus and P.~Kerfriden, \emph{A cutfem method for two-phase flow problems},
  Computer Methods in Applied Mechanics and Engineering \textbf{348} (2019),
  {\magenta{185--206}}.

\bibitem{CoFaGiSpre15}
P.~Colli, M.~Farshbaf-Shaker, G.~Gilardi, and J.~Sprekels, \emph{{Optimal
  Boundary Control of a Viscous Cahn--Hilliard System with Dynamic Boundary
  Condition and Double Obstacle Potentials}}, {SIAM Journal on Control and
  Optimization} \textbf{53} (2015), no.~4, 2696--2721.

\bibitem{DeGroot62}
S.~De~Groot and P.~Mazur, \emph{{Non-equilibrium thermodynamics}}, {(1962),
  Dover edition}, 2013.

\bibitem{Dumon20111387}
A.~Dumon, C.~Allery, and A.~Ammar, \emph{{Proper general decomposition (PGD)
  for the resolution of Navier--Stokes equations}}, Journal of Computational
  Physics \textbf{230} (2011), no.~4, 1387--1407.

\bibitem{ElioLar92}
C.~Elliott and S.~Larsson, \emph{{{Error estimates with smooth and nonsmooth
  data for a finite element method for the Cahn-Hilliard equation}}}, Math.
  Comp. \textbf{58} (1992), no.~S33-S36, 603--630.

\bibitem{Elliottt1989}
C.~M. Elliott, \emph{{{The Cahn-Hilliard Model for the Kinetics of Phase
  Separation}}}, pp.~35--73, Birkh{\"a}user Basel, Basel, 1989.

\bibitem{Elliott1989}
C.~M. Elliott, D.~A. French, and F.~A. Milner, \emph{{{A second order splitting
  method for the Cahn-Hilliard equation}}}, Numerische Mathematik \textbf{54}
  (1989), no.~5, 575--590.

\bibitem{Elliott1986}
C.~M. Elliott and Z.~Songmu, \emph{{{On the Cahn-Hilliard equation}}}, Archive
  for Rational Mechanics and Analysis \textbf{96} (1986), no.~4, 339--357.

\bibitem{Nitsche_4thorder_EmDoHa2010}
A.~Embar, J.~Dolbow, and I.~Harari, \emph{Imposing {{\magenta{{D}irichlet}}}
  boundary conditions with {Nitsche}'s method and spline-based finite
  elements}, International Journal for Numerical Methods in Engineering
  \textbf{83} (2010), no.~7, 877--898.

\bibitem{eyre_1998}
D.~J. Eyre, \emph{{\magenta{Unconditionally Gradient Stable Time Marching the
  Cahn-Hilliard Equation}}}, {\magenta{MRS Proceedings}}
  \textbf{{\magenta{529}}} ({\magenta{1998}}), {\magenta{39}}.

\bibitem{FuKovacsLarssonLindgren18}
D.~Furihata, M.~Kov{\`{a}}cs, S.~Larsson, and F.~Lindgren, \emph{{{Strong
  Convergence of a Fully Discrete Finite Element Approximation of the
  Stochastic Cahn--Hilliard Equation}}}, {{SIAM Journal on Numerical Analysis}}
  \textbf{56} (2018), no.~2, 708--731.

\bibitem{CrossGoudenege2012}
L.~Gouden{\`{e}}ge, D.~Martin, and G.~Vial, \emph{{{High Order Finite Element
  Calculations for the Cahn-Hilliard Equation}}}, Journal of Scientific
  Computing \textbf{52} (2012), no.~2, 294--321.

\bibitem{convectiveOptimal_Hinze17}
C.~Gr{\"a}{\ss}le, M.~Hinze, and N.~Scharmacher, \emph{{{POD for Optimal
  Control of the Cahn-Hilliard System Using Spatially Adapted Snapshots}}},
  Numerical Mathematics and Advanced Applications ENUMATH 2017 (F.~A. Radu,
  K.~Kumar, I.~Berre, J.~M. Nordbotten, and I.~S. Pop, eds.), Springer
  International Publishing, 2019, pp.~703--711.

\bibitem{Grepl2007}
M.~Grepl, Y.~Maday, N.~Nguyen, and A.~Patera, \emph{{Efficient reduced-basis
  treatment of nonaffine and nonlinear partial differential equations}}, ESAIM:
  M2AN \textbf{41} (2007), no.~3, 575--605.

\bibitem{grepl2005}
M.~Grepl and A.~Patera, \emph{{A posteriori error bounds for reduced-basis
  approximations of parametrized parabolic partial differential equations}},
  ESAIM: M2AN \textbf{39} (2005), no.~1, 157--181.

\bibitem{GURTIN96}
M.~E. Gurtin, D.~Polignone, and J.~Vinals, \emph{Two-phase binary fluids and
  immiscible fluids described by an order parameter}, Mathematical Models and
  Methods in Applied Sciences \textbf{06} (1996), no.~06, 815--831.

\bibitem{Haasdonk2008}
B.~Haasdonk and M.~Ohlberger, \emph{{Reduced basis method for finite volume
  approximations of parametrized linear evolution equations}}, Mathematical
  Modelling and Numerical Analysis \textbf{42} (2008), no.~2, 277--302.

\bibitem{etna_vol32_pp145-161}
B.~Haasdonk, M.~Ohlberger, and G.~Rozza, \emph{A reduced basis method for
  evolution schemes with parameter-dependent explicit operators}, Electron.
  Trans. Numer. Anal. \textbf{32} (2008), 145--161.

\bibitem{Harari_Nitsche15}
I.~Harari and E.~Grosu, \emph{A unified approach for embedded boundary
  conditions for fourth-order elliptic problems}, International Journal for
  Numerical Methods in Engineering \textbf{104} (2015), no.~7, 655--675.

\bibitem{HeRoSta16}
J.~Hesthaven, G.~Rozza, and B.~Stamm, \emph{{Certified Reduced Basis Methods
  for Parametrized Partial Differential Equations}}, SpringerBriefs in
  Mathematics, Springer International Publishing, 2016.

\bibitem{HintermullerHinze2013Ellipse}
M.~Hinterm\"{u}ller, M.~Hinze, and C.~Kahle, \emph{{An Adaptive Finite Element
  Moreau-Yosida-based Solver for a Coupled Cahn-Hilliard/Navier-Stokes
  System}}, J. Comput. Phys. \textbf{235} (2013), no.~C, 810--827.

\bibitem{HiKeWe17}
M.~Hinterm{\"{u}}ller, T.~Keil, and D.~Wegner, \emph{{Optimal Control of a
  Semidiscrete Cahn--Hilliard--Navier--Stokes System with Nonmatched Fluid
  Densities}}, {SIAM Journal on Control and Optimization} \textbf{55} (2017),
  no.~3, 1954--1989.

\bibitem{controldistrCHiWe12}
M.~{\magenta{Hinterm\"{u}ller}} and D.~Wegner, \emph{{Distributed Optimal
  Control of the Cahn--Hilliard System Including the Case of a Double-Obstacle
  Homogeneous Free Energy Density}}, SIAM Journal on Control and Optimization
  \textbf{50} (2012), no.~1, 388--418.

\bibitem{HinzeNS13}
M.~Hinze and C.~Kahle, \emph{{{A Nonlinear Model Predictive Concept for Control
  of Two-Phase Flows Governed by the Cahn-Hilliard Navier-Stokes System}}},
  System Modeling and Optimization (Berlin, Heidelberg) (D.~H{\"o}mberg and
  F.~Tr{\"o}ltzsch, eds.), Springer Berlin Heidelberg, 2013, pp.~348--357.

\bibitem{Israelachvili11}
J.~N. Israelachvili, \emph{{Intermolecular and Surface Forces}}, Elsevier,
  2011.

\bibitem{Jeong2015}
D.~Jeong and J.~Kim, \emph{{{Microphase separation patterns in diblock
  copolymers on curved surfaces using a nonlocal Cahn-Hilliard equation}}}, The
  European Physical Journal E \textbf{38} (2015), no.~11, 117.

\bibitem{Junseok16}
K.~Junseok, L.~Seunggyu, C.~Yongho, L.~Seok-Min, and J.~Darae, \emph{{Basic
  Principles and Practical Applications of the Cahn--Hilliard Equation}},
  Mathematical Problems in Engineering (2016), no.~1, 79--141.

\bibitem{Kalashnikova_ROMcomprohtua}
I.~Kalashnikova and M.~F. Barone, \emph{{On the stability and convergence of a
  Galerkin reduced order model (ROM) of compressible flow with solid wall and
  far-field boundary treatment}}, International Journal for Numerical Methods
  in Engineering \textbf{83} (2010), no.~10, 1345--1375.

\bibitem{KaraliNagase14}
G.~Karali and Y.~Nagase, \emph{{{On the existence of solution for a
  Cahn--Hilliard/Allen--Cahn equation}}}, Discrete and Continuous Dynamical
  Systems - S \textbf{7} (2014), no.~1937-1632\_2014\_1\_127, 127.

\bibitem{KaBaRO18}
E.~N. Karatzas, F.~Ballarin, and G.~Rozza, \emph{Projection-based reduced order
  models for a cut finite element method in parametrized domains}, Computers \&
  Mathematics with Applications \textbf{79} (2020), no.~3,
  {\magenta{833--851}}.

\bibitem{KaNoBaRo20}
E.~N. Karatzas, M.~Nonino, F.~Ballarin, and G.~Rozza, \emph{{A Reduced order
  cut finite element basis for stationary and evolutionary geometrically
  parameterized Navier--Stokes systems}}, {\magenta{Accepted for publication in
  Computers \& Mathemetics with Applications, preprint at arXiv:2010.04953}},
  2021.

\bibitem{KaratzasStabileAtallahScovazziRozza2018}
E.~N. Karatzas, G.~Stabile, N.~Atallah, G.~Scovazzi, and G.~Rozza, \emph{{A
  Reduced Order Approach for the Embedded Shifted Boundary {FEM} and a Heat
  Exchange System on Parametrized Geometries}}, In: Fehr J., Haasdonk B. (eds)
  IUTAM Symposium on Model Order Reduction of Coupled Systems, Stuttgart,
  Germany, May 22--25, 2018. IUTAM Bookseries, vol 36. Springer, Cham (2020).

\bibitem{KaratzasStabileNouveauScovazziRozza2018}
E.~N. Karatzas, G.~Stabile, L.~Nouveau, G.~Scovazzi, and G.~Rozza, \emph{{{A
  reduced basis approach for PDEs on parametrized geometries based on the
  shifted boundary finite element method and application to a Stokes flow}}},
  Computer Methods in Applied Mechanics and Engineering \textbf{347} (2019),
  {\magenta{568--587}}.

\bibitem{KaratzasStabileNouveauScovazziRozzaNS2019}
E.~N. Karatzas, G.~Stabile, L.~Nouveau, G.~Scovazzi, and G.~Rozza, \emph{A
  reduced-order shifted boundary method for parametrized incompressible
  {N}avier--{S}tokes equations}, {\magenta{Computer Methods in Applied
  Mechanics and Engineering}} \textbf{{\magenta{370}}} (2020),
  {\magenta{113--273}}.

\bibitem{KaKaTra21}
G.~Katsouleas, E.~N. Karatzas, and F.~Travlopanos, \emph{{Discrete Empirical
  Interpolation and unfitted mesh FEMs: application in PDE-constrained
  optimization (2021)}}, {Submitted, arXiv:2010.09059}.

\bibitem{KaKaTra20}
\bysame, \emph{{{Cut finite element error estimates for a class of nonlinear
  elliptic PDEs}}}, Loughborough University, doi:
  {10.17028/rd.lboro.12154854.v1}, extended version at arXiv:2003.06489, 2020,
  pp.~1--6.

\bibitem{Kunisch2002492}
K.~Kunisch and S.~Volkwein, \emph{{Galerkin proper orthogonal decomposition
  methods for a general equation in fluid dynamics}}, SIAM Journal on Numerical
  Analysis \textbf{40} (2002), no.~2, 492--515.

\bibitem{Lehrenfeld2018L2errorAO}
C.~Lehrenfeld and A.~Reusken, \emph{{L2-error analysis of an isoparametric
  unfitted finite element method for elliptic interface problems}}, vol.~2,
  Journal of Numerical Mathematics, 2019, pp.~85--99.

\bibitem{LI20182100}
C.~Li, R.~Qin, J.~Ming, and Z.~Wang, \emph{{{A discontinuous Galerkin method
  for stochastic Cahn--Hilliard equations}}}, Computers \& Mathematics with
  Applications \textbf{75} (2018), no.~6, {\magenta{2100--2114}}, 2nd Annual
  Meeting of SIAM Central States Section, September 30-October 2, 2016.

\bibitem{LI2013102}
Y.~Li, D.~Jeong, J.~Shin, and J.~Kim, \emph{{{A conservative numerical method
  for the Cahn--Hilliard equation with Dirichlet boundary conditions in complex
  domains}}}, Computers \& Mathematics with Applications \textbf{65} (2013),
  no.~1, {\magenta{102--115}}.

\bibitem{pseudo2}
M.~Luby, \emph{{{\magenta{Pseudorandomness and Cryptographic Applications}}}},
  {\magenta{ISBN 9780691025469}}, {\magenta{Princeton University Press}},
  {\magenta{1996}}.

\bibitem{pseudo1}
M.~Matsumoto and T.~Nishimura, \emph{{\magenta{Mersenne twister: a
  623-dimensionally equi-distributed uniform pseudo-random number generator}}},
  {\magenta{ACM Transactions on Modeling and Computer Simulation}}
  \textbf{{\magenta{8 (1)}}} ({\magenta{1998}}), {\magenta{3--30}}.

\bibitem{NOVICKCOHEN1984277}
A.~Novick-Cohen and L.~A. Segel, \emph{{{Nonlinear aspects of the Cahn-Hilliard
  equation}}}, Physica D: Nonlinear Phenomena \textbf{10} (1984), no.~3,
  277--298.

\bibitem{quarteroniRB2016}
A.~Quarteroni, A.~Manzoni, and F.~Negri, \emph{Reduced basis methods for
  partial differential equations}, vol.~92, UNITEXT/La Matematica per il 3+2
  book series, Springer International Publishing, 2016.

\bibitem{REGAZZONI18}
F.~Regazzoni, N.~Parolini, and M.~Verani, \emph{Topology optimization of
  multiple anisotropic materials, with application to self-assembling diblock
  copolymers}, Computer Methods in Applied Mechanics and Engineering
  \textbf{338} (2018), 562--596.

\bibitem{Reshma14}
S.~Reshma and H.~Hansa, J~Thattil, \emph{{Inpainting of Binary Images Using the
  Cahn-Hilliard Equation}}, International Journal of Computer Science
  Engineering and Technology \textbf{4} (2014), no.~11, 296--300.

\bibitem{RoccaSprekels15conv_control}
E.~Rocca and J.~Sprekels, \emph{{Optimal Distributed Control of a Nonlocal
  Convective Cahn--Hilliard Equation by the Velocity in Three Dimensions}},
  SIAM Journal on Control and Optimization \textbf{53} (2015), no.~3,
  1654--1680.

\bibitem{RokhzadiPhD18}
A.~Rokhzadi, \emph{{IMEX and Semi-Implicit Runge-Kutta Schemes for CFD
  Simulations}}, Ph.D. thesis, {Civil Engineering Department, Faculty of
  Engineering, University of Ottawa}, 2018.

\bibitem{Rozza2004ReducedBM}
G.~Rozza, \emph{{Reduced {\magenta{basis methods for elliptic equations in
  subdomains with a-posteriori error bounds and adaptivity}}}}, App. Num. Math.
  \textbf{55} (2005), no.~4, 403--424.

\bibitem{Rozza2009}
G.~Rozza, \emph{{Reduced basis methods for Stokes equations in domains with
  non-affine parameter dependence}}, Computing and Visualization in Science
  \textbf{12} (2009), no.~1, 23--35.

\bibitem{Rozza2008}
G.~Rozza, D.~Huynh, and A.~Patera, \emph{Reduced basis approximation and a
  posteriori error estimation for affinely parametrized elliptic coercive
  partial differential equations: {A}pplication to transport and continuum
  mechanics}, Archives of Computational Methods in Engineering \textbf{15}
  (2008), no.~3, 229--275.

\bibitem{RoHuMa13}
G.~Rozza, D.~B.~P. Huynh, and A.~Manzoni, \emph{{Reduced basis approximation
  and a posteriori error estimation for Stokes flows in parametrized
  geometries: Roles of the inf-sup stability constants}}, Numerische Mathematik
  \textbf{125} (2013), no.~1, 115--152.

\bibitem{RoVe07}
G.~Rozza and K.~Veroy, \emph{{On the stability of the reduced basis method for
  Stokes equations in parametrized domains}}, Computer Methods in Applied
  Mechanics and Engineering \textbf{196} (2007), no.~7, 1244--1260.

\bibitem{schoeberl}
J.~Sch\"oberl, A.~Arnold, J.~Erb, J.~M. Melenk, and T.~P. Wihler,
  \emph{{\magenta{C++11 implementation of finite elements in {NGS}olve}}},
  Tech. report, {\magenta{Institute for Analysis and Scientific Computing,
  Vienna University of Technology, ASC Report 30/2014}}, {\magenta{2014}}.

\bibitem{Schott2016}
B.~Schott, \emph{{Stabilized Cut Finite Element Methods for Complex Interface
  Coupled Flow Problems}}, Ph.D. thesis, {Technische Universit\"{a}t
  M\"{u}nchen (TUM)}, 2016.

\bibitem{Shenyang04}
H.~Shenyang, \emph{{Phase-field Models of Microstructure Evolution in a System
  with Elastic Inhomogeneity and Defects}}, Ph.D. thesis, {Pennsylvania State
  University, Department of Materials Science and Engineering}, 2004.

\bibitem{Veroy2003}
K.~Veroy, C.~Prud'homme, and A.~Patera, \emph{{Reduced-basis approximation of
  the viscous Burgers equation: rigorous a posteriori error bounds}}, Comptes
  Rendus Mathematique \textbf{337} (2003), no.~9, 619--624.

\bibitem{WELLS2006860}
G.~N. Wells, E.~Kuhl, and K.~Garikipati, \emph{{{A discontinuous Galerkin
  method for the Cahn--Hilliard equation}}}, Journal of Computational Physics
  \textbf{218} (2006), no.~2, 860--877.

\bibitem{Signori2019}
G.~Welper, \emph{{Optimal treatment for a phase field system of Cahn-Hilliard
  type modeling tumor growth by asymptotic scheme}}, arXiv:1902.01079v2, 2019.

\bibitem{WODO20116037}
O.~Wodo and B.~Ganapathysubramanian, \emph{{{Computationally efficient solution
  to the Cahn--Hilliard equation: Adaptive implicit time schemes, mesh
  sensitivity analysis and the 3D isoperimetric problem}}}, Journal of
  Computational Physics \textbf{230} (2011), no.~15, 6037--6060.

\bibitem{CrossXU2019524}
M.~Xu, H.~Guo, and Q.~Zou, \emph{{{Hessian recovery based finite element
  methods for the Cahn-Hilliard equation}}}, {{Journal of Computational
  Physics}} \textbf{386} (2019), 524--540.

\bibitem{Zhang2018}
X.~Zhang, H.~Li, and C.~Liu, \emph{{{Optimal Control Problem for the
  Cahn--Hilliard/Allen--Cahn Equation with State Constraint}}}, Applied
  Mathematics {\&} Optimization (2018).

\bibitem{conCHZhao2014}
X.~Zhao and C.~Liu, \emph{{{Optimal Control for the Convective Cahn--Hilliard
  Equation in 2D Case}}}, Applied Mathematics {\&} Optimization \textbf{70}
  (2014), no.~1, 61--82.

\bibitem{ZHAO2017177}
Y.~Zhao, D.~Schillinger, and B.-X. Xu, \emph{{{Variational boundary conditions
  based on the Nitsche method for fitted and unfitted isogeometric
  discretizations of the mechanically coupled Cahn-Hilliard equation}}},
  Journal of Computational Physics \textbf{340} (2017), 177--199.

\bibitem{Zhou2006}
S.~Zhou and M.~Y. Wang, \emph{{{Multimaterial structural topology optimization
  with a generalized Cahn--Hilliard model of multiphase transition}}},
  Structural and Multidisciplinary Optimization \textbf{33} (2006), no.~2, 89.

\end{thebibliography}
\end{document}